\documentclass[11pt]{article}

\usepackage[T1]{fontenc}
\usepackage{lmodern}
\usepackage{amsmath,amssymb,amsthm,mathtools}
\usepackage[a4paper,margin=30mm]{geometry}
\usepackage{microtype}
\usepackage{enumitem}
\usepackage{pgfplots}
\pgfplotsset{compat=1.18}
\usepackage[colorlinks=true,linkcolor=blue,citecolor=blue,urlcolor=blue]{hyperref}

\newtheorem{theorem}{Theorem}[section]
\newtheorem{proposition}[theorem]{Proposition}
\newtheorem{corollary}[theorem]{Corollary}
\newtheorem{lemma}[theorem]{Lemma}
\theoremstyle{remark}
\newtheorem{remark}[theorem]{Remark}

\newcommand{\Will}{\mathcal W}
\newcommand{\R}{\mathbb R}
\newcommand{\Sph}{\mathbb S}
\newcommand{\dd}{\,\mathrm d}
\newcommand{\Area}{\operatorname{Area}}
\newcommand{\spt}{\operatorname{spt}}
\newcommand{\BL}{\mathrm{BL}}
\newcommand{\Var}{\mathbf V}
\newcommand{\TV}{\mathrm{TV}}

\newcommand{\miwfZ}{\mathbb Z}
\newcommand{\miwfE}{\mathcal E}
\newcommand{\miwfD}{\mathcal D}
\newcommand{\miwfnorm}[1]{\left\lVert#1\right\rVert}

\newcommand{\HH}{\mathbb H}
\newcommand{\TT}{\mathrm{TT}}
\newcommand{\rank}{\operatorname{rank}}
\newcommand{\Diff}{\operatorname{Diff}}
\newcommand{\Mob}{\operatorname{Mob}}

\numberwithin{equation}{section}

\title{Geometric Endpoints, Global Trajectories, and Teichm\"uller Projection
for the M\"obius-Invariant Willmore Flow}
\author{
Mohameden Ahmedou\thanks{Mathematisches Institut, Justus-Liebig-Universit\"at
Gie\ss en, Arndtstrasse 2, 35392 Gie\ss en, Germany.
E-mail: \texttt{Mohameden.Ahmedou@math.uni-giessen.de}.}
\and
Ruben Jakob\thanks{Mathematics Department, Technion--Israel Institute of
Technology, 3200003 Haifa, Israel.
E-mail: \texttt{rubenj@technion.ac.il}.}}
\date{November 8, 2026}

\hypersetup{
 pdfauthor={Mohameden Ahmedou and Ruben Jakob},
 pdftitle={Geometric Endpoints, Global Trajectories, and Teichm\"uller Projection
for the M\"obius-Invariant Willmore Flow}
}

\begin{document}
\maketitle
\begin{abstract}
We study endpoints and conformal projection of the M\"obius-invariant
Willmore flow in the round three-sphere. Under uniform nonumbilicity or
an integrable weighted condition, image-measure transport yields a
unique finite-time integral-varifold limit. Reduced Hopf trajectories
have unique anchored finite-time profile and marked conformal endpoints
at arbitrary energy. Every smooth global reduced trajectory converges
smoothly modulo reparametrization to an elastic profile, with quantitative
rates and finite dissipation length. We compute the differential of the
conformal-class projection along reduced trajectories and prove that
its maximal rank persists at every finite time when it holds initially. Near
the Clifford torus, Hopf tori with the same marked conformal type can
have noncollinear projected velocities, ruling out a common autonomous
Teichm\"uller equation even after positive time recalibrations.
\end{abstract}


\section{Introduction and main results}

\subsection{The Willmore functional and the classical flow}

A recurring issue in the analysis of geometric flows is the gap
between subsequential compactness and convergence of the complete
trajectory.  Energy and curvature bounds may provide geometric cluster
points as time approaches an endpoint, but they do not by themselves
exclude different time sequences from selecting different limits.  In
parametric surface flows this problem is compounded by the freedom of
reparametrization and by the possible noncompactness of the ambient
symmetry group.  We study this endpoint question together with the conformal projection
of the M\"obius-invariant Willmore flow. The geometric estimates give
convergence, while explicit Hopf examples determine what information
the conformal class alone cannot retain.

For a closed immersed surface $F\colon\Sigma\to\Sph^3$, where
$\Sph^3\subset\R^4$ is the unit round three-sphere with metric
$g_{\Sph^3}$, the Willmore functional in the normalization used below is
\[
 \Will(F)=\int_\Sigma\left(1+\frac14|H|^2\right)\dd\mu.
\]
Here $H$ is the trace mean-curvature vector in $\Sph^3$ and
$\mu=\mu_{F^*g_{\Sph^3}}$ is the induced area measure.
Equivalently, after stereographic projection, this is the Euclidean
Willmore energy.  Its invariance under conformal transformations of the
ambient sphere is the source of both its geometric rigidity and its
noncompactness.  The Li--Yau inequality \cite{LiYau} shows, in
particular, that the threshold $8\pi$ excludes self-intersections for
closed immersions of smaller energy.  For embedded closed surfaces of positive genus, the theorem of
Marques and Neves \cite[Theorem~A]{MarquesNeves} gives the lower bound
$2\pi^2$, with equality precisely for conformal images of the Clifford
torus. The same bound and equality statement hold for immersions:
a nonembedded immersion has multiplicity at least two at some point,
so Li--Yau gives $\Will\geq8\pi>2\pi^2$.

The negative $L^2$-gradient flow of $\Will$ is the classical Willmore
flow.  Early analytic treatments include the small-energy theory of
Kuwert--Sch\"atzle \cite{KuwertSchaetzleSmall} and the existence and
stability analysis near spheres of Simonett \cite{Simonett}.
Kuwert and Sch\"atzle developed the general fourth-order evolution and
its concentration analysis in
\cite{KuwertSchaetzleGradient}; Kuwert--Scheuer \cite{KuwertScheuer}
later obtained quantitative stability estimates for the barycenter,
quadratic moment, and, in codimension one, enclosed volume and averaged
mean curvature in the small-energy regime.
Embeddedness need not be preserved \cite{MayerSimonett}, and examples of singular behavior are
known \cite{Blatt}.  Dall'Acqua, M\"uller, Sch\"atzle, and Spener
\cite{DallAcquaMuellerSchaetzleSpener} obtained a detailed analysis for
tori of revolution, including degenerating examples.
Palmurella and Rivi\`ere
\cite{PalmurellaRiviereApproach,PalmurellaRiviereFlow} developed a
parametric framework for the classical flow and proved short-time
existence estimates in terms of geometric data.
In this framework, energy concentration, possible branch points, and
degeneration of the conformal class are central endpoint issues; see
the subsection on Willmore flow in \cite{LanMartinoRiviere} for their discussion and for
a broader account of the classical flow and Willmore theory.

Compactness for sequences with bounded Willmore energy is also built
on the weak theory of conformal immersions.  Kuwert--Li \cite{KuwertLi}
developed compactness for $W^{2,2}$ conformal immersions.
Rivi\`ere \cite{RiviereAnalysis} established a weak formulation and
regularity theory for the Willmore equation, together with compactness
results for Willmore surfaces.  The energy bounds for degenerating
conformal classes in \cite{KuwertLi,RiviereDegenerating} and the
conformal-factor estimates in \cite{SchaetzleConformalFactor} clarify
the role of parametrization and conformal structure.
Energy quantization for sequences of Willmore surfaces does not, by
itself, apply to arbitrary time slices of a flow.
At the level of geometric images, Allard's compactness
theorem \cite{Allard} provides the natural integral-varifold framework.

\subsection{The M\"obius-invariant Willmore flow}

The M\"obius-invariant Willmore flow (MIWF) is
\[
 \partial_tF=-|A^0|^{-4}\nabla_{L^2}\Will(F),
\]
where $A^0$ denotes the trace-free second fundamental form.  The
weight $|A^0|^{-4}$ corrects the conformal scaling of the classical
gradient and makes the evolution covariant under fixed M\"obius
transformations without a rescaling of time.  The price is a singular
coefficient at umbilic points.  An umbilic-free closed oriented surface in $\Sph^3$ has a smooth
principal line field. On the double cover orienting this line field it
gives a nowhere-zero tangent vector field, so the Euler characteristic
of the cover, and hence of the surface, is zero. Thus the connected
closed domain in this setting is a torus.

The analytic development of the MIWF began with Jakob's short-time
existence and uniqueness theorem \cite{JakobShortTime}; functional
analytic properties and regularity were established in
\cite{JakobRegularity}.  The Hopf fibration relates distinguished
toroidal trajectories to fourth-order curve flows on $\Sph^2$.  The
corresponding curve geometry and reduction mechanisms were developed
in Jakob's work on Willmore flow of Hopf tori \cite{JakobHopf}, building
on Pinkall's description of Hopf tori as inverse images of closed
curves under the Hopf map \cite{Pinkall}.  For the related classical
Hopf-torus endpoint problem, see also \cite{Ahmedou.JakobHopf}.
We use only the equivariance
and geometric reduction formulas that are unaffected by the later
corrections in \cite{JakobHopfCorrection}; in particular, we do not use
the corrected assertions in parts \textup{(II)}--\textup{(III)} of
Theorem~1 or in the second part of Proposition~6 of \cite{JakobHopf}.
Jakob proved global existence, full convergence, and stability results
for the MIWF in \cite{JakobGlobal}, including stability of smoothly
convergent trajectories whose limits are umbilic-free local minimizers
of the Willmore functional.

The endpoint compactness theory is markedly less rigid.  Jakob's
subsequent theorem \cite[Theorem~1.1]{Jakob} treats toroidal MIWF
trajectories in $\Sph^3$ with initial energy at most $8\pi$.  From
every sequence of times approaching the maximal time it extracts a
subsequence whose multiplicity-one image measures converge weakly to
the weight of an integral rectifiable varifold.  The limit is zero, or
its support is a closed embedded orientable Lipschitz surface of genus
zero or one; in the genus-one case the support admits a uniformly
conformal bi-Lipschitz parametrization in
$W^{2,2}\cap W^{1,\infty}$.  For trajectories in the Hopf sector every
nonzero subsequential support is an embedded $C^1$ Hopf torus
\cite[Theorem~1.2]{Jakob}.  These results identify all possible
subsequential endpoints, but they do not show that two sequences of
times along the same trajectory have the same limit.
In our earlier preprint \cite{AhmedouJakobUnique}, we proved
finite-time varifold uniqueness under uniform nonumbilicity, together
with the low-energy anchored Hopf-profile endpoint, convergence of
its Pinkall modulus \cite[Corollary~6.6]{AhmedouJakobUnique}, and an
infinite-time criterion under finite dissipation length. We retain
these results and their proofs here and develop the all-energy Hopf
endpoint, global reduced trajectories, and conformal projection
questions in the sections that follow.

\subsection{Finite-time geometric endpoints}

The approach of \cite{AhmedouJakobUnique} estimates the evolution of the
geometric image measure itself. The resulting estimate is independent of the choice
of coordinates on the torus and does not require a time-dependent
M\"obius normalization.  Under a uniform quantitative nonumbilicity
hypothesis it makes the image measures Cauchy at every finite endpoint;
integral compactness then upgrades the measure limit to a unique
varifold limit.  Our principal result is the following.

\begin{theorem}[Unique finite-time limit]
\label{thm:main}
Let $F\colon[0,T)\times\Sigma\to\Sph^3$ be a smooth MIWF of a compact
torus, where $0<T<\infty$, and write $F_t:=F(t,\cdot)$.  Let $A_t^0$
denote the trace-free second fundamental form of $F_t$.  Put
\[
 W(t):=\Will(F_t),\qquad W_0:=W(0),\qquad
 \mu_t:=\mu_{F_t^*g_{\Sph^3}}.
\]
Assume that
for some $c_0>0$,
\[
 \inf_{0<t<T}\min_{x\in\Sigma}|A_t^0(x)|^2\geq c_0.
\]
Let $\Var_t$ denote the multiplicity-counting integral two-varifold in
$\R^4$ induced by $F_t$, and set
$\lambda_t:=\lVert\Var_t\rVert=(F_t)_\#\mu_t$.
Then there exists a unique integral two-varifold
$\Var_T$ in $\R^4$, possibly the zero varifold, such that, with
$\rightharpoonup$ denoting weak varifold convergence,
\begin{equation}
 \label{eq:main-convergence}
 \Var_t\rightharpoonup\Var_T
 \qquad\text{as }t\nearrow T.
\end{equation}
No reparametrizations of $\Sigma$ and no M\"obius transformations of
$\Sph^3$ are required.

More precisely, there are finite nonnegative Radon measures $\mu_T$ on
$\Sigma$ and $\lambda_T=\lVert\Var_T\rVert$ on $\Sph^3$ such that
$\mu_t\to\mu_T$ in total variation and
$\lambda_t\to\lambda_T$ in the bounded-Lipschitz distance
$d_{\BL}$ defined in \eqref{eq:bl-distance}.  Setting
\[
 W_T:=\lim_{t\nearrow T}W(t),
\]
there holds, for every $0<s<T$,
\begin{align}
 \label{eq:main-endpoint-tv}
 \lVert\mu_T-\mu_s\rVert_{\TV}
 &\leq \frac{2\sqrt{W(s)}}{c_0}
 \sqrt{(T-s)\bigl(W(s)-W_T\bigr)},\\
 \label{eq:main-endpoint-bl}
 d_{\BL}(\lambda_T,\lambda_s)
 &\leq \frac{2\sqrt{W(s)}}{c_0}
 \sqrt{(T-s)\bigl(W(s)-W_T\bigr)}.
\end{align}
Moreover, $\Var_T$ has generalized Euclidean mean-curvature vector
$H_{\Var_T}^{\R^4}\in L^2(\lVert\Var_T\rVert;\R^4)$ and
\begin{equation}
 \label{eq:main-curvature-lsc}
 \frac14\int_{\R^4}|H_{\Var_T}^{\R^4}|^2
 \dd\lVert\Var_T\rVert\leq W_T.
\end{equation}
\end{theorem}

When $T=T_{\max}(F_0)<\infty$, we call the resulting maximal
trajectory a \emph{singular flow line}, consistently with
\cite[Definition~2.1]{Jakob}.  The theorem itself is slightly more
general, since it applies on every finite interval on which the stated
uniform nonumbilicity bound holds.

The image-measure estimate follows both the motion of the immersion and
the variation of its induced area. Uniform mass and first-variation
bounds give integral cluster varifolds; their common weight determines
their almost-everywhere tangent planes and makes the cluster unique.
Proposition~\ref{prop:weighted-endpoint} replaces uniform nonumbilicity
by an integrable weighted transport condition. No low-energy hypothesis
is needed for either uniqueness theorem. For finite maximal trajectories
with initial energy at most $8\pi$, the compactness theorem of
\cite{Jakob} then makes the zero, genus-zero, and genus-one endpoint
alternatives sequence independent. The limiting parametrization in the
toroidal case remains a weak $W^{2,2}$ object unless an additional
energy identity or regularity estimate is available.

\subsection{Finite and infinite Hopf endpoints}
For a Hopf profile of geodesic curvature $\kappa$, the surface identity
$|A^0|^2=2(1+\kappa^2)\geq2$ makes quantitative nonumbilicity automatic.
At arbitrary finite initial energy, Theorem~\ref{thm:all-energy-profile}
gives a unique regular anchored endpoint, with full weak $W^{2,2}$ and
strong $C^{1,\alpha}$ convergence for every $0<\alpha<1/2$.
Its Hopf lift is the actual multiplicity-counting varifold endpoint.
Under the additional surface-energy bound $8\pi$ and finite maximality,
the limiting profile and Hopf support are embedded and counted once.
If dissipation is bounded along an endpoint sequence,
Corollary~\ref{miwf:endpoint-W42} gives a $W^{4,2}$ profile, equality of
its energy with the limiting energy, and full strong $W^{2,2}$ convergence.

For a smooth global reduced trajectory, Theorem~\ref{miwf:full} proves
smooth convergence modulo reparametrization to one regular closed
elastic curve. The proof establishes an analytic gradient inequality
in an $H^2$ normal graph chart, compares its dual norm with weighted
dissipation, and uses an energy-plus-$L^2$ estimate to prevent exit.
No nondegeneracy of the limit or uniform curvature bound along the
trajectory is assumed. The theorem gives energy and $W^{2,2}$ rates
and finite integral of the square root of dissipation. In the Hopf
sector this also gives a unique nonzero varifold limit at infinite
time, including every fixed finite torus covering. The conclusion assumes that the reduced trajectory exists for all time.

\subsection{Conformal endpoints, rank, and projectability}
Pinkall's lattice gives the standard marked Hopf modulus
$\tau=(\mathcal A+iL)/(4\pi)$, where $\mathcal A$ is a continuously transported real
spherical-area lift and $L$ is profile length. Theorem
\ref{thm:all-energy-modulus} bounds its hyperbolic speed directly by
dissipation, independently of the rank of the conformal-class map.
The joint estimate gives a unique marked endpoint at finite time.
A sharp scalar estimate also controls the variation of the square
root of profile length. For global reduced
trajectories, Theorem~\ref{miwf:full} supplies the required finite
dissipation length and hence also a marked infinite-time limit.

The rank question has a different role. Proposition~\ref{tp:prop:Hopf-rank}
expresses the Hopf differential through a two-by-two curvature Gram
matrix and computes its exact least-norm inverse. Rank is two precisely
when profile curvature is nonconstant. Its quantitative inverse can
diverge near a covered circle. By applying Jakob's positive-time
analyticity theorem \cite[Theorem~2.3(ii)]{JakobRegularity} to a
reparametrization-invariant scalar measure of the deviation from
constant mean curvature,
Theorem~\ref{tp:thm:no-circle} excludes finite circular slices from
initially noncircular trajectories. Asymptotic convergence to a
covered equator is compatible with this rigidity.

The projected velocity also retains profile information within a
fixed conformal fibre.
Theorem~\ref{tp:thm:noncollinear} constructs embedded noncircular Hopf
tori, arbitrarily close to the Clifford torus and with energy below
$8\pi$, that have the same marked modulus and noncollinear projected
velocities. Thus the projected velocity is not determined by the
marked conformal type alone. Corollary~\ref{tp:cor:no-time-change}
excludes a common autonomous equation $\dot\tau=b(\tau)$ even after
separate recalibrations of time with positive derivative. Section~\ref{tp:sec:coupled} gives local coordinates
$(\tau,k)$, with an additional profile variable $k$, and the coupled
evolution identities $\dot\tau=b(\tau,k)$, $\partial_t k=c(\tau,k)$
along smooth trajectories. Compactness and variation results in
\cite{KuwertSchaetzleBounds,KuwertSchaetzleConformal,BernardRiviere}
provide context for the distinctions between conformal control,
critical-surface compactness, and evolution considered here.

\subsection{Organization}
Sections~\ref{sec:setting}--\ref{sec:topology} treat surface transport,
endpoint uniqueness, and the finite maximal-time alternatives.
Section~\ref{sec:profiles} proves the finite-time Hopf profile results.
Section~\ref{sec:miwf-compactness} develops dissipation compactness,
local regularization, and full global reduced convergence.
Section~\ref{sec:conformal} treats marked conformal endpoints, rank,
and circular rigidity. Section~\ref{tp:sec:nonprojectability} proves the
projection obstructions. Section~\ref{tp:sec:coupled} gives the local coupled coordinates
and discusses endpoint regularity and equatorial multiplicity.

\section{Setting, conventions, and assumptions}\label{sec:setting}

Let $\Sigma$ be a smooth compact torus and let
\[
 F\colon [0,T)\times\Sigma\longrightarrow \Sph^3\subset\R^4,
 \qquad F_t:=F(t,\cdot),
\]
be a smooth solution of the M\"obius-invariant Willmore flow on its
interval of existence.  Each $F_t$, including $F_0$, is an
umbilic-free immersion.  The local well-posedness and regularity theory
for this evolution is developed in
\cite{JakobShortTime,JakobRegularity}.  Here $\Sph^3$ is the unit
sphere, $g_{\Sph^3}$ is its round metric, and $g_{\mathrm{euc}}$ is the
Euclidean metric of $\R^4$.  When the interval is maximal, we denote
its right endpoint by $T_{\max}(F_0)$.  We write
\[
 g_t:=F_t^*g_{\Sph^3}=F_t^*g_{\mathrm{euc}},
 \qquad \mu_t:=\mu_{g_t},
 \qquad \dd\mu_t:=\dd\mu_{g_t},
\]
and denote by $A_t$, $A_t^0$, and $H_t$ respectively the second
fundamental form, its trace-free part, and the trace mean-curvature
vector of $F_t$ in $\Sph^3$.  We set
\[
 \Area(F_t):=\mu_t(\Sigma).
\]
Unless an ambient space is displayed explicitly, norms and inner
products of tangent and normal tensors are taken with respect to
$g_t$ and $g_{\Sph^3}$.  We write $\nabla^{\Sph^3}$ for the spherical
gradient.
Our normalization of the
Willmore functional is
\begin{equation}
 \label{eq:willmore}
 \Will(F_t)
 :=\int_\Sigma
 \left(1+\frac14\lvert H_t\rvert^2\right)\dd\mu_t.
\end{equation}
If $H_t^{\R^4}$ denotes the Euclidean trace mean-curvature vector, then
\[
 H_t^{\R^4}=H_t-2F_t,
 \qquad
 \lvert H_t^{\R^4}\rvert^2
 =\lvert H_t\rvert^2+4,
\]
and consequently
\begin{equation}
 \label{eq:euclidean-willmore}
 \Will(F_t)=\frac14\int_\Sigma
 \lvert H_t^{\R^4}\rvert^2\dd\mu_t.
\end{equation}
Throughout the paper we write
\[
 W(t):=\Will(F_t),
 \qquad
 W_0:=W(0)=\Will(F_0).
\]

Let
\[
 G_t:=\nabla_{L^2}\Will(F_t)
\]
be the $L^2(\Sigma,\dd\mu_t)$-gradient of $\Will$, normalized by
\[
 D\Will(F_t)[\phi]
 =\int_\Sigma\langle G_t,\phi\rangle\dd\mu_t
\]
for every smooth normal variation field
$\phi\in\Gamma(N_{F_t}\Sigma)\subset\Gamma(F_t^*T\Sph^3)$.
Here $D\Will(F_t)$ is the first variation of $\Will$,
the fibre of $N_{F_t}\Sigma$ at $x$ is the orthogonal complement
of $DF_t(T_x\Sigma)$ in $T_{F_t(x)}\Sph^3$,
and $\Gamma(E)$ denotes the space of smooth sections of a bundle $E$.
The MIWF equation is
\begin{equation*}
 \partial_tF_t
 =-\frac{1}{\lvert A_t^0\rvert^4}G_t.
 \tag{MIWF}\label{eq:miwf}
\end{equation*}
Equivalently, with
\[
 a_t:=\lvert A_t^0\rvert^2,
\]
one has $\partial_tF_t=-a_t^{-2}G_t$.  These conventions give the
dissipation identity
\begin{equation}
 \label{eq:dissipation}
 -W'(t)
 =\int_\Sigma a_t^{-2}\lvert G_t\rvert^2\dd\mu_t.
\end{equation}
Pairing \eqref{eq:miwf} with $G_t$ gives this identity. In our
trace convention the spherical Willmore gradient is
\[
 G_t=\tfrac12\bigl(\Delta_{F_t}^{\perp,\Sph^3}H_t
                         +Q(A_t^0)(H_t)\bigr).
\]
Here
\[
 Q(A_t^0)(H):=g_t^{ik}g_t^{j\ell}
 \langle (A_t^0)_{ij},H\rangle(A_t^0)_{k\ell}.
\]
The operator $\Delta_{F_t}^{\perp,\Sph^3}$ is the normal-connection Laplacian
along $F_t$, and $(g_t^{ij})$ is the inverse matrix of $g_t$ in local
coordinates.
For example, write $H_t=H_{\rm sc}N$ for a chosen unit normal.
With our sign convention a variation $fN$ gives
$\delta H_{\rm sc}=\Delta f+(|A|^2+2)f$ and
$\delta\dd\mu=-H_{\rm sc}f\dd\mu$. Differentiating
\eqref{eq:willmore} and integrating $H_{\rm sc}\Delta f$ by parts
gives the scalar gradient
$\tfrac12(\Delta H_{\rm sc}+|A^0|^2 H_{\rm sc})$.
In Weiner's first variation \cite[p.~25, equation~(9)]{Weiner},
the mean-curvature vector is half the trace
\cite[p.~22]{Weiner}. On a closed surface the boundary terms vanish,
giving the normalization in \cite[equation~(4)]{Jakob}.

The additional assumption is the uniform nonumbilicity condition
\begin{equation*}
 \boxed{
 \inf_{0<t<T}\min_{x\in\Sigma}
 \lvert A_t^0(x)\rvert^2\geq c_0>0.}
 \tag{UN}\label{eq:un}
\end{equation*}
Condition~\eqref{eq:un} is imposed only when explicitly invoked.
This is stronger than requiring every individual immersion $F_t$ to
be umbilic-free.  Notice also that \eqref{eq:un} is an analytic
condition and is not the same as the ``nondegenerate case'' in
\cite{Jakob}, which means that a nonzero limiting surface has genus
one.

Throughout, ``Radon measure'' has its standard meaning and is scalar and
nonnegative unless it is explicitly described as signed or vector-valued.
In particular, every Radon measure on the compact spaces $\Sigma$ and
$\Sph^3$ is finite.

Let $G_2(\R^4):=\R^4\times G(4,2)$ be the Grassmann bundle of
unoriented two-planes over $\R^4$, and let $C_c(X)$ denote the space of
continuous, compactly supported real-valued functions on $X$.
Writing $\dd F_t$ for the differential of $F_t$, let $\Var_t$ be the
integral two-varifold in $\R^4$ associated with the immersion $F_t$,
defined by
\begin{equation}
 \label{eq:varifold-definition}
 \Var_t(\Phi)
 :=\int_\Sigma
 \Phi\bigl(F_t(x),\dd F_t(T_x\Sigma)\bigr)\dd\mu_t(x),
 \qquad \Phi\in C_c\bigl(G_2(\R^4)\bigr).
\end{equation}
For a Borel map $f$ and a Radon measure $\eta$, define the push-forward
by $(f_\#\eta)(B):=\eta(f^{-1}(B))$ for every Borel set $B$.
The weight of the immersion varifold $\Var_t$ is
\begin{equation}
 \label{eq:weight-definition}
 \lambda_t:=\lVert\Var_t\rVert=(F_t)_\#\mu_t.
\end{equation}
For a two-varifold $V$ in $\R^4$, its first variation is the linear
functional
\begin{equation}
 \label{eq:first-variation-definition}
 \delta V(X)
 :=\int_{G_2(\R^4)}\operatorname{div}_S X(x)\dd V(x,S),
 \qquad X\in C_c^1(\R^4;\R^4),
\end{equation}
where $\operatorname{div}_S X$ is the trace of $DX$ restricted to the
plane $S$.  In general $\delta V$ is a distribution.  If it is locally
bounded with respect to the supremum norm on test fields, the
Riesz--Markov theorem, applied componentwise, represents it by an
$\R^4$-valued Radon measure, still denoted by $\delta V$; its scalar
total-variation measure is denoted by $\lVert\delta V\rVert$.  For such
a varifold, the generalized mean-curvature vector $H_V^{\R^4}$ is minus the
$\R^4$-valued Radon--Nikodym derivative of the part of $\delta V$ that
is absolutely continuous with respect to $\lVert V\rVert$.  In the
case considered below, the singular part vanishes and therefore
\begin{equation}
 \label{eq:generalized-mean-curvature-definition}
 \delta V(X)
 =-\int_{\R^4}\langle H_V^{\R^4},X\rangle\dd\lVert V\rVert.
\end{equation}
For any varifold $\Var$ we use the standard abbreviation
$\spt\Var:=\spt\lVert\Var\rVert$.
We regard $\lambda_t$ as a Radon measure on $\Sph^3$ in the transport
calculation and, via the inclusion $\Sph^3\hookrightarrow\R^4$, as the
weight measure of $\Var_t$ in the compactness argument.
This definition keeps the natural multiplicity when $F_t$ is not
injective.  For a countably two-rectifiable set $M$ and a multiplicity
$\theta$, let $v(M,\theta)$ denote the rectifiable varifold carried by
$M$ with multiplicity $\theta$; let $\mathcal H^2$ denote
two-dimensional Hausdorff measure and $\mathcal H^2\llcorner M$ its
restriction to $M$.  If $F_t$ is an embedding, the area formula gives
\begin{equation}
 \label{eq:embedded-varifold}
 \Var_t=v(F_t(\Sigma),1),
 \qquad
 \lambda_t=\mathcal H^2\llcorner F_t(\Sigma).
\end{equation}

For $\varphi\in C^1(\Sph^3)$, put
\[
 \lVert\varphi\rVert_{C^1(\Sph^3)}
 :=\lVert\varphi\rVert_{L^\infty(\Sph^3)}
   +\lVert\nabla^{\Sph^3}\varphi\rVert_{L^\infty(\Sph^3)}.
\]
For positive Radon measures $\alpha$ and $\beta$ on $\Sph^3$, we use
the bounded-Lipschitz distance
\begin{equation}
 \label{eq:bl-distance}
 d_{\BL}(\alpha,\beta)
 :=\sup\left\{
 \left|\int\varphi\dd\alpha-\int\varphi\dd\beta\right|:
 \varphi\in C^1(\Sph^3),\
 \lVert\varphi\rVert_{C^1(\Sph^3)}\leq1
 \right\}.
\end{equation}
For nonempty compact sets $K,L\subset\R^4$, the Hausdorff distance
used below is the Euclidean Hausdorff distance
\begin{equation}
 \label{eq:hausdorff-distance}
 d_{\mathcal H}(K,L)
 :=\max\left\{
 \sup_{x\in K}\inf_{y\in L}|x-y|,
 \sup_{y\in L}\inf_{x\in K}|x-y|
 \right\}.
\end{equation}
Finally, $\rightharpoonup$ denotes weak convergence of Radon measures
or weak varifold convergence, according to context; the latter means
convergence against every test function in $C_c(G_2(\R^4))$.
We abbreviate $L^2(\Sigma,\dd\mu_t)$ to $L^2(\dd\mu_t)$ whenever the
domain is clear.

\subsection{Marked conformal geometry and gauges}

For $\phi\in\Diff^+(\Sigma)$, the immersion $F\circ\phi$ has the
same unparametrized image and energy. Diffeomorphisms isotopic to the
identity preserve the marked Teichm\"uller point; arbitrary ones can
change it by an element of the mapping-class group. All such
diffeomorphisms preserve its unmarked moduli point.

A fixed $M\in\Mob(\Sph^3)$ preserves the induced marked conformal
class and sends a solution of \eqref{eq:miwf} to another solution.
For time-dependent $M_t$ and $\phi_t$, however,
\begin{equation}\label{tp:eq:moving-gauge}
 \partial_t(M_t\circ F_t\circ\phi_t)
 =DM_t(F_t\circ\phi_t)
  \bigl((\partial_tF_t)\circ\phi_t+DF_t(\phi_t)\partial_t\phi_t\bigr)
  +(\partial_tM_t)(F_t\circ\phi_t).
\end{equation}
The source term involving $\partial_t\phi_t$ is tangential.
The target term involving $\partial_t M_t$ may have a normal component.
For local uniqueness in $\Sph^3$, the compact image of any smooth
slice misses a stereographic pole and continues to miss it on a short
time interval. Fixed conformal covariance then allows the Euclidean
uniqueness theorem \cite{JakobShortTime} to be applied in that chart.

The source and target actions commute. In the arguments below,
source reparametrizations fix the parametrization gauge, while any
ambient M\"obius transformation is held fixed in time.

We write
\[
 \mathcal T_1\simeq\HH=\{z\in\mathbb C:\Im z>0\},\qquad
 \mathcal M_1=\HH/\mathrm{PSL}_2(\miwfZ).
\]
For example, reparametrizations representing Dehn twists can change a
marked modulus from $\tau$ to $\tau+n$, while leaving image, energy,
and unmarked conformal class fixed. Thus geometric compactness and
moduli compactness do not, without marking control, give compactness
of a specified lift to $\mathcal T_1$.

The compactness inputs used below are applied with their stated
energy and conformal hypotheses. In particular, the marked conformal
convergence required in \cite[Proposition~2.2]{KuwertSchaetzleConformal}
is distinct from its ambient M\"obius normalization.

\section{Intrinsic and image-measure estimates}\label{sec:transport}

\begin{lemma}[Control of the geometric velocity]
\label{lem:velocity}
Assume \eqref{eq:un}.  Then, for every $0<t<T$,
\begin{equation}
 \label{eq:velocity}
 \lVert\partial_tF_t\rVert_{L^2(\Sigma,\dd\mu_t)}^2
 \leq c_0^{-2}
 \bigl(-W'(t)\bigr).
\end{equation}
\end{lemma}

\begin{proof}
By \eqref{eq:miwf},
\[
 \begin{aligned}
 \lVert\partial_tF_t\rVert_{L^2(\dd\mu_t)}^2
 &=\int_\Sigma a_t^{-4}\lvert G_t\rvert^2\dd\mu_t\\
 &=\int_\Sigma a_t^{-2}
       \bigl(a_t^{-2}\lvert G_t\rvert^2\bigr)\dd\mu_t.
 \end{aligned}
\]
Since $a_t\geq c_0$, one has $a_t^{-2}\leq c_0^{-2}$, and the
conclusion follows from \eqref{eq:dissipation}.
\end{proof}

\begin{lemma}[Evolution of the area and its push-forward]
\label{lem:transport}
Let $X_t:=\partial_tF_t$.  Then the induced area measure on the fixed
domain satisfies
\begin{equation}
 \label{eq:domain-evolution}
 \partial_t\dd\mu_t
 =-\langle H_t,X_t\rangle\dd\mu_t
 =a_t^{-2}\langle H_t,G_t\rangle\dd\mu_t.
\end{equation}
Moreover, for every $\varphi\in C^1(\Sph^3)$,
\begin{equation}
 \label{eq:pushforward-evolution}
 \begin{aligned}
 \frac{\dd}{\dd t}\int_{\Sph^3}\varphi\dd\lambda_t
 &=\int_\Sigma
 \left[
  \big\langle\nabla^{\Sph^3}\varphi(F_t),X_t\big\rangle
  -\varphi(F_t)\langle H_t,X_t\rangle
 \right]\dd\mu_t\\
 &=\int_\Sigma a_t^{-2}
 \left[
  -\big\langle\nabla^{\Sph^3}\varphi(F_t),G_t\big\rangle
  +\varphi(F_t)\langle H_t,G_t\rangle
 \right]\dd\mu_t.
 \end{aligned}
\end{equation}
\end{lemma}

\begin{proof}
The MIWF velocity is normal to the immersed surface inside $\Sph^3$.
Thus $\partial_t(g_t)_{ij}=-2\langle(A_t)_{ij},X_t\rangle$.
Taking half the metric trace gives
$\partial_t\dd\mu_t=-\langle H_t,X_t\rangle\dd\mu_t$
pointwise. Substituting $X_t=-a_t^{-2}G_t$ proves
\eqref{eq:domain-evolution}.  By the definition of the
push-forward measure,
\[
 \int_{\Sph^3}\varphi\dd\lambda_t
 =\int_\Sigma\varphi(F_t)\dd\mu_t.
\]
Differentiating this identity, applying the chain rule to
$\varphi(F_t)$, and using \eqref{eq:domain-evolution} gives the first
line of \eqref{eq:pushforward-evolution}.  Substituting the MIWF
equation gives the second line.  Thus both the motion of the immersion
and the variation of its induced area measure are accounted for.
\end{proof}

We regard $\mu_t$ as the area measure on the fixed domain, so that
\[
 \mu_t(B)=\int_B\dd\mu_t,
 \qquad B\subset\Sigma\ \text{Borel}.
\]
For a signed Radon measure $\sigma$ on $\Sigma$, we use the convention
\[
 \lVert\sigma\rVert_{\TV}
 :=\lvert\sigma\rvert(\Sigma)
 =\sup_{\substack{\psi\in C^0(\Sigma)\\
                   \lVert\psi\rVert_{L^\infty}\leq1}}
   \left|\int_\Sigma\psi\dd\sigma\right|.
\]

The area-evolution and energy-dissipation argument underlying the following estimate also appears in the proof of Lemma D.6 of \cite{DallAcquaMuellerSchaetzleSpener} for the classical Willmore flow. 
Here we adapt this argument to the M\"obius-invariant Willmore flow under condition \eqref{eq:un} and estimate the induced measures in total variation.
\begin{proposition}[Total-variation estimate on the fixed domain]
\label{prop:fixed-tv}
Assume \eqref{eq:un}, with the constant $c_0>0$ fixed there.  For every
$0<s<t<T$,
\begin{equation}
 \label{eq:fixed-tv}
 \lVert\mu_t-\mu_s\rVert_{\TV}
 \leq
 \frac{2\sqrt{W(s)}}{c_0}
 \sqrt{(t-s)\bigl(W(s)-W(t)\bigr)}.
\end{equation}
In particular, for every Borel set $B\subset\Sigma$,
\begin{equation}
 \label{eq:borel-tv}
 \lvert\mu_t(B)-\mu_s(B)\rvert
 \leq
 \frac{2\sqrt{W(s)}}{c_0}
 \sqrt{(t-s)\bigl(W(s)-W(t)\bigr)}.
\end{equation}
If $W_0\leq8\pi$, the coefficient $2\sqrt{W(s)}/c_0$ in
\eqref{eq:fixed-tv}--\eqref{eq:borel-tv} is at most
$\sqrt{32\pi}/c_0$.
\end{proposition}

\begin{proof}
Fix a smooth positive reference measure $\dd\bar\mu$ on $\Sigma$ and
write $\dd\mu_r=\rho_r\dd\bar\mu$.  From
\eqref{eq:domain-evolution} and the
fundamental theorem of calculus,
\[
 \begin{aligned}
 \lVert\mu_t-\mu_s\rVert_{\TV}
 &=\int_\Sigma\lvert\rho_t-\rho_s\rvert\dd\bar\mu\\
 &\leq\int_s^t\int_\Sigma
 a_r^{-2}\lvert H_r\rvert\lvert G_r\rvert
 \dd\mu_r\dd r.
 \end{aligned}
\]
For each $r$, Cauchy--Schwarz, \eqref{eq:un}, and the identities
\[
 \int_\Sigma\lvert H_r\rvert^2\dd\mu_r
 =4\bigl(W(r)-\Area(F_r)\bigr)\leq4W(r),
 \qquad
 \int_\Sigma a_r^{-2}\lvert G_r\rvert^2\dd\mu_r=-W'(r),
\]
give
\[
 \begin{aligned}
 \int_\Sigma a_r^{-2}\lvert H_r\rvert\lvert G_r\rvert
 \dd\mu_r
 &\leq
 \left(\int_\Sigma a_r^{-2}\lvert H_r\rvert^2\dd\mu_r\right)^{1/2}
 \left(\int_\Sigma a_r^{-2}\lvert G_r\rvert^2\dd\mu_r\right)^{1/2}\\
 &\leq
 \frac{2\sqrt{W(r)}}{c_0}
 \sqrt{-W'(r)}\\
 &\leq
 \frac{2\sqrt{W(s)}}{c_0}\sqrt{-W'(r)}.
 \end{aligned}
\]
Here we used $W(r)\leq W(s)$ for $r\in[s,t]$.
A second application of Cauchy--Schwarz, now in time, yields
\[
 \begin{aligned}
 \lVert\mu_t-\mu_s\rVert_{\TV}
 &\leq\frac{2\sqrt{W(s)}}{c_0}
 \int_s^t\sqrt{-W'(r)}\dd r\\
 &\leq\frac{2\sqrt{W(s)}}{c_0}
 \sqrt{t-s}
 \sqrt{\int_s^t-W'(r)\dd r},
 \end{aligned}
\]
which is \eqref{eq:fixed-tv}.  The estimate for a Borel set follows from
$\lvert(\mu_t-\mu_s)(B)\rvert\leq
\lVert\mu_t-\mu_s\rVert_{\TV}$.
\end{proof}

The preceding estimate controls measures on the fixed domain.  To
control their push-forwards one must also estimate the motion of the
maps $F_t$.

\begin{proposition}[Bounded-Lipschitz estimate for the image measures]
\label{prop:image-bl}
Under \eqref{eq:un}, for every $0<s<t<T$ and every
$\varphi\in C^1(\Sph^3)$,
\begin{equation}
 \label{eq:image-bl}
 \left|\int_{\Sph^3}\varphi\dd\lambda_t
       -\int_{\Sph^3}\varphi\dd\lambda_s\right|
 \leq
 \frac{2\sqrt{W(s)}}{c_0}\,
 \lVert\varphi\rVert_{C^1(\Sph^3)}
 \sqrt{(t-s)\bigl(W(s)-W(t)\bigr)},
\end{equation}
where the $C^1(\Sph^3)$-norm is the one fixed above
\eqref{eq:bl-distance}.
\end{proposition}

\begin{proof}
Fix $0<s<t<T$. For each intermediate time $r\in[s,t]$,
the transport identity \eqref{eq:pushforward-evolution} and
Cauchy--Schwarz give
\[
 \begin{aligned}
 \left|\frac{\dd}{\dd r}
       \int_{\Sph^3}\varphi\dd\lambda_r\right|
 &\leq
 \Bigl(
   \lVert\nabla^{\Sph^3}\varphi\rVert_{L^\infty}
    \Area(F_r)^{1/2}
   +\lVert\varphi\rVert_{L^\infty}
    \lVert H_r\rVert_{L^2(\dd\mu_r)}
 \Bigr)
 \lVert X_r\rVert_{L^2(\dd\mu_r)}.
 \end{aligned}
\]
Since
\[
 \Area(F_r)\leq W(r),
 \qquad
 \lVert H_r\rVert_{L^2(\dd\mu_r)}
 =2\sqrt{W(r)-\Area(F_r)}\leq2\sqrt{W(r)},
\]
the velocity estimate \eqref{eq:velocity} and energy monotonicity
$W(r)\leq W(s)$ yield
\[
 \left|\frac{\dd}{\dd r}
       \int_{\Sph^3}\varphi\dd\lambda_r\right|
 \leq
 \frac{2\sqrt{W(s)}}{c_0}
 \lVert\varphi\rVert_{C^1(\Sph^3)}\sqrt{-W'(r)}.
\]
Integration over $[s,t]$ and Cauchy--Schwarz in time prove
\eqref{eq:image-bl}.
\end{proof}

\begin{remark}[Intrinsic and image measures]
The measure $\mu_t$ lives on the fixed parameter torus, whereas
$\lambda_t=(F_t)_\#\mu_t$ also records the motion of $F_t$.
Formula~\eqref{eq:pushforward-evolution} controls this additional
transport term. The estimates use $\Area(F_t)\leq W(t)\leq W_0$
directly in $\Sph^3$; its Euclidean inclusion is used for varifold
compactness in Section~\ref{sec:unique}.
\end{remark}

\begin{remark}[Finite and infinite time]
On a finite time interval, \eqref{eq:fixed-tv} and
\eqref{eq:image-bl} make the two measure families Cauchy at the
endpoint. For infinite time, the corresponding criterion is finite
dissipation length, as stated in
Proposition~\ref{prop:infinite-finite-length}.
\end{remark}

\section{Proof of the unique-limit theorem}\label{sec:unique}

\begin{lemma}[A rectifiable varifold is determined by its weight]
\label{lem:weight-determines-varifold}
Let $V_1$ and $V_2$ be rectifiable two-varifolds in $\R^4$.  If
$\lVert V_1\rVert=\lVert V_2\rVert$, then $V_1=V_2$.
\end{lemma}

\begin{proof}
Let $\lambda$ be the common weight. In a rectifiable representation,
partition a carrier of $\lambda$, up to a null set, into measurable
pieces $E$ contained in $C^1$ surfaces. At $\lambda$-almost every
$x\in E$, differentiation on that surface shows that the rescalings
$r^{-2}(y\mapsto(y-x)/r)_\#(\lambda\llcorner E)$ converge locally
to a positive multiple of plane measure on the surface tangent at $x$.
Moreover, $\lambda(E\cap B_r(x))/\lambda(B_r(x))\to1$ by differentiation
with respect to $\lambda$. Since $r^{-2}\lambda(E\cap B_r(x))$ has a
positive finite limit, the complementary mass in $B_r(x)$ is $o(r^2)$.
Thus the rescaling limit of the full weight determines the same plane
$P(x)$, independently of the rectifiable representation. Each varifold
therefore acts by
\[
 V_i(\Phi)=\int\Phi(x,P(x))\,\dd\lambda(x),\qquad i=1,2,
\]
for every $\Phi\in C_c(G_2(\R^4))$. Thus $V_1=V_2$; compare the
rectifiable-varifold representation in \cite{Allard}.
\end{proof}

\begin{proof}[Proof of Theorem~\ref{thm:main}]
Proposition~\ref{prop:fixed-tv} and the finiteness of $T$ show that
$(\mu_t)_{t<T}$ is Cauchy in total variation.  The space of signed
Radon measures on $\Sigma$ is complete in the total-variation norm,
and its nonnegative cone is closed.  Hence there is a Radon
measure $\mu_T$ such that $\mu_t\to\mu_T$ in total variation.
Letting $t\nearrow T$ in \eqref{eq:fixed-tv} gives
\eqref{eq:main-endpoint-tv}.

Proposition~\ref{prop:image-bl} implies
\begin{equation}
 \label{eq:bl-cauchy}
 d_{\BL}(\lambda_s,\lambda_t)
 \leq\frac{2\sqrt{W(s)}}{c_0}
 \sqrt{(t-s)\bigl(W(s)-W(t)\bigr)}.
\end{equation}
Because $T<\infty$, the right-hand side tends to zero as
$s,t\nearrow T$.  Moreover,
\[
 \lambda_t(\Sph^3)=\Area(F_t)\leq W_0.
\]
Since $\Sph^3$ is compact, the Radon measures $\lambda$ with
$\lambda(\Sph^3)\leq W_0$ form a sequentially compact set for weak-*
convergence in $C(\Sph^3)^*$.  Hence every sequence $t_j\nearrow T$
has a subsequence along which $\lambda_{t_j}$ converges weak-* in
$C(\Sph^3)^*$.  On this uniformly bounded family the bounded-Lipschitz
distance metrizes this weak-* convergence, and estimate
\eqref{eq:bl-cauchy} forces
any two such cluster measures to coincide.  It follows that there is a
unique Radon measure $\lambda_T$ and that the full family
satisfies
\begin{equation}
 \label{eq:weight-limit}
 \lambda_t\rightharpoonup\lambda_T
 \qquad\text{as }t\nearrow T.
\end{equation}
In fact the convergence holds in $d_{\BL}$, and letting
$t\nearrow T$ in \eqref{eq:bl-cauchy} proves
\eqref{eq:main-endpoint-bl}.

We next use varifold compactness.  For every $t$, $\Var_t$ is an
integral varifold supported in $\Sph^3$ and
\[
 \lVert\Var_t\rVert(\R^4)=\Area(F_t)\leq W_0.
\]
Since $\Var_t$ is induced by a smooth closed immersion, its first
variation is represented by the Euclidean mean-curvature vector and
satisfies
\begin{equation}
 \label{eq:first-variation}
 \begin{aligned}
 \lVert\delta\Var_t\rVert(\R^4)
 &\leq\int_\Sigma\lvert H_t^{\R^4}\rvert\dd\mu_t\\
 &\leq
 \lVert H_t^{\R^4}\rVert_{L^2(\dd\mu_t)}
 \Area(F_t)^{1/2}
 \leq2W_0.
 \end{aligned}
\end{equation}
The area formula gives positive integer density
at $\lVert\Var_t\rVert$-almost every point, also for immersed slices
with overlapping sheets. Allard's integral compactness theorem
\cite{Allard}, in the form recorded in
\cite[Chapter~8, Theorem~5.8 and Remark~5.9]{SimonGMT},
therefore shows
that every sequence $t_j\nearrow T$ admits a subsequence $t_{j_k}$ and
an integral two-varifold $\Var$ such that
\begin{equation}
 \label{eq:cluster-varifold}
 \Var_{t_{j_k}}\rightharpoonup\Var.
\end{equation}
Taking weight measures in \eqref{eq:cluster-varifold} and using
\eqref{eq:weight-limit} gives
\begin{equation}
 \label{eq:cluster-weight}
 \lVert\Var\rVert=\lambda_T.
\end{equation}

If $\lambda_T=0$, set $\Var_T:=0$; then
\eqref{eq:cluster-weight} makes every cluster varifold the zero
varifold.  If $\lambda_T\neq0$, all cluster varifolds are integral,
and hence rectifiable, by Allard compactness and have the same weight
$\lambda_T$ by \eqref{eq:cluster-weight}.  Lemma
\ref{lem:weight-determines-varifold} therefore shows that they
coincide; denote their common value by $\Var_T$.

If full convergence failed, some sequence $t_j\nearrow T$ would
remain outside a fixed neighborhood of $\Var_T$. Integral compactness
provides a convergent subsequence, whose limit must be $\Var_T$ by
the uniqueness just proved. This contradiction proves
\eqref{eq:main-convergence}.

It remains to prove the endpoint curvature assertion.  For every
$X\in C_c^1(\R^4;\R^4)$, the definition
\eqref{eq:first-variation-definition} and varifold convergence give
\[
 \delta\Var_t(X)\longrightarrow\delta\Var_T(X),
\]
because $(x,S)\mapsto\operatorname{div}_S X(x)$ is continuous and
compactly supported on $G_2(\R^4)$.  Since
\[
 \delta\Var_t(X)
 =-\int_\Sigma
 \langle H_t^{\R^4},X(F_t)\rangle\dd\mu_t,
\]
Cauchy--Schwarz, \eqref{eq:euclidean-willmore},
\eqref{eq:weight-limit}, and $W(t)\to W_T$ imply
\begin{align*}
 |\delta\Var_T(X)|
 &\leq
 \limsup_{t\nearrow T}
 \left(\int_\Sigma|H_t^{\R^4}|^2\dd\mu_t\right)^{1/2}
 \left(\int_{\Sph^3}|X|^2\dd\lambda_t\right)^{1/2}\\
 &=2\sqrt{W_T}
 \left(\int_{\Sph^3}|X|^2\dd\lambda_T\right)^{1/2}.
\end{align*}
Since $C_c^1(\R^4;\R^4)$ is dense in
$L^2(\lambda_T;\R^4)$, the preceding estimate shows that
$\delta\Var_T$ extends uniquely to a bounded linear functional on this
space.  The Riesz representation theorem for this Hilbert space
supplies a vector field
$H_{\Var_T}^{\R^4}\in L^2(\lambda_T;\R^4)$ satisfying
\[
 \delta\Var_T(X)
 =-\int_{\R^4}\langle H_{\Var_T}^{\R^4},X\rangle\dd\lambda_T,
 \qquad
 \lVert H_{\Var_T}^{\R^4}\rVert_{L^2(\lambda_T)}
 \leq2\sqrt{W_T}.
\]
By \eqref{eq:generalized-mean-curvature-definition}, this vector field
is the generalized Euclidean mean curvature of $\Var_T$, the singular
part of $\delta\Var_T$ vanishes, and the norm bound is precisely
\eqref{eq:main-curvature-lsc}.  This is the standard
lower-semicontinuity mechanism for weak mean curvature under varifold
convergence; compare \cite[(1.1) and the introduction]
{SchaetzleLowerSemicontinuity}.
\end{proof}

\subsection{Relaxing uniform nonumbilicity}

The uniform constant in \eqref{eq:un} can be replaced by a weighted
transport condition.  This extension does not require a low-energy
assumption and applies to a prescribed smooth trajectory.

\begin{proposition}[Weighted endpoint criterion]
\label{prop:weighted-endpoint}
Let $F_t$, $0\leq t<T<\infty$, be a smooth umbilic-free MIWF, and put
\[
 D(t)=-W'(t),\qquad a_t=|A_t^0|^2,\qquad
 B(t)=\int_\Sigma\frac{1+|H_t|^2}{a_t^2}\dd\mu_t.
\]
If
\begin{equation}\label{eq:weighted-integrability}
 \int_0^T\sqrt{D(t)B(t)}\dd t<\infty,
\end{equation}
then all convergence and endpoint mean-curvature conclusions of
Theorem~\ref{thm:main} hold, with the quantitative estimates
\begin{align}
 \|\mu_T-\mu_s\|_{\TV}
 &\leq\int_s^T\sqrt{DB}\dd r,\label{eq:weighted-TV}\\
 d_{\BL}(\lambda_T,\lambda_s)
 &\leq\int_s^T\sqrt{DB}\dd r.
 \label{eq:weighted-BL}
\end{align}
A sufficient condition is
\begin{equation}\label{eq:integrable-minimum}
 \int_0^T m(t)^{-2}\dd t<\infty,
 \qquad m(t)=\min_\Sigma|A_t^0|^2.
\end{equation}
In particular this permits $m(t)\geq c(T-t)^\alpha$ near $T$, with
$c>0$ and $0<\alpha<1/2$.
\end{proposition}
\begin{proof}
Writing $V_t=\partial_tF_t$, the dissipation identity is
$D=\int a^2|V|^2\dd\mu$.  The area identity in
Lemma~\ref{lem:transport}, expressed relative to a fixed smooth
positive density, gives
\[
 \|\mu_t-\mu_s\|_{\TV}
 \leq\int_s^t\int_\Sigma|H_r||V_r|\dd\mu_r\dd r
 \leq\int_s^t\sqrt{DB}\dd r.
\]
For every test function admissible in \eqref{eq:bl-distance}, the
push-forward identity in that lemma and weighted Cauchy--Schwarz give,
with $b=\|\nabla\varphi\|_\infty$ and $c=\|\varphi\|_\infty$,
\[
 \left|\frac{\dd}{\dd r}\int\varphi\dd\lambda_r\right|
 \leq\sqrt D\left(\int_\Sigma
                  \frac{(b+c|H|)^2}{a^2}\dd\mu\right)^{1/2}
 \leq\sqrt{DB}.
\]
Indeed $b+c\leq1$, and
$(b+c|H|)^2\leq(b^2+c^2)(1+|H|^2)\leq1+|H|^2$.
Thus both measure families are Cauchy under
\eqref{eq:weighted-integrability}.  Completeness in total variation
gives $\mu_T$, and compactness of the sphere and the mass bound
$\lambda_t(\Sph^3)\leq W_0$ give the unique limit $\lambda_T$ in
$d_{\BL}$.  Letting $t\uparrow T$ proves the two tail estimates.

The compactness part of the proof of Theorem~\ref{thm:main} used
only the bounds $\|\Var_t\|(\R^4)\leq W_0$ and
$\|\delta\Var_t\|(\R^4)\leq2W_0$, which still hold here.
Consequently every time sequence has an integral-varifold cluster
limit.  Each cluster has weight $\lambda_T$, so
Lemma~\ref{lem:weight-determines-varifold} makes the cluster unique
and gives convergence of the full family.  Finally, for
$X\in C_c^1(\R^4;\R^4)$, first variation and Cauchy--Schwarz give
\[
 |\delta\Var_T(X)|\leq2\sqrt{W_T}\|X\|_{L^2(\lambda_T)}.
\]
The Hilbert-space Riesz theorem therefore supplies the generalized
mean curvature with squared norm at most $4W_T$, exactly as asserted
in \eqref{eq:main-curvature-lsc}.

For the sufficient condition, energy decrease yields
\[
 B(t)\leq m(t)^{-2}
       \left(\Area(F_t)+\int_\Sigma|H_t|^2\dd\mu_t\right)
       \leq4W_0m(t)^{-2}.
\]
Since $\int_0^T D\dd t=W_0-W_T<\infty$, Cauchy--Schwarz in time
proves \eqref{eq:weighted-integrability} from
\eqref{eq:integrable-minimum}.  The final power bound makes
$m^{-2}\leq c^{-2}(T-t)^{-2\alpha}$ near $T$; away from $T$,
smoothness and nonumbilicity bound $m^{-2}$.  It is therefore integrable.
\end{proof}

\begin{remark}[Weighted transport]
Under \eqref{eq:weighted-integrability}, the uniqueness argument in
Corollary~\ref{cor:classification} applies to the low-energy
subsequential alternatives even without a uniform positive lower
bound for $|A^0|^2$. Endpoint energy identity is a separate issue;
the reduced-flow criterion is given in
Corollary~\ref{miwf:endpoint-W42}.
\end{remark}

\subsection{A finite-dissipation-length criterion at infinite time}

The finite-time hypothesis in Theorem~\ref{thm:main} can be replaced
by an explicit integrability condition on the dissipated energy.

\begin{proposition}[Unique limit under finite dissipation length]
\label{prop:infinite-finite-length}
Let $F\colon[0,\infty)\times\Sigma\to\Sph^3$ be a smooth MIWF of a
compact torus.  Assume that there is one constant $c_0>0$ such that
\[
 \inf_{t>0}\min_{x\in\Sigma}|A_{F_t}^0(x)|^2\geq c_0,
\]
and, in addition, that
\begin{equation}
 \label{eq:finite-dissipation-length}
 \int_0^\infty\sqrt{-W'(r)}\dd r<\infty.
\end{equation}
Then there are nonnegative Radon measures $\mu_\infty$ on $\Sigma$ and
$\lambda_\infty$ on $\Sph^3$ and a unique integral two-varifold
$\Var_\infty$ in $\R^4$, possibly the zero varifold, such that
\[
 \mu_t\longrightarrow\mu_\infty\quad\text{in total variation},
 \qquad
 \lambda_t\longrightarrow\lambda_\infty
 \quad\text{in }d_{\BL},
 \qquad
 \Var_t\rightharpoonup\Var_\infty
\]
as $t\to\infty$, and $\lambda_\infty=\lVert\Var_\infty\rVert$.
Writing $W_\infty:=\lim_{t\to\infty}W(t)$, one also has
\[
 H_{\Var_\infty}^{\R^4}
 \in L^2(\lVert\Var_\infty\rVert;\R^4),
 \qquad
 \frac14\int_{\R^4}|H_{\Var_\infty}^{\R^4}|^2
 \dd\lVert\Var_\infty\rVert\leq W_\infty.
\]
\end{proposition}

\begin{proof}
Before the final Cauchy--Schwarz inequality in time, the proofs of
Propositions~\ref{prop:fixed-tv} and~\ref{prop:image-bl} give, for
$0<s<t$,
\begin{align*}
 \lVert\mu_t-\mu_s\rVert_{\TV}
 &\leq\frac{2\sqrt{W(s)}}{c_0}
 \int_s^t\sqrt{-W'(r)}\dd r,\\
 d_{\BL}(\lambda_t,\lambda_s)
 &\leq\frac{2\sqrt{W(s)}}{c_0}
 \int_s^t\sqrt{-W'(r)}\dd r.
\end{align*}
Assumption \eqref{eq:finite-dissipation-length} makes both right-hand
sides tend to zero as $s,t\to\infty$.  Completeness of the signed
Radon measures on $\Sigma$ in total variation gives $\mu_\infty$.
Compactness of $\Sph^3$, the uniform mass bound
$\lambda_t(\Sph^3)\leq W_0$, and the second Cauchy estimate give the
unique limit $\lambda_\infty$ in $d_{\BL}$.

The mass bound $\lVert\Var_t\rVert(\R^4)\leq W_0$ and the
first-variation bound \eqref{eq:first-variation} are uniform for all
$t>0$.  Allard
compactness therefore makes every sequence $t_j\to\infty$ admit an
integral-varifold cluster limit.  Each such limit has weight
$\lambda_\infty$, so Lemma~\ref{lem:weight-determines-varifold} shows
that all cluster varifolds coincide.  The usual contradiction argument
then gives convergence of the full family to a unique integral
varifold $\Var_\infty$ with weight $\lambda_\infty$.

Finally, $W(t)$ decreases to $W_\infty$.  Applying
\eqref{eq:first-variation-definition} along an arbitrary sequence
$t_j\to\infty$ and repeating the $L^2(\lambda_\infty)$ estimate in the
last part of the proof of Theorem~\ref{thm:main} gives
\[
 |\delta\Var_\infty(X)|
 \leq2\sqrt{W_\infty}
 \lVert X\rVert_{L^2(\lambda_\infty)}.
\]
The Hilbert-space Riesz theorem and
\eqref{eq:generalized-mean-curvature-definition} prove the asserted
mean-curvature statement and its bound.
\end{proof}

\begin{remark}[Dissipation length in the Hopf sector]
For a global Hopf trajectory with fixed additional covering degree
$d$, Lemma~\ref{tp:lem:Hopf-reconstruction} and the pointwise identity
$|A^0|^2=2(1+\kappa^2)$ in \eqref{eq:sheetwise-hopf-gradient}, together with
Theorem~\ref{miwf:full}, verify the hypotheses of
Proposition~\ref{prop:infinite-finite-length}.
\end{remark}

\section{Unique endpoint alternatives at finite maximal time}\label{sec:topology}

The results in this section concern the maximal time of existence,
because this is the endpoint to which the compactness theorems of
\cite{Jakob} apply.  We combine Theorem~\ref{thm:main} directly with those results
in $\Sph^3$; no stereographic projection is introduced.

\begin{lemma}[Umbilic event versus uniform nonumbilicity]
\label{lem:umbilic-alternative}
Let $F\colon[0,T)\times\Sigma\to\Sph^3$ be a smooth MIWF with
$T<\infty$.  Then exactly one of the following alternatives holds:
\begin{enumerate}[label=\textup{(\roman*)}]
 \item there are $t_j\nearrow T$ and $x_j\in\Sigma$ such that
 \begin{equation}
  \label{eq:umbilic-event}
  |A^0_{F_{t_j}}(x_j)|^2\longrightarrow0;
 \end{equation}
 \item there is $c_0>0$ such that
 \begin{equation}
  \label{eq:uniform-nonumbilic-alternative}
  \inf_{0<t<T}\min_{x\in\Sigma}|A^0_{F_t}(x)|^2\geq c_0.
 \end{equation}
\end{enumerate}
\end{lemma}

\begin{proof}
Set
\[
 m(t):=\min_{x\in\Sigma}|A^0_{F_t}(x)|^2.
\]
The function $m$ is continuous and positive on $[0,T)$, since the
flow is smooth and each $F_t$ is umbilic-free.  Consequently, for every
$0<\varepsilon<T$ it has a positive minimum on the compact interval
$[0,T-\varepsilon]$.  If
\eqref{eq:uniform-nonumbilic-alternative} fails, choose
$\tau_j\in(0,T)$ with $m(\tau_j)<1/j$.  The preceding compact-interval
bound forces $\tau_j\to T$.  After passing to an increasing subsequence
and choosing $x_j$ at which $m(\tau_j)$ is attained, one obtains
\eqref{eq:umbilic-event}.  Conversely,
\eqref{eq:uniform-nonumbilic-alternative} plainly excludes
\eqref{eq:umbilic-event}.
\end{proof}

\begin{corollary}[Sequence-independent endpoint topology]
\label{cor:classification}
Assume the hypotheses of Theorem~\ref{thm:main} with
\begin{equation}
 \label{eq:maximal-time}
 T=T_{\max}(F_0)<\infty
 \qquad\text{and}\qquad
 W(0)\leq8\pi.
\end{equation}
Then exactly one of the following mutually exclusive alternatives
holds:
\begin{enumerate}[label=\textup{(\alph*)}]
 \item $\Var_T=0$;
 \item $\Var_T\neq0$ and $\spt\Var_T$ is a closed embedded orientable
 Lipschitz surface of genus zero;
 \item $\Var_T\neq0$ and $\spt\Var_T$ is a closed embedded orientable
 Lipschitz surface of genus one.
\end{enumerate}
In the nonzero alternatives, $\Var_T$ has unit density.
The alternative and the limiting varifold are independent of the
sequence of times tending to $T$.
\end{corollary}

\begin{proof}
If $W(t_0)=8\pi$ for some $t_0>0$, monotonicity and $W(0)\leq8\pi$
give $W(t)=8\pi$ on $[0,t_0]$.  The dissipation identity then implies
$G_t=0$ for every $t\in(0,t_0)$, and hence $F_t$ is stationary on that
interval.  By continuity it is stationary up to $t=t_0$.  The stationary
solution is a smooth global solution, and uniqueness in the short-time
theory for the MIWF \cite{JakobShortTime} therefore makes it the
maximal continuation of $F_{t_0}$.  This
contradicts $T_{\max}(F_0)<\infty$.
Consequently,
\[
 W(t)<8\pi\qquad\text{for every }0<t<T,
\]
and the Li--Yau inequality \cite{LiYau}, as used in
\cite[proof of Theorem~1.1]{Jakob}, shows that every $F_t$, $t>0$, is
an embedding.
Thus \eqref{eq:embedded-varifold} holds.

Given any sequence $t_j\nearrow T$, Theorem~1.1(1) of \cite{Jakob}
provides a subsequence and an integral two-varifold
$\widetilde{\Var}$ such that the multiplicity-one image measures
converge weakly to $\lVert\widetilde{\Var}\rVert$.  If this limit is
nonzero, the same theorem gives unit density and says that its support is
a closed embedded orientable Lipschitz surface of genus zero or one.
Since these image measures are precisely $\lambda_t$, one has
$\lVert\widetilde{\Var}\rVert=\lambda_T
=\lVert\Var_T\rVert$.  Both varifolds are rectifiable, so
Lemma~\ref{lem:weight-determines-varifold} gives
$\widetilde{\Var}=\Var_T$.
This proves the endpoint trichotomy and its sequence independence.
\end{proof}

\begin{corollary}[Exhaustive finite-time breakdown alternative]
\label{cor:four-way}
Let $F\colon[0,T_{\max})\times\Sigma\to\Sph^3$ be a maximal smooth
MIWF with
\[
 T_{\max}=T_{\max}(F_0)<\infty,
 \qquad W(0)\leq8\pi.
\]
Then exactly one of the following mutually exclusive alternatives
holds:
\begin{enumerate}[label=\textup{(\Alph*)}]
 \item there are $t_j\nearrow T_{\max}$ and $x_j\in\Sigma$ such that
 \[
  |A^0_{F_{t_j}}(x_j)|^2\longrightarrow0;
 \]
 \item uniform nonumbilicity holds and the unique limit is zero:
 $\Var_{T_{\max}}=0$;
 \item uniform nonumbilicity holds, $\Var_{T_{\max}}\neq0$, and
 $\spt\Var_{T_{\max}}$ is a closed embedded orientable Lipschitz
 surface of genus zero;
 \item uniform nonumbilicity holds, $\Var_{T_{\max}}\neq0$, and
 $\spt\Var_{T_{\max}}$ is a closed embedded orientable Lipschitz
 surface of genus one.
\end{enumerate}
In alternatives \textup{(C)} and \textup{(D)}, the limit has unit
density.
\end{corollary}

\begin{proof}
Lemma~\ref{lem:umbilic-alternative} gives the mutually exclusive
alternative between \textup{(A)} and uniform nonumbilicity.  In the
latter case Theorem~\ref{thm:main} and
Corollary~\ref{cor:classification} give exactly \textup{(B)},
\textup{(C)}, or \textup{(D)}, including the unit-density statement.
\end{proof}

\begin{corollary}[Full Hausdorff convergence in the nonzero case]
\label{cor:hausdorff}
Under the hypotheses of Corollary~\ref{cor:classification}, if
$\Var_T\neq0$, then
\begin{equation}
 \label{eq:hausdorff}
 F_t(\Sigma)\longrightarrow\spt\Var_T
 \qquad\text{in Hausdorff distance as }t\nearrow T.
\end{equation}
\end{corollary}

\begin{proof}
If \eqref{eq:hausdorff} failed, there would exist $\varepsilon>0$ and a
sequence $t_j\nearrow T$ such that
\[
 d_{\mathcal H}\bigl(F_{t_j}(\Sigma),\spt\Var_T\bigr)
 \geq\varepsilon.
\]
After passing to a subsequence, Theorem~1.1(1) of \cite{Jakob} gives
weak convergence of the image measures and, because the limit is
nonzero, Hausdorff convergence to the support of the limiting
varifold.  The latter is $\Var_T$ by Theorem~\ref{thm:main}, a contradiction.
\end{proof}

\begin{corollary}[Parametrization of the unique toroidal limit]
\label{cor:parametrization}
Assume the hypotheses of Corollary~\ref{cor:classification} and suppose that
$\Var_T\neq0$ and $\spt\Var_T$ has genus one.  Then there exists a
uniformly conformal bi-Lipschitz homeomorphism
\[
 f\in W^{2,2}(\Sigma,\R^4)\cap W^{1,\infty}(\Sigma,\R^4),
 \qquad
 f\colon\Sigma\longrightarrow M:=\spt\Var_T,
\]
and there are a smooth zero-scalar-curvature, unit-volume metric
$g_{\mathrm{poin}}$ on $\Sigma$ and a function $u\in L^\infty(\Sigma)$
such that
\begin{equation}
 \label{eq:conformal-limit}
 f^*g_{\mathrm{euc}}=e^{2u}g_{\mathrm{poin}}.
\end{equation}
This identity is understood almost everywhere on \(\Sigma\).
Moreover,
\begin{equation}
 \label{eq:unit-density-limit}
 \Var_T=v(M,1),
\end{equation}
and
\begin{equation}
 \label{eq:limit-measure}
 \lVert\Var_T\rVert
 =f_\#\mu_{f^*g_{\mathrm{euc}}}
 =\mathcal H^2\llcorner M.
\end{equation}
The measure identity \eqref{eq:limit-measure} is independent of the
sequence used to construct the parametrization $f$.
\end{corollary}

\begin{proof}
Apply Theorem~1.1(2) and formula (9) of \cite{Jakob} to any sequence
realizing the nonzero genus-one limit.  They provide $f$, the conformal
identity \eqref{eq:conformal-limit}, and the measure identity
\eqref{eq:limit-measure}.
Since the limit has unit density, the corresponding full rectifiable
varifold is $v(M,1)$.  Equivalently, it is obtained by integrating over
$\Sigma$ with the area measure of $f^*g_{\mathrm{euc}}$ and with the
approximate tangent plane carried by $f$.  Theorem~\ref{thm:main} makes this
geometric varifold independent of the sequence.
\end{proof}

\begin{remark}[Geometric uniqueness and parametrization]
Theorem~\ref{thm:main} and Corollary~\ref{cor:parametrization} prove uniqueness of the unparametrized
geometric limit.  They do not produce a unique limiting map
$F_T\colon\Sigma\to\Sph^3$.  The parametrizations in
\cite[Theorem~1.1(2)]{Jakob} are obtained after sequence-dependent
diffeomorphisms of the domain.  Different choices can therefore give
different parametrizations of the same uniquely determined varifold.
\end{remark}

\begin{remark}[Nonumbilicity at a weak endpoint]
\label{rem:nonumbilic-transfer}
The uniform lower bound along the smooth trajectory does not, by
itself, imply
\[
 \operatorname*{ess\,inf}_{x\in\Sigma}|A^0_f(x)|^2>0
\]
for the parametrization supplied by
Corollary~\ref{cor:parametrization}, where $A_f^0$ denotes the weak
trace-free second fundamental form of $f$.  The convergence in
\cite[Theorem~1.1(2)]{Jakob} is only weak in $W^{2,2}$ and weak-* in
$W^{1,\infty}$ after sequence-dependent reparametrizations.  Since
$A^0$ depends on second derivatives and on the induced metric, neither
this weak convergence nor the accompanying $C^0$ convergence transmits
a pointwise lower bound.  Such a conclusion would require an
additional no-defect argument, for example a strong $W^{2,2}$
compactness statement with sufficient control of the induced metrics.
The stronger convergence in \cite[Theorem~1.1(3)]{Jakob} requires an
energy identity which is not among the hypotheses of
Corollary~\ref{cor:parametrization}.
\end{remark}

\subsection{Endpoint regularity and continuation}
The genus-one endpoint in Corollary~\ref{cor:four-way} has the
$W^{2,2}\cap W^{1,\infty}$ parametrization of
Corollary~\ref{cor:parametrization}. Smooth continuation requires
control in a Cauchy trace space and compatibility with the incoming
trajectory. The maximal-regularity theory in
\cite{JakobRegularity} uses
\[
 W^{1,p}(0,T;L^p)\cap L^p(0,T;W^{4,p}),\qquad
 F_0\in W^{4-4/p,p},\qquad p>3.
\]
To exclude alternative \textup{(D)} by continuation, one would need a
gauge for the incoming trajectory with an umbilic-free trace at
$T_{\max}$ and a local solution from that trace which extends the
trajectory. Subsequence-dependent reparametrizations give only the
weak endpoint of Corollary~\ref{cor:parametrization}. Quantitative
nonumbilicity must also pass to the limit
(Remark~\ref{rem:nonumbilic-transfer}). A restart from a merely
$W^{2,2}$ trace would require a lower-regularity Cauchy theory;
compare \cite[Remarks~7.1--7.3]{Jakob}.

\section{Finite-time reduced Hopf profiles}\label{sec:profiles}

Let \(\pi\colon\Sph^3\to\Sph^2\) denote the standard Hopf fibration.
For a regular closed parametrized curve $\gamma$, its pullback circle
bundle is
\[
 P_\gamma:=\{(x,z)\in\Sph^1\times\Sph^3:\pi(z)=\gamma(x)\},
 \qquad X_\gamma(x,z):=z.
\]
It is a smooth torus: a horizontal lift corrected by its terminal
fibre phase gives a global section and trivializes this circle bundle.
The map $X_\gamma$ is an immersion.  We call $X_\gamma$, and its composition
with a finite covering of $P_\gamma$, a Hopf immersion.  This definition
retains self-intersections and repeated traversals of $\gamma$.
Write $L(\gamma)=\int\dd s_\gamma$ and
$\miwfE(\gamma)=\int(1+\kappa_\gamma^2)\dd s_\gamma$, where
$\kappa_\gamma$ is the signed geodesic curvature on the unit sphere.
If $F=X_\gamma\circ q$ and the covering $q$ has degree $d\geq1$, then
the fibre integration and local curvature formulas give
\begin{equation}\label{eq:hopf-cover-energy}
 \Area(F)=d\pi L(\gamma),\qquad
 \Will(F)=d\pi \miwfE(\gamma).
\end{equation}
Here the length and energy of $\gamma$ already count its traversals;
$d$ counts the sheets of the additional covering of the parameter torus.
Indeed, the horizontal dilation is $2$, each fibre has length $2\pi$, and
the trace mean curvature on each local sheet is $2\kappa_\gamma$.
These are local identities from \cite[pp.~380 and 385]{Pinkall} and
\cite[Proposition~4.3]{Jakob}; the degree enters only on integration.
For marked conformal classes we give the covering torus the pullback
orientation and use a positively oriented homology basis.
We use this convention also for the global Hopf trajectories below.
In the finite-time, low-energy statements, embeddedness at positive
times ensures that the torus and the profile are each counted once.

\begin{equation}\label{tp:eq:Hopf-geometry}
 A=\begin{pmatrix}2\kappa&1\\1&0\end{pmatrix},\qquad
 A^0=\begin{pmatrix}\kappa&1\\1&-\kappa\end{pmatrix},\qquad H=2\kappa.
\end{equation}
These matrices use the parallel frame consisting of the horizontal
profile tangent followed by the positive vertical tangent, and the
horizontal unit normal $\widetilde\nu$ lifting the profile normal.
Here $A$ and $A^0$ are the scalar second fundamental forms
relative to $\widetilde\nu$, and $H=2\kappa$ is the signed trace of $A$;
the mean-curvature vector is $2\kappa\widetilde\nu$.
The induced metric is flat; the sign of the off-diagonal entries fixes
the normal convention consistently with the Hopf calculation.
For the standard pullback torus, choose one real area lift $\mathcal A_t$ and put
$\tau_t=(\mathcal A_t+iL_t)/(4\pi)\in\HH$. We use the complete metric
$h=(dx^2+dy^2)/y^2$ on $\HH$. The lattice and holonomy conventions
are justified in Section~\ref{sec:conformal}. A fixed finite torus
cover corresponds to its fixed finite-index sublattice.

\subsection{Curve geometry and the reduced equation}

All curves are parametrized closed immersions of $\R/\miwfZ$ into the
unit sphere $\Sph^2\subset\R^3$; multiplicity is retained. Let $s$ denote arclength,
$T=\partial_s\gamma$, $\nu=\gamma\times T$, and $\kappa$ the signed geodesic
curvature. In the frame equations and their differentiated identities below,
derivatives of vector-valued maps and fields are ordinary derivatives in
the ambient space $\R^3$. Covariant differentiation on the unit sphere
is denoted by $\nabla_s$. The frame equations are
\begin{equation}\label{miwf:frame}
 \gamma_s=T,\qquad T_s=-\gamma+\kappa\nu,
 \qquad \nu_s=-\kappa T.
\end{equation}
In particular, $\nabla_sT=\kappa\nu$, whereas the ordinary derivative
$T_s$ also contains the radial term $-\gamma$.
We use the normalization of \cite[Proposition 4.2 and (91)]{Jakob}:
\begin{equation}\label{miwf:flow}
 \miwfE(\gamma)=\int_\gamma(1+\kappa^2)\dd s,\qquad
 G=G_\gamma=2\kappa_{ss}+\kappa^3+\kappa,\qquad
 (\partial_t\gamma)^\perp=-\frac{G}{(1+\kappa^2)^2}\nu.
\end{equation}
Here $\perp$ denotes projection onto the line spanned by $\nu$ in the
tangent plane of the sphere. A reduced trajectory may be represented
with a tangential velocity. Its normal gauge is the equation
$\partial_t\gamma=-G_\gamma\nu/(1+\kappa^2)^2$ on the fixed
parameter circle. Explicitly, if
$\partial_t\gamma=f\nu+\zeta T$, solve
\[
 \partial_t\rho_t(x)
 =-\frac{\zeta(t,\rho_t(x))}
         {|\partial_x\gamma(t,\rho_t(x))|},\qquad \rho_0(x)=x.
\]
On every compact smooth time interval this equation gives circle
diffeomorphisms, and $\gamma_t\circ\rho_t$ has velocity
$(f\nu)\circ\rho_t$.
Throughout, $T_{\max}(\gamma_0)$ denotes maximal smooth continuation
of the normal-gauge solution. Statements about energy, length,
curvature, and convergence modulo reparametrization are independent
of the chosen tangential representative.

For a variation with velocity $f\nu$, the variation formulas
\eqref{tp:eq:curve-variation} and periodic integration by parts give
\[
 \delta\miwfE
 =\int\bigl(2\kappa f_{ss}+2\kappa(1+\kappa^2)f
             -\kappa(1+\kappa^2)f\bigr)\dd s
 =\int Gf\dd s.
\]
Tangential reparametrizations do not affect the energy identity
\begin{equation}\label{miwf:dissipation}
 -\frac{\mathrm d}{\mathrm dt}\miwfE(\gamma_t)
 =\miwfD(\gamma_t),\qquad
 \miwfD(\gamma)=\int_\gamma\frac{G^2}{(1+\kappa^2)^2}\dd s.
\end{equation}
An elastic curve is a regular closed solution of $G=0$.

In the normal gauge, the scalar variation identities are
\begin{equation}\label{tp:eq:curve-variation}
 \partial_t\dd s=-\kappa f\dd s,\qquad
 \partial_t\kappa=f_{ss}+(1+\kappa^2)f,
\end{equation}
and consequently
\begin{equation}\label{tp:eq:curvature-PDE}
 \partial_t\kappa
 =-\partial_s^2\left(\frac{2\kappa_{ss}+\kappa^3+\kappa}
                             {(1+\kappa^2)^2}\right)
  -\frac{2\kappa_{ss}+\kappa^3+\kappa}{1+\kappa^2}.
\end{equation}
To verify the curvature formula, covariantly differentiate the unit
tangent along a normal variation. Its normal time derivative is
$f_s\nu$, the arclength commutator is
$[\partial_t,\partial_s]=\kappa f\partial_s$, and the unit-sphere
curvature term is $f\nu$. Differentiating $\nabla_sT=\kappa\nu$
and taking its normal component therefore gives
$\kappa_t=f_{ss}+\kappa^2f+f$. For the covariant formulation of
curvature evolution in space forms, compare \cite{DallAcquaSphere}.
\begin{lemma}\label{miwf:length}
If $\miwfE(\gamma)\leq E_0$ with $0<E_0<\infty$, then its length $L$ satisfies
\begin{equation}\label{miwf:lengthbound}
 \frac{4\pi^2}{E_0}\leq L\leq E_0.
\end{equation}
\end{lemma}
\begin{proof}
The upper bound follows from the definition of $\miwfE$.
In arclength coordinates on $\R/L\miwfZ$, closedness gives
$\int_0^L T\dd s=\gamma(L)-\gamma(0)=0$ and $|T|=1$.
By \eqref{miwf:frame},
$\int_0^L|T_s|^2\dd s=\miwfE(\gamma)$.
The periodic Wirtinger inequality consequently gives
\[
 \miwfE(\gamma)=\int_0^L|T_s|^2\dd s
 \geq\frac{4\pi^2}{L^2}\int_0^L|T|^2\dd s
 =\frac{4\pi^2}{L}.
\]
The same inequality applies to periodic $W^{1,2}$ tangents, including
those of multiply covered curves.
\end{proof}

\subsection{Reconstruction on a fixed parameter torus}
\begin{lemma}[Normal Hopf reconstruction]\label{tp:lem:Hopf-reconstruction}
Let $\gamma_t$, $0\leq t<T\leq\infty$, be a smooth reduced
MIWF of regular closed profiles on a fixed parameter circle.
There is a smooth normal Hopf MIWF $F_t$ on one fixed parameter
torus representing their pullback immersions. Each fixed finite
covering of that torus gives a normal MIWF as well. For an additional
covering of degree $d$, with $D_F(t):=-W'(t)$,
\[
 \Will(F_t)=d\pi\miwfE(\gamma_t),\qquad
 D_F(t)=d\pi\miwfD(\gamma_t).
\]
\end{lemma}
\begin{proof}
Use a parameter circle of period $2\pi$. Lift $t\mapsto\gamma_t(0)$
horizontally to choose a smooth $z_t(0)\in\Sph^3$. For each $t$,
horizontally lift $\theta\mapsto\gamma_t(\theta)$ from that point,
obtaining $z_t(\theta)$ smoothly in both variables. Periodicity and
equivariance of the connection give
$z_t(\theta+2\pi)=q(t)z_t(\theta)$ for a smooth phase $q(t)\in\Sph^1$.
Choose a smooth real lift $\alpha(t)$ with $q(t)=e^{i\alpha(t)}$;
it exists on the time interval, for instance by integrating
$\alpha'=\operatorname{Im}(q'/q)$ after choosing $\alpha(0)$.
Then
\[
 \sigma_t(\theta)=e^{-i\alpha(t)\theta/(2\pi)}z_t(\theta)
\]
is a smooth periodic section of the pulled-back Hopf bundle.
Consequently $(\theta,\varphi)\mapsto e^{i\varphi}\sigma_t(\theta)$
is a smooth family of immersions on one fixed torus. All lifting
ODEs exist on every compact time subinterval: their coefficients
are smooth and the Hopf fibres are compact. Denote this family by
$\widehat F_t$.

We now verify its normal velocity without applying a global reduction
theorem to a self-intersecting curve. Fix a time and a parameter value.
After restricting to a sufficiently short interval, the corresponding
profile arc is embedded and contained in a coordinate disc. A smaller
arc extends to a smooth embedded closed curve in that disc: in local
coordinates take a smooth return arc on one side of the original arc,
joining its endpoints with matching jets. The Hopf preimage of this
closed curve is a simply parametrized smooth Hopf torus. Its immersion
and the pullback immersion under consideration agree locally, up to a
smooth change of surface coordinates and a fibre rotation. In
particular their local geometric jets agree. Thus the local formulas
of \cite[Proposition~4.3]{Jakob}, whose simple-torus hypotheses hold
for the extended curve, apply to the original parameter sheet. This
argument is pointwise in the parameter domain, so it also applies
where distinct sheets have the same image or the profile repeats.

To compute the gradient in the normalization of
Section~\ref{sec:setting}, let $e_1$ be the
horizontal unit lift of $T$, $e_2$ the unit vertical vector, and
$N$ the horizontal unit lift of $\nu$. Since $D\pi(e_1)=2T$,
$D\pi(N)=2\nu$, and $D\pi(e_2)=0$, the local shape operator is
\eqref{tp:eq:Hopf-geometry}. Horizontal lifting over arclength and
fibre rotation provide local coordinates $(s,\varphi)$ in which
$e_1=2\partial_s$ and $e_2=\partial_\varphi$; the induced metric is
$\tfrac14\dd s^2+\dd\varphi^2$. In codimension one the normal
connection of $N$ vanishes, and $\kappa$ is constant along fibres.
Consequently
\[
 \Delta^{\perp,\Sph^3}H=8\kappa_{ss}N,\qquad
 Q(A^0)(H)=4\kappa(1+\kappa^2)N.
\]
The gradient formula in Section~\ref{sec:setting} therefore gives
\begin{equation}\label{eq:sheetwise-hopf-gradient}
 \nabla_{L^2}\Will(\widehat F_t)
 =2\bigl(2\kappa_{ss}+\kappa^3+\kappa\bigr)N=2GN,
 \qquad |A^0|^2=2(1+\kappa^2).
\end{equation}
These are pointwise identities on the parameter torus; no identification
of distinct sheets is involved.

The equality $\pi\circ\widehat F_t=\gamma_t\circ p$, where $p$ is
the fixed base projection, implies after differentiation that
$D\pi(\partial_t\widehat F_t)=\partial_t\gamma_t\circ p$.
If $f=\langle\partial_t\gamma_t,\nu\rangle$, its normal surface
component is therefore $(f/2)N$, even if the profile has tangential
velocity. By \eqref{miwf:flow} and
\eqref{eq:sheetwise-hopf-gradient},
\[
 (\partial_t\widehat F_t)^\perp
 =-\frac{G}{2(1+\kappa^2)^2}N
 =-|A^0|^{-4}\nabla_{L^2}\Will(\widehat F_t).
\]
The remaining tangential velocity is a smooth vector field
$D\widehat F_t Z_t$ on the parameter torus. Solving
$\partial_t\phi_t=-Z_t\circ\phi_t$, $\phi_0=\mathrm{id}$,
removes this velocity. On every compact time subinterval $Z_t$ is
smooth on a compact torus, so its ODE gives smooth diffeomorphisms
throughout that subinterval. Uniqueness makes these solutions
compatible. Thus $F_t:=\widehat F_t\circ\phi_t$ is defined on all
of $[0,T)$ and solves the normal MIWF on the fixed torus.
A fixed torus covering preserves the equation because the operator
is local and normal.

The profile normal speed $f$ is twice the surface normal speed.
Using $a=2(1+\kappa^2)$ and fibre integration gives
\[
 D_F=d\pi\int(1+\kappa^2)^2 f^2\dd s
     =d\pi\int\frac{G_\gamma^2}{(1+\kappa^2)^2}\dd s.
\]
The energy identity follows from \eqref{eq:hopf-cover-energy}.
\end{proof}

\begin{corollary}[Finite-time Hopf torus trajectories]
\label{cor:hopf}
Suppose that the MIWF starts from a smooth parametrization \(F_0\) of
a smooth Hopf torus in the sense of Pinkall \cite{Pinkall}, that
\(W(0)\leq8\pi\), and that
\(T=T_{\max}(F_0)<\infty\).  Then there exists a unique unrenormalized
limiting integral varifold \(\Var_T\).  Its support is an embedded
\(C^1\) Hopf torus, and the convergence is both weak varifold
convergence and Hausdorff convergence.  Moreover,
\eqref{eq:conformal-limit}--\eqref{eq:limit-measure} hold
for a uniformly conformal bi-Lipschitz parametrization of the limiting
Hopf torus.  In particular, among the alternatives of
Corollary~\ref{cor:four-way}, every such finite-time Hopf trajectory
lies in alternative \textup{(D)}.
\end{corollary}

\begin{proof}
First identify the Hopf action on the fixed domain. If the initial
parametrization has an additional covering degree $d$, then
$W(0)=d\pi\miwfE(\gamma_0)$ by \eqref{eq:hopf-cover-energy}.
Lemma~\ref{miwf:length} gives $\miwfE(\gamma_0)\geq 2\pi$,
so $W(0)\leq8\pi$ forces $d=1$. In the standard Hopf
parametrization, fibre rotation induces a free circle action
$r_\theta$ on the parameter torus with
$F_0\circ r_\theta=e^{i\theta}F_0$. The two normal MIWF solutions
$F_t\circ r_\theta$ and $e^{i\theta}F_t$ agree initially.
Uniqueness gives equality throughout $[0,T)$, by restarting at
any common smooth slice. Hence $\pi\circ F_t$ descends to a smooth
curve $\gamma_t$ on the fixed quotient circle. Its differential is
nonzero because $DF_t$ has rank two and sends the fibre direction
to the Hopf direction. The sheetwise identity
\eqref{eq:sheetwise-hopf-gradient} gives the reduced equation
\eqref{miwf:flow} for this profile.
For a Hopf torus with profile curve \(\gamma_t\),
Proposition~4.3 of \cite{Jakob} gives
\begin{equation}
 \label{eq:hopf-nonumbilic}
 |A_t^0|^2=2\bigl(1+\kappa_{\gamma_t}^2\bigr)\geq2 .
\end{equation}
Here \(\kappa_{\gamma_t}\) is the signed geodesic curvature of
\(\gamma_t\subset\Sph^2\).  Thus \eqref{eq:un} holds with \(c_0=2\).
Theorem~\ref{thm:main} yields uniqueness of the varifold limit.
We give a direct argument for the regularity and topology of its support.
The proof of Corollary~\ref{cor:classification} gives embeddedness of
$F_t$ for $t>0$.  Choose the reduced profiles to traverse their
projected embedded curves once.  The Hopf area and energy identities give
\[
 \Area(F_t)=\pi L(\gamma_t),\qquad
 E_t:=\int_{\gamma_t}(1+\kappa_{\gamma_t}^2)\dd s
       =\Will(F_t)/\pi\leq8.
\]
The Fourier estimate in Lemma~\ref{miwf:length} gives
$4\pi^2/8\leq L(\gamma_t)\leq8$.  In particular the limiting mass is
positive, so $\Var_T\neq0$.  Corollaries~\ref{cor:classification}
and~\ref{cor:hausdorff} provide an embedded closed Lipschitz surface
$M:=\spt\Var_T$ and full Hausdorff convergence to $M$.
Each $F_t(\Sigma)$ is invariant under the Hopf fibre action.  Passing
this invariance to the Hausdorff limit gives
$M=\pi^{-1}(\pi(M))$.

For any $t_j\nearrow T$, reparametrize the profiles with constant
speed on $\R/\mathbb Z$.  We use explicitly the second identity in
\eqref{miwf:frame}, namely $T_s=-\gamma+\kappa\nu$.  Since
$\partial_x\gamma=L T$ in this parametrization, it gives
\[
 \partial_x^2\gamma=L^2(-\gamma+\kappa\nu),\qquad
 \|\partial_x^2\gamma_{t_j}\|_{L^2(\R/\mathbb Z)}^2
 =L(\gamma_{t_j})^3 E_{t_j}\leq8^4.
\]
Together with the length bounds this gives a subsequence converging
in $C^1$ to a regular closed curve $\beta$ with constant positive
speed.  Hausdorff convergence and continuity of $\pi$ then give
$B:=\beta(\Sph^1)=\pi(M)$, and hence $M=\pi^{-1}(B)$.
We now show that $B$ is embedded, using the manifold property of $M$.
Fix $z\in B$ and a parameter at which $\beta$ equals $z$.  On a small
open parameter interval its image is an embedded $C^1$ arc $J$
through $z$.  A Hopf bundle chart over this arc gives an embedded
$C^1$ surface patch over $J$ around each point of the fibre
$\pi^{-1}(z)$.  Each patch lies in $M$ and is open in $M$ by
invariance of domain: its inclusion is an injective map between
two-dimensional topological manifolds, to which
\cite[Theorem~2B.3]{Hatcher} applies in coordinate charts.
The union $O$ of finitely many such patches is open in $M$, contains
the entire compact fibre $\pi^{-1}(z)$, and satisfies $\pi(O)\subset J$.
The compact set $\pi(M\setminus O)$ does not contain $z$.
Consequently there is a neighborhood $U$ of $z$ in $\Sph^2$ with
$B\cap U\subset J$; shrinking $U$ inside an arc chart gives
$B\cap U=J\cap U$.  Thus $B$ is locally an embedded regular $C^1$
arc at every point.  It is compact and connected and has no boundary,
so it is an embedded $C^1$ circle.
Its Hopf preimage is an embedded $C^1$ torus: horizontal lifting of
one traversal and a fibre rotation closing its terminal phase give
a $C^1$ product parametrization by $\Sph^1\times\Sph^1$.
This proves the asserted topology and regularity directly.
Corollary~\ref{cor:parametrization} gives the remaining claims.
\end{proof}

We next fix the profile gauge.  Identify \(\Sph^1\) with
\(\R/\mathbb Z\).  Use the smooth oriented profile on the quotient circle in
Corollary~\ref{cor:hopf}, in the normal gauge of \eqref{miwf:flow}.
Fix a material point \(0\in\Sph^1\) and set
\begin{equation}
 \label{eq:elastic-energy}
 \miwfE(\gamma):=\int_{\Sph^1}(1+\kappa_\gamma^2)\dd s_\gamma,
 \qquad L(\gamma):=\int_{\Sph^1}\dd s_\gamma .
\end{equation}
Thus \(\partial_t\gamma_t\) is normal to the curve in \(\Sph^2\).
For \(x\in[0,1]\), extended equivariantly to \(\R\), define
\begin{equation}
 \label{eq:anchored-gauge}
 \phi_t(x):=\frac1{L(\gamma_t)}
       \int_0^x|\partial_z\gamma_t(z)|\dd z,\qquad
 \psi_t:=\phi_t^{-1},\qquad
 \widetilde\gamma_t:=\gamma_t\circ\psi_t .
\end{equation}
Then \(\phi_t(0)=\psi_t(0)=0\) and
\begin{equation}
 \label{eq:constant-profile-speed}
 |\partial_x\widetilde\gamma_t|=L(\gamma_t)
 \quad\text{on }\Sph^1 .
\end{equation}
Since \(\widetilde\gamma_t\) is an orientation-preserving
reparametrization of \(\gamma_t\), one also has
\[
 L(\widetilde\gamma_t)=L(\gamma_t),
 \qquad \miwfE(\widetilde\gamma_t)=\miwfE(\gamma_t).
\]
The material anchor \(0\) fixes the phase of this
constant-speed reparametrization; its spatial image \(\gamma_t(0)\)
evolves with the normal-gauge profile.

\begin{lemma}[Anchored constant-speed estimate for normal families]
\label{lem:anchored-gauge}
Let \(I\subset\R\) be an interval and let
\(f\colon I\times\Sph^1\to\R^d\) be a smooth family of immersed
closed curves whose velocity is everywhere normal, namely
\begin{equation}
 \label{eq:normal-curve-family}
 \langle\partial_tf,\partial_xf\rangle=0.
\end{equation}
Let \(\widetilde f\) be the anchored constant-speed parametrization
defined as in \eqref{eq:anchored-gauge}.  Write
\[
 L(f)=\int_{\Sph^1}\dd s_f,\qquad
 E_{\R^d}(f)=\frac12\int_{\Sph^1}|\vec\kappa_f|^2\dd s_f .
\]
Here $\vec\kappa_f:=\partial_{s_f}^2f$ is the Euclidean curvature
vector, with $s_f$ denoting arclength along $f$.
Then, at every time,
\begin{equation}
 \label{eq:closed-gauge-estimate}
 \|\partial_t\widetilde f\|_{L^2(\dd x)}
 \leq
 \left(\frac2{L(f)}+16E_{\R^d}(f)\right)^{1/2}
 \|\partial_tf\|_{L^2(\dd s_f)} .
\end{equation}
\end{lemma}

\begin{proof}
Fix a time and abbreviate \(V:=\partial_tf\), \(L:=L(f)\), and
\(q:=\langle V,\vec\kappa_f\rangle\).  If \(s=s(x)\) denotes
arc length from the anchored point, the normality condition
\eqref{eq:normal-curve-family} gives
\[
 \partial_t\dd s_f=-q\,\dd s_f.
\]
Consequently, differentiation of
\[
 \phi(x)=L^{-1}\int_0^x|\partial_zf|\dd z,\qquad
 \psi=\phi^{-1},
\]
and the first variation formula for the length of a closed curve give
\[
 L'=-\int_0^Lq\,\dd s
\]
and
\[
 |(\partial_xf\circ\psi)\,\partial_t\psi|\circ\phi
 =
 \left|
   -\frac{s}{L}\int_0^Lq\,\dd s+\int_0^sq\,\dd r
 \right| .
\]
Since \(\dd(\phi(x))=\dd s/L\), a change of variables yields
\[
 \int_{\Sph^1}|V\circ\psi|^2\dd x
 =\frac1L\|V\|_{L^2(\dd s_f)}^2.
\]
Moreover,
\begin{align*}
 \int_{\Sph^1}|(\partial_xf\circ\psi)\partial_t\psi|^2\dd x
 &\leq
 4\left(\int_{\Sph^1}|q|\dd s_f\right)^2\\
 &\leq
 4\|V\|_{L^2(\dd s_f)}^2
   \|\vec\kappa_f\|_{L^2(\dd s_f)}^2\\
 &=8E_{\R^d}(f)\|V\|_{L^2(\dd s_f)}^2 .
\end{align*}
Finally,
\(\partial_t\widetilde f=V\circ\psi+
 (\partial_xf\circ\psi)\partial_t\psi\).
The inequality \(|X+Y|^2\leq2|X|^2+2|Y|^2\) proves
\eqref{eq:closed-gauge-estimate}.  This is the closed-curve,
normal-velocity version of
\cite[Lemma~4.10 and Remark~4.11]{RuppSpener}.
\end{proof}

\begin{lemma}[Profile and surface speeds]
\label{lem:profile-surface-speed}
Let $F_t$ be a smooth MIWF of Hopf immersions and let $\gamma_t$ be its
corresponding normal-gauge profile evolution, retaining profile
traversals.  If the additional torus covering degree in
\eqref{eq:hopf-cover-energy} is $d$, then
\begin{equation}
 \label{eq:profile-surface-speed}
 \|\partial_t\gamma_t\|_{L^2(\dd s_{\gamma_t})}
 =\frac2{\sqrt{d\pi}}
 \|\partial_tF_t\|_{L^2(\Sigma,\dd\mu_t)} .
\end{equation}
No embeddedness hypothesis is required.  In the embedded,
once-traversed case one has $d=1$.
\end{lemma}

\begin{proof}
The unit normal $N$ of the Hopf immersion is horizontal. The Hopf
map has horizontal dilation
\begin{equation}
 \label{eq:hopf-horizontal-dilation}
 |D\pi(v)|=2|v| \qquad\text{for every horizontal vector }v.
\end{equation}
By \eqref{eq:sheetwise-hopf-gradient},
$\nabla_{L^2}\Will(F_t)=2G_{\gamma_t}N$ and
$|A_t^0|^2=2(1+\kappa_{\gamma_t}^2)$ on each sheet. Thus in the
normal gauges $|\partial_tF_t|=\tfrac12|\partial_t\gamma_t|$.
The fibre length $2\pi$ and \eqref{eq:hopf-horizontal-dilation} give
\[
 \int_\Sigma h\,\dd\mu_t=d\pi\int_{\Sph^1}h\,\dd s_{\gamma_t}
\]
for each function $h$ pulled back from the profile. Taking
$h=|\partial_t\gamma_t|^2/4$ proves \eqref{eq:profile-surface-speed}.
The arclength integral counts every profile traversal, and $d$ counts
the additional torus sheets.
\end{proof}

\subsection{Anchored finite-time endpoints at arbitrary energy}

The energy restriction in the embedded-endpoint theorem is not needed
for convergence of the parametrized profiles.  The quantities $\miwfE$, $G$, and $\miwfD$ use the normalization
\eqref{miwf:flow}--\eqref{miwf:dissipation}.

\begin{theorem}[All-energy finite-time profile endpoint]
\label{thm:all-energy-profile}
Let $\gamma_t$, $0\leq t<T<\infty$, be a smooth normal solution of
the reduced MIWF with a smooth regular closed initial curve.  Retain
parametrized multiplicity and use one fixed material anchor in
\eqref{eq:anchored-gauge}.  Put
$E_t=\miwfE(\gamma_t)$, $E_0=\miwfE(\gamma_0)$, and
$E_T=\lim_{t\uparrow T}E_t$.  With
\[
 C(E_0)=\left(\frac{E_0}{2\pi^2}+8E_0\right)^{1/2},
\]
one has, for $0\leq s<t<T$,
\begin{equation}\label{eq:all-energy-profile-tail}
 \|\widetilde\gamma_t-\widetilde\gamma_s\|_{L^2(\dd x)}
 \leq C(E_0)\sqrt{(t-s)(E_s-E_t)}.
\end{equation}
There is a unique regular constant-speed endpoint
$\gamma_T\in W^{2,2}(\R/\mathbb Z,\Sph^2)$ such that
\[
 \widetilde\gamma_t\rightharpoonup\gamma_T\text{ in }W^{2,2},
 \qquad
 \widetilde\gamma_t\longrightarrow\gamma_T
       \text{ in }C^{1,\alpha}\quad(0<\alpha<1/2),
\]
and
\begin{equation}\label{eq:all-energy-endpoint-bounds}
 \frac{4\pi^2}{E_0}\leq L_T=|\partial_x\gamma_T|\leq E_0,
 \qquad \miwfE(\gamma_T)\leq E_T.
\end{equation}
The multiplicity-counting Hopf varifolds induced by these profiles
converge to the varifold of $X_{\gamma_T}$.  For every fixed additional
torus covering of degree $d$, the limit is $d\Var(X_{\gamma_T})$,
with total mass at least $4d\pi^3/E_0$.  Neither embeddedness nor a
low-energy hypothesis is required.
\end{theorem}
\begin{proof}
Lemma~\ref{miwf:length} gives
$4\pi^2/E_0\leq L_t\leq E_t\leq E_0$.  For a spherical curve its
Euclidean bending energy in Lemma~\ref{lem:anchored-gauge} is $E_t/2$.
The reduced speed $f=-G/(1+\kappa^2)^2$ satisfies
\[
 \|\partial_t\gamma_t\|_{L^2(\dd s_t)}^2
 =\int\frac{G^2}{(1+\kappa^2)^4}\dd s_t\leq\miwfD(\gamma_t).
\]
The anchored estimate \eqref{eq:closed-gauge-estimate}, with
$2/L_t+8E_t\leq E_0/(2\pi^2)+8E_0=C(E_0)^2$, gives
$\|\partial_t\widetilde\gamma_t\|_{L^2(\dd x)}
 \leq C(E_0)\sqrt{\miwfD(\gamma_t)}$.
Integration and time Cauchy--Schwarz, using $E'=-\miwfD$, prove
\eqref{eq:all-energy-profile-tail}.  Finite time makes the profiles
Cauchy in $L^2$.

In constant-speed coordinates, the same identity
$T_s=-\gamma+\kappa\nu$ from \eqref{miwf:frame} gives
\[
 \|\partial_x^2\widetilde\gamma_t\|_2^2=L_t^3E_t\leq E_0^4.
\]
Together with $|\widetilde\gamma_t|=1$ and
$|\partial_x\widetilde\gamma_t|=L_t\leq E_0$, this is a uniform
$W^{2,2}$ bound.  Every time sequence admits a subsequence converging
weakly in $W^{2,2}$ and strongly in $C^{1,\alpha}$ for each
$\alpha<1/2$.  Its $L^2$ limit is the already fixed $\gamma_T$.
A contradiction argument for each of these topologies gives
convergence of the full family.  A diagonal subsequence over
$\alpha\uparrow1/2$ verifies the assertions simultaneously for all
smaller exponents.  Uniform convergence of the derivatives gives
$L_t\to L_T$ and $|\partial_x\gamma_T|=L_T>0$.
Weak lower semicontinuity in the last display and $L_t\to L_T$
give $\miwfE(\gamma_T)\leq E_T$.

We justify the varifold assertion on the parameter bundle, so that
repeated traversals are retained.  Choose a fixed smooth local section
of the Hopf bundle over a neighborhood of $\gamma_T(0)$, and evaluate
it at $\widetilde\gamma_t(0)$ for all sufficiently late times.  This
gives initial lifts $z_t(0)\to z_T(0)$ in $\Sph^3$.
Horizontally lift the parameter paths from these points to obtain
$z_t(x)$, $0\leq x\leq1$.
Their coefficients converge uniformly by the $C^1$ profile convergence.
Continuous dependence in this compact-fibre ODE, followed by its
defining equation for $z_t'$, gives $z_t\to z_T$ in $C^1$.
Write $z_t(1)=q_tz_t(0)$ with $q_t\in\Sph^1$.  Since $q_t\to q_T$,
choose compatible real phases $q_t=e^{i\alpha_t}$ with
$\alpha_t\to\alpha_T$.  Equivariance and periodicity give
$z_t(x+1)=q_tz_t(x)$ and the same identity for their derivatives.
Consequently
\[
 \sigma_t(x)=e^{-i\alpha_t x}z_t(x),\qquad
 Y_t(x,\vartheta)=e^{i\vartheta}\sigma_t(x)
\]
are periodic, $Y_t$ represents $X_{\gamma_t}$ on one fixed torus,
and $Y_t\to Y_T$ in $C^1$.  The limit is an immersion because its
profile speed is positive and its fibre derivative is a unit vertical
vector.  Its induced area densities and tangent planes converge
uniformly.  Testing the varifold definition therefore gives
$\Var(Y_t)\rightharpoonup\Var(Y_T)$.
On any fixed finite torus cover the area formula multiplies this
varifold by its degree $d$; the identity of the limit does not depend
on the chosen bundle trivialization.  Finally
$\Area(Y_t\circ q)=d\pi L_t$ proves the mass bound.
\end{proof}

\subsection{The embedded low-energy endpoint}

\begin{theorem}[Full convergence of a fixed anchored finite-time Hopf profile]
\label{thm:anchored-profile}
Under the hypotheses of Corollary~\ref{cor:hopf}, fix the oriented
normal gauge and the material anchor introduced above, and let
\(\widetilde\gamma_t\) be the anchored constant-speed profile
\eqref{eq:anchored-gauge}.  Set
\begin{equation}
 \label{eq:gauge-constant}
 C_{\mathrm g}:=\left(\frac4{\pi^2}+64\right)^{1/2}.
\end{equation}
Then, for \(0\leq s<t<T\),
\begin{equation}
 \label{eq:anchored-profile-cauchy}
 \|\widetilde\gamma_t-\widetilde\gamma_s\|_{L^2(\Sph^1)}
 \leq
 \frac{C_{\mathrm g}}{\sqrt\pi}
 \sqrt{(t-s)\bigl(W(s)-W(t)\bigr)} .
\end{equation}
There are an embedded regular constant-speed curve
\(\gamma_T\in W^{2,2}(\Sph^1,\Sph^2)\cap C^1(\Sph^1,\Sph^2)\)
and a number \(L_T\in[\pi^2/2,8]\) such that, as \(t\nearrow T\),
\begin{align}
 \widetilde\gamma_t&\rightharpoonup\gamma_T
 &&\text{weakly in }W^{2,2}(\Sph^1,\R^3),\label{eq:profile-weak}\\
 \widetilde\gamma_t&\longrightarrow\gamma_T
 &&\text{strongly in }W^{1,2}(\Sph^1,\R^3),\label{eq:profile-w12}\\
 \widetilde\gamma_t&\longrightarrow\gamma_T
 &&\text{strongly in }C^{1,\alpha}(\Sph^1,\R^3),
   \quad 0<\alpha<\tfrac12,\label{eq:profile-strong}\\
 L(\widetilde\gamma_t)&\longrightarrow L_T=L(\gamma_T),
 \qquad |\partial_x\gamma_T|=L_T .\label{eq:length-limit}
\end{align}
The trace of \(\gamma_T\) is the embedded \(C^1\) profile
\(\pi(\spt\Var_T)\) of the unique limiting Hopf torus.
\end{theorem}

\begin{proof}
The quotient profile in Corollary~\ref{cor:hopf} is a smooth normal
solution of the reduced flow and
$E_0=\miwfE(\gamma_0)=W(0)/\pi\leq8$. Apply
Theorem~\ref{thm:all-energy-profile} with the same material anchor.
Its constant satisfies
\[
 C(E_0)^2=\frac{E_0}{2\pi^2}+8E_0
 \leq\frac4{\pi^2}+64=C_{\mathrm g}^2.
\]
Since $W(t)=\pi\miwfE(\gamma_t)$,
\eqref{eq:all-energy-profile-tail} gives
\eqref{eq:anchored-profile-cauchy}. The same theorem supplies the
full weak $W^{2,2}$ and strong $C^{1,\alpha}$ convergence, hence also
strong $W^{1,2}$ convergence, and
\[
 L_T\in[\pi^2/2,8],\qquad
 \miwfE(\gamma_T)\leq E_T\leq8.
\]
It also identifies $\Var_T$ with the Hopf varifold of $\gamma_T$.
Thus $\gamma_T(\Sph^1)=\pi(\spt\Var_T)$, the embedded $C^1$
profile obtained in Corollary~\ref{cor:hopf}.

The regular map $\gamma_T$ is a finite covering of this circle, of
some positive degree $m$ after orienting the image by $\gamma_T$.
Because it has constant speed, it is the $m$-fold traversal of a
constant-speed once-traversed regular curve $\beta\in W^{2,2}\cap C^1$.
The periodic Wirtinger argument of Lemma~\ref{miwf:length} applies to
its $W^{1,2}$ tangent and gives $\miwfE(\beta)\geq2\pi$. Therefore
\[
 2\pi m\leq m\miwfE(\beta)=\miwfE(\gamma_T)\leq8<4\pi.
\]
Hence $m=1$, and $\gamma_T$ is an embedding.
\end{proof}

\begin{remark}[Endpoint energy and curvature]
\label{rem:no-strong-profile}
For the anchored constant-speed profiles,
$\|\partial_x^2\widetilde\gamma_t\|_2^2=L_t^3E_t$.
Thus weak $W^{2,2}$ convergence and $L_t\to L_T>0$ give strong
$W^{2,2}$ convergence exactly when $E_t\to\miwfE(\gamma_T)$.
Corollary~\ref{miwf:endpoint-W42} proves this energy identity if
dissipation is bounded along an endpoint sequence.
\end{remark}

\begin{remark}[Continuation from a profile endpoint]
Smooth continuation requires a regular umbilic-free initial
immersion in a Cauchy trace space. The first-order endpoint in
Theorem~\ref{thm:all-energy-profile} and the marked conformal
endpoint in Section~\ref{sec:conformal} are geometric limits.
Their relation to stronger trace regularity is addressed by
Corollary~\ref{miwf:endpoint-W42}.
\end{remark}

\begin{corollary}[A nonzero finite Hopf varifold endpoint]
\label{tp:cor:Hopf-varifold-endpoint}
For the reduced trajectory in Theorem~\ref{thm:all-energy-profile}
with $T<\infty$, its reconstructed Hopf MIWF and every fixed finite
covering have a unique nonzero integral-varifold endpoint, counting
multiplicity. In the normal surface gauge their area measures also
converge in total variation. For an additional covering of degree $d$,
the limiting total mass is at least $4d\pi^3/E_0$, where
$E_0=\miwfE(\gamma_0)$.
\end{corollary}
\begin{proof}
Lemma~\ref{tp:lem:Hopf-reconstruction} supplies the fixed normal
surface gauge and each additional covering on $[0,T)$.
Formula \eqref{tp:eq:Hopf-geometry} gives $a=|A^0|^2\geq2$.
Thus $\int_0^T m(t)^{-2}\dd t\leq T/4$, and Proposition
\ref{prop:weighted-endpoint} gives one integral-varifold endpoint and the
stated source-measure convergence. Finally
$\Area(F_t)=d\pi L_t\geq4d\pi^3/E_0$.
Weak convergence on the compact sphere preserves total mass, so
the same lower bound holds at the endpoint.
\end{proof}

\section{Dissipation compactness and global reduced trajectories}
\label{sec:miwf-compactness}

\subsection{Quantitative curvature estimates}

\begin{proposition}\label{miwf:coercivity}
For each $0<E_0<\infty$ there is a constant $C=C(E_0)$ such that every
smooth regular closed curve with $\miwfE\leq E_0$ satisfies
\begin{align}
 \miwfnorm{\kappa_s}_{L^2}&\leq C(1+\sqrt{\miwfD}),\label{miwf:first}\\
 \miwfnorm{\kappa}_{L^\infty}^2&\leq C(1+\sqrt{\miwfD}),\label{miwf:max}\\
 \miwfnorm{\kappa_{ss}}_{L^2}&\leq C(1+\miwfD).\label{miwf:second}
\end{align}
All norms in this proposition are with respect to arclength.
\end{proposition}
\begin{proof}
Write
\[
 K=\int\kappa^2\dd s,\quad X=\miwfnorm{\kappa_s}_2,
 \quad M=\miwfnorm\kappa_\infty,\quad
 f=\frac{2\kappa_{ss}+\kappa^3+\kappa}{1+\kappa^2}.
\]
Thus $K\leq E_0$ and $\|f\|_2=\sqrt{\miwfD}$. Since the average of
$\kappa^2$ is $K/L$, there is $s_0$ with $\kappa(s_0)^2\leq K/L$.
For any $s$, integrate $(\kappa^2)_s=2\kappa\kappa_s$ along an arc
from $s_0$ to $s$ and apply Cauchy--Schwarz on that arc, then enlarge
the integrals to the full circle. Taking the maximum over $s$ gives
\begin{equation}\label{miwf:elementary}
 M^2\leq\frac KL+2\sqrt K X.
\end{equation}
Multiply $2\kappa_{ss}+\kappa^3+\kappa=(1+\kappa^2)f$ by $\kappa$
and integrate over the periodic arclength circle. Periodicity gives
$\int 2\kappa\kappa_{ss}\,ds=-2\int\kappa_s^2\,ds=-2X^2$;
rearranging yields the exact identity
\[
 2X^2=\int\kappa^4\dd s+K-\int f\kappa(1+\kappa^2)\dd s.
\]
Cauchy--Schwarz and $|\kappa|\leq M$ give
$\int\kappa^4\,ds\leq KM^2$ and
\[
 \left|\int f\kappa(1+\kappa^2)\,ds\right|
 \leq\|f\|_2\|\kappa(1+\kappa^2)\|_2
 \leq\sqrt{\miwfD K}(1+M^2).
\]
Consequently,
\begin{align*}
 2X^2
 &\leq KM^2+K+\sqrt{\miwfD K}(1+M^2)\\
 &=(K+\sqrt{\miwfD K})(1+M^2)\\
 &\leq(K+\sqrt{\miwfD K})
       \left(1+\frac KL+2\sqrt K X\right).
\end{align*}
By \eqref{miwf:lengthbound}, $K/L\leq E_0^2/(4\pi^2)$.
Consequently the last expression is bounded by
$C(E_0)(1+\sqrt\miwfD)(1+X)$. To make the absorption explicit, put
$B=C(E_0)(1+\sqrt\miwfD)\geq1$. From $2X^2\leq B(1+X)$, after
enlarging $B$ if needed, Young's inequality gives
$BX\leq X^2+\frac14B^2$ and hence
$X^2\leq C B^2\leq C(E_0)(1+\miwfD)$. This proves
\eqref{miwf:first}. Inserting
$X\leq C(E_0)(1+\sqrt\miwfD)$ into \eqref{miwf:elementary} gives
$\|\kappa\|_\infty^2=M^2\leq C(E_0)(1+\sqrt\miwfD)$, which is
\eqref{miwf:max}. Finally, the defining equation for $f$ implies
\begin{align*}
 2\miwfnorm{\kappa_{ss}}_2
 &\leq (1+M^2)\sqrt\miwfD+\miwfnorm{\kappa^3+\kappa}_2\\
 &\leq(1+M^2)(\sqrt\miwfD+\sqrt K)
 \leq C(E_0)(1+\sqrt\miwfD)^2
 \leq C(E_0)(1+\miwfD).
\end{align*}
This proves \eqref{miwf:second}.
\end{proof}

\begin{proposition}[Compactness at bounded-dissipation times]
\label{miwf:bounded-dissipation-compactness}
Let $\gamma_j$ be smooth regular closed curves with
$\miwfE(\gamma_j)\leq E_0$ and $\miwfD(\gamma_j)\leq D_0$ for fixed
$E_0,D_0<\infty$.  After constant-speed reparametrization on
$\R/\miwfZ$ with a fixed material anchor, a subsequence converges weakly
in $W^{4,2}$ and strongly
in $C^{3,\alpha}$ for every $0<\alpha<\frac12$ to a regular closed
$W^{4,2}$ curve.  The limit need not be elastic.
\end{proposition}
\begin{proof}
Write $s_j$ for arclength along $\gamma_j$. The three estimates
\eqref{miwf:first}, \eqref{miwf:max}, and
\eqref{miwf:second}, proved in Proposition~\ref{miwf:coercivity},
give, respectively,
\begin{align*}
 \|\partial_{s_j}\kappa_j\|_{L^2(ds_j)}
 &\leq C(E_0)(1+\sqrt{D_0}),\\
 \|\kappa_j\|_{L^\infty}^2
 &\leq C(E_0)(1+\sqrt{D_0}),\\
 \|\partial_{s_j}^2\kappa_j\|_{L^2(ds_j)}
 &\leq C(E_0)(1+D_0).
\end{align*}
The energy bound also gives
$\|\kappa_j\|_{L^2(ds_j)}\leq E_0^{1/2}$. Write
the constant-speed representatives on $\R/\miwfZ$ as
$|\partial_x\gamma_j|=L_j$. The length bounds put $L_j$ in the
fixed compact interval $[4\pi^2/E_0,E_0]$. For each fixed $j$, choose
arclength $s_j=L_jx$ starting at the anchor. Since $L_j$ is independent
of $x$, arclength differentiation is simply
\[
 \partial_{s_j}=L_j^{-1}\partial_x.
\]
All these vector-valued derivatives are ordinary derivatives in $\R^3$.
The frame equations \eqref{miwf:frame} can be written as
\[
 \partial_x\gamma_j=L_jT_j,\qquad
 \partial_x^2\gamma_j=L_j^2(-\gamma_j+\kappa_j\nu_j),
 \quad \nu_j=\gamma_j\times T_j.
\]
Differentiating
$\partial_{s_j}^2\gamma_j=-\gamma_j+\kappa_j\nu_j$
once and twice with respect to arclength, and using
$\partial_{s_j}T_j=-\gamma_j+\kappa_j\nu_j$ and
$\partial_{s_j}\nu_j=-\kappa_jT_j$, yields
\begin{align}
 \partial_{s_j}^3\gamma_j
 &=-(1+\kappa_j^2)T_j+(\partial_{s_j}\kappa_j)\nu_j,
 \label{miwf:thirdder}\\
 \partial_{s_j}^4\gamma_j
 &=(1+\kappa_j^2)\gamma_j
 -3\kappa_j\partial_{s_j}\kappa_jT_j
 +(\partial_{s_j}^2\kappa_j-\kappa_j-\kappa_j^3)\nu_j.
 \label{miwf:fourthder}
\end{align}
Now combine these two formulas with the elementary constant-speed
identities
\[
 \partial_x^3\gamma_j=L_j^3\partial_{s_j}^3\gamma_j,\qquad
 \partial_x^4\gamma_j=L_j^4\partial_{s_j}^4\gamma_j.
\]
The bounds \eqref{miwf:first}, \eqref{miwf:max}, and
\eqref{miwf:second}, together with \eqref{miwf:lengthbound}, show
that both derivatives are uniformly bounded in $L^2(\R/\miwfZ)$.
Indeed, the middle estimate bounds $\kappa_j$ in $L^\infty$, while
the first and third bound its first and second arclength derivatives
in $L^2$. Every term in \eqref{miwf:thirdder}--\eqref{miwf:fourthder}
is therefore controlled, since $|\gamma_j|=|T_j|=|\nu_j|=1$.
The conversion of the norms uses $\dd s_j=L_j\,\dd x$ and gives
\[
 \|\partial_x^k\gamma_j\|_{L^2(dx)}
 =L_j^{k-1/2}\|\partial_{s_j}^k\gamma_j\|_{L^2(ds_j)}
 \quad (1\leq k\leq4).
\]
Moreover, the first and third curvature estimates bound
$\partial_{s_j}\kappa_j$ in $W^{1,2}$. The one-dimensional Morrey
embedding $W^{1,2}\hookrightarrow C^0$, with uniform constants on
these circles by \eqref{miwf:lengthbound}, also bounds
$\partial_{s_j}\kappa_j$ in $C^0$. Together with the middle
curvature estimate and the frame formulas, this shows that all
$\partial_x^k\gamma_j$ through order three are uniformly bounded
in $C^0(\R/\miwfZ)$.
Together with $|\gamma_j|=1$ and the length bounds, this gives
$\|\gamma_j\|_{W^{4,2}(\R/\miwfZ;\R^3)}\leq C(E_0,D_0)$.
Weak compactness and the compact embedding
$W^{4,2}\Subset C^{3,\alpha}$ for $0<\alpha<1/2$ give the
asserted convergence. A diagonal subsequence over exponents increasing
to $1/2$ gives one subsequence for every such $\alpha$. Passing also
to a subsequence with $L_j\to L_\infty>0$, strong $C^1$
convergence gives $|\partial_x\gamma_\infty|=L_\infty$, so the
limit is regular. Uniform convergence also preserves
$|\gamma_\infty|=1$. The anchor fixes the phase on the parameter circle.
For a constant sequence of latitudes of geodesic curvature $1$,
$G=2$ and $\miwfD=L<\infty$, but the limit is not elastic.
\end{proof}

\begin{corollary}[Finite-interval spacetime curvature bound]
\label{miwf:spacetime-curvature}
Let $\gamma_t$ be a smooth solution of \eqref{miwf:flow} on
$[0,T)$, where $T<\infty$, and suppose
$\miwfE(\gamma_t)\leq E_0<\infty$.  Put
$E_T=\lim_{t\uparrow T}\miwfE(\gamma_t)$.  There is a constant
$C(E_0)$ such that
\begin{equation}\label{miwf:spacetime-L6}
 \int_0^T\int_{\gamma_t}(1+\kappa_t^2)^3\,\dd s_t\,\dd t
 \leq C(E_0)\bigl(T+E_0-E_T\bigr).
\end{equation}
For a Hopf lift with additional torus covering degree $d$, profile
multiplicity being included in $\dd s_t$, one has exactly
\[
 \int_0^T\int_\Sigma |A_t^0|^6\dd\mu_t\dd t
 =8d\pi\int_0^T\int_{\gamma_t}(1+\kappa_t^2)^3\dd s_t\dd t.
\]
Thus \eqref{miwf:spacetime-L6} gives a finite spacetime $L^6$ bound
for every fixed finite covering.
\end{corollary}
\begin{proof}
Write $M_t=\|\kappa_t\|_{L^\infty}$. Squaring the middle curvature
estimate \eqref{miwf:max} gives
\[
 (1+M_t^2)^2\leq C(E_0)(1+\miwfD(\gamma_t)).
\]
Consequently,
\[
 \int_{\gamma_t}(1+\kappa_t^2)^3\,\dd s_t
 \leq (1+M_t^2)^2\miwfE(\gamma_t)
 \leq C(E_0)(1+\miwfD(\gamma_t)).
\]
Integrating and using
$\int_0^T\miwfD(\gamma_t)\,\dd t
=\miwfE(\gamma_0)-E_T\leq E_0-E_T$
proves \eqref{miwf:spacetime-L6}.  The Hopf-lift statement follows
from the displayed geometric identities.
\end{proof}

\subsection{Compactness at times of vanishing dissipation}

\begin{proposition}\label{miwf:ps}
Let $\gamma_j$ be smooth regular closed curves with
$\miwfE(\gamma_j)\leq E_0$ and $\miwfD(\gamma_j)\to0$. After
constant-speed reparametrization on $\R/\miwfZ$ and passage to a subsequence,
there is a smooth regular closed elastic curve $\gamma_\infty$ such that
\begin{equation}\label{miwf:strong}
 \gamma_j\longrightarrow\gamma_\infty\quad\text{strongly in }W^{4,2}
 \text{ and in }C^{3,\alpha}\quad(0<\alpha<\tfrac12).
\end{equation}
\end{proposition}
\begin{proof}
Reparametrize each curve at constant speed on the fixed circle
$\R/\miwfZ$, retaining the chosen material anchor, and write
$|\partial_x\gamma_j|=L_j$ and $s_j=L_jx$ for arclength from that anchor.
By \eqref{miwf:lengthbound}, after
passing to a subsequence,
\[
 L_j\longrightarrow L_\infty,\qquad
 0<4\pi^2/E_0\leq L_\infty\leq E_0.
\]
Since $\miwfD(\gamma_j)\to0$, this sequence of dissipations is bounded.
The three estimates \eqref{miwf:first}, \eqref{miwf:max}, and
\eqref{miwf:second} therefore bound, respectively,
$\partial_{s_j}\kappa_j$ in $L^2(ds_j)$, $\kappa_j$ in
$L^\infty$, and $\partial_{s_j}^2\kappa_j$ in $L^2(ds_j)$.
As proved in Proposition~\ref{miwf:bounded-dissipation-compactness},
these bounds, the length bounds \eqref{miwf:lengthbound}, and the
ordinary Frenet identities \eqref{miwf:thirdder}--\eqref{miwf:fourthder},
combined with the elementary conversions
$\partial_x^k=L_j^k\partial_{s_j}^k$ and $\dd s_j=L_j\dd x$,
give a uniform $W^{4,2}$ bound on the fixed parameter circle.
Applying that proposition and extracting a subsequence gives
\[
 \gamma_j\rightharpoonup\gamma_\infty\quad\hbox{in }W^{4,2},
 \qquad
 \gamma_j\longrightarrow\gamma_\infty\quad\hbox{in }C^{3,\alpha}
 \quad\text{for every }0<\alpha<\tfrac12.
\]
The same subsequence works for every such $\alpha$, as in that proposition.
The constant-speed identity passes to the limit, so
$|\partial_x\gamma_\infty|=L_\infty>0$ and the limit is regular.
Since $\gamma_j$ take values in $\Sph^2$, their uniform convergence
also gives $\gamma_\infty(\R/\miwfZ)\subset\Sph^2$.

We now use the vanishing dissipation to identify the limit as elastic
and prove strong $W^{4,2}$ convergence.
The curvature is given on this parametrization by
\[
 \kappa_j=L_j^{-2}
 \langle\partial_x^2\gamma_j,\nu_j\rangle,\qquad
 T_j=L_j^{-1}\partial_x\gamma_j,\qquad
 \nu_j=\gamma_j\times T_j.
\]
Thus the strong $C^3$ convergence and $L_j\to L_\infty$ imply
$\kappa_j\to\kappa_\infty$ in $C^1$. Let
$M_j=\|\kappa_j\|_\infty$. From the definition of the dissipation,
\begin{equation}\label{miwf:residual}
 \|G_j\|_{L^2(ds_j)}
 \leq(1+M_j^2)\sqrt{\miwfD(\gamma_j)}\longrightarrow0.
\end{equation}
Because $ds_j=L_j\,dx$ and the lengths are bounded below, this also
implies $G_j\to0$ in $L^2(dx)$. The identity
$G_j=2\partial_{s_j}^2\kappa_j+\kappa_j^3+\kappa_j$, together with
$\partial_{s_j}=L_j^{-1}\partial_x$, becomes
\[
 2L_j^{-2}\partial_x^2\kappa_j+\kappa_j^3+\kappa_j=G_j.
\]
It follows that
\[
 \partial_x^2\kappa_j
 =\frac{L_j^2}{2}\bigl(G_j-\kappa_j^3-\kappa_j\bigr)
 \longrightarrow
 -\frac{L_\infty^2}{2}
       (\kappa_\infty^3+\kappa_\infty)
 \quad\text{strongly in }L^2(dx).
\]
Since $\kappa_j\to\kappa_\infty$ in $C^1$, the distributional limit
of $\partial_x^2\kappa_j$ is $\partial_x^2\kappa_\infty$.
Identifying it with the strong $L^2$ limit above gives the elastic equation
\[
 2\partial_{s_\infty}^2\kappa_\infty
       +\kappa_\infty^3+\kappa_\infty=0,
 \qquad \partial_{s_\infty}=L_\infty^{-1}\partial_x.
\]
The right-hand side of the corresponding equation for
$\partial_x^2\kappa_\infty$ is a smooth function of
$\kappa_\infty$. The one-dimensional ODE therefore bootstraps
$\kappa_\infty$ to $C^\infty$. The Frenet system
\eqref{miwf:frame} then bootstraps $\gamma_\infty$ to $C^\infty$.

For strong $W^{4,2}$ convergence, use the fourth-derivative identity
\eqref{miwf:fourthder} in the fixed coordinate:
\[
 \partial_x^4\gamma_j=L_j^4\left[
 (1+\kappa_j^2)\gamma_j
 -3\kappa_j(\partial_{s_j}\kappa_j)T_j
 +(\partial_{s_j}^2\kappa_j-\kappa_j-\kappa_j^3)\nu_j
 \right].
\]
We already have strong convergence of $\gamma_j,T_j,\nu_j,\kappa_j$
and $\partial_{s_j}\kappa_j$ in $C^0$, and
$\partial_{s_j}^2\kappa_j=L_j^{-2}\partial_x^2\kappa_j$
converges strongly in $L^2$. Hence the right-hand side converges
strongly in $L^2$, so $\partial_x^4\gamma_j\to
\partial_x^4\gamma_\infty$ strongly in $L^2$. The derivatives through
order three already converge uniformly. Thus
$\gamma_j\to\gamma_\infty$ strongly in $W^{4,2}$, as claimed.
\end{proof}

\begin{corollary}\label{miwf:goodtimes}
Every global solution of \eqref{miwf:flow} admits times $t_j\to\infty$
and constant-speed reparametrizations converging as in \eqref{miwf:strong}
to a smooth regular elastic curve.
\end{corollary}
\begin{proof}
Energy monotonicity bounds $\miwfE(\gamma_t)$ by $E_0=\miwfE(\gamma_0)$ and gives
\[
 \int_0^\infty\miwfD(\gamma_t)\dd t\leq E_0.
\]
For each integer $j\geq1$, choose $t_j\in(j,j+1)$ with
$\miwfD(\gamma_{t_j})\leq2\int_j^{j+1}\miwfD(\gamma_t)\dd t$;
if the integral vanishes, any interior time works. The integrals tend to
zero. Apply Proposition \ref{miwf:ps} and retain a subsequence if necessary.
\end{proof}

\subsection{Uniform local regularization}

We use the linear parabolic isomorphisms of
\cite[Theorem 1.23 and Remark 1.31]{Amann}, applied on a fixed parameter circle.
The construction below is local near a smooth elastic curve and requires
no uniform parabolicity bound along the entire trajectory. It is
consistent with the maximal-regularity approach to the surface MIWF in
\cite{JakobRegularity}; the scalar application is justified directly below.

\begin{lemma}\label{miwf:smoothing}
Let $\eta$ be a fixed smooth regular closed elastic curve. Suppose smooth
regular curves $\eta_j$ converge to $\eta$ in $W^{4,2}$, after
orientation-preserving reparametrization. For each prescribed finite
$2\tau>0$, the normal-gauge solutions of \eqref{miwf:flow} starting at
$\eta_j$ exist on $[0,2\tau]$ for all sufficiently large $j$. After
smooth reparametrization, they converge to $\eta$ in every $C^k$ norm,
uniformly for times in $[\tau,2\tau]$.
\end{lemma}
\begin{proof}
For $L>0$ write $\mathcal C_L=\R/L\miwfZ$. Use the arclength
parameter of $\eta$ on the fixed circle $\mathcal C_{L_\eta}$. Set $T_*=\eta_s$, $\nu_*=\eta\times T_*$, and let
$\kappa_*$ be its signed geodesic curvature. Nearby parametrized curves
can be written, after a diffeomorphism close to the identity, as
\[
 \Gamma(u)=\cos(u)\eta+\sin(u)\nu_*.
\]
For an immersed curve this construction takes place over the parameter
circle. Define the normal-coordinate map
\[
 \Xi(s,z)=\cos z\,\eta(s)+\sin z\,\nu_*(s).
\]
At $(s,0)$ its two derivatives are $T_*(s)$ and $\nu_*(s)$,
an orthonormal basis of $T_{\eta(s)}\Sph^2$. On sufficiently small
parameter rectangles it therefore has smooth local inverses
$(P_\alpha,R_\alpha)=\Xi_\alpha^{-1}$. Their values are a base
parameter and a scalar normal coordinate, respectively. For each
$\eta_j(s)$ choose the branch near $(s,0)$ and set locally
\[
 \sigma_j(s)=P_\alpha(\eta_j(s)),\qquad
 r_j(s)=R_\alpha(\eta_j(s)).
\]
Uniqueness near the same parameter makes these coordinates agree on
overlaps, with base parameters understood modulo $L_\eta$.
The smooth inverse branches and the chain rule through order four
give $\sigma_j\to\mathrm{id}$ and $r_j\to0$ in $W^{4,2}$.
Since $\sigma_j'>0$ for large $j$, $\sigma_j$ is a circle
diffeomorphism. The graph height is the scalar function
\[
 h_j=r_j\circ\sigma_j^{-1},\qquad
 \eta_j=\Gamma(h_j)\circ\sigma_j.
\]
Composition and inversion are continuous in $W^{4,2}$ for these
diffeomorphisms: the differentiated formulas through order four use
only products of lower derivatives, the top derivatives in $L^2$,
and denominators bounded away from zero by the bound on $\sigma_j'$.
Thus $h_j\to0$ in $W^{4,2}$. In this construction $P_\alpha$ and
$R_\alpha$ extract coordinates from points on a chosen local branch,
whereas $h\mapsto\Gamma(h)$ constructs a parametrized curve from a
scalar height. This construction uses local branches over the parameter circle
and therefore accommodates self-intersections and multiple covers.
For the related surface construction, see the proof of
\cite[Theorem~1.3]{JakobGlobal}, which uses Skorzinski's normal-graph
theorem \cite{Skorzinski}; see also \cite{Jakob}.
The \L ojasiewicz--Simon approach to convergence of the classical
Willmore flow is developed in Chill, Fa\v{s}angov\'a, and
Sch\"atzle~\cite{ChillFasangovaSchaetzle}. Here the
chart is used to construct the scalar initial-value problem near an arbitrary
elastic immersion and to prove uniform positive-time regularization.
The explicit inverse branches and the chain rule give the
$W^{4,2}$ continuity required here for any smooth elastic reference
immersion.

Define
\[
 a=\cos u-\kappa_*\sin u,\quad c=\sin u+\kappa_*\cos u,\quad
 p=u_s,
\]
\[
 b=-\sin u\,\eta+\cos u\,\nu_*,\quad
 v=(a^2+p^2)^{1/2}.
\]
Then $\Gamma_s=aT_*+u_s b$, and the unit normal within $\Sph^2$ is
$N=(ab-u_sT_*)/v$. In particular the normal velocity of a graph is
$(a/v)u_t$. Direct differentiation using \eqref{miwf:frame} gives
\[
 \Gamma_{ss}=-v^2\Gamma+(-2cp-\kappa_*'\sin u)T_*
                         +(ac+u_{ss})b.
\]
Consequently its signed geodesic curvature is
\begin{equation}\label{miwf:graphcurvature}
 \kappa[u]=\frac{\langle\Gamma_{ss},N\rangle}{v^2}
 =\frac{a u_{ss}+c(a^2+2p^2)+\kappa_*'\sin u\,p}{v^3}.
\end{equation}
Arclength differentiation along the graph is $v^{-1}\partial_s$.
Substituting \eqref{miwf:graphcurvature} into \eqref{miwf:flow} gives
\begin{equation}\label{miwf:graphflow}
 u_t+a_4(s,u,u_s,u_{ss})u_{ssss}
 =Q(s,u,u_s,u_{ss},u_{sss}),\qquad
 a_4=\frac{2}{v^4(1+\kappa[u]^2)^2}.
\end{equation}
For the lower-order coefficients, write
\[
 C(s,y,p)=\frac{a(s,y)}{v(s,y,p)^3},\qquad
 B(s,y,p)=\frac{c(s,y)(a(s,y)^2+2p^2)
                    +\kappa_*'(s)\sin y\,p}{v(s,y,p)^3}.
\]
Thus $\kappa[u]=K(s,u,u_s,u_{ss})$, where
$K(s,y,p,z)=C(s,y,p)z+B(s,y,p)$. To keep partial derivatives of
coefficient functions distinct from derivatives along a graph, set
\[
 \mathcal D_0=\partial_s+p\partial_y+z\partial_p,\qquad
 F=\mathcal D_0K,\qquad V=\mathcal D_0v.
\]
For example, $\partial_s\kappa[u]=C u_{sss}+F$ after evaluation
at $(s,u,u_s,u_{ss})$. Two differentiations, including the factors
$v^{-1}$ in arclength differentiation along the graph, now give
\begin{align*}
 \partial_{s_\Gamma}^2\kappa[u]
 &=\frac{a}{v^5}u_{ssss}+B_0+B_1u_{sss},\\
 B_0&=v^{-2}\mathcal D_0F-v^{-3}VF,\\
 B_1&=v^{-2}(\mathcal D_0C+F_z)-v^{-3}VC.
\end{align*}
All coefficients on the right are evaluated at $(s,u,u_s,u_{ss})$.
In particular the right-hand side in \eqref{miwf:graphflow} is
\[
 Q=q_0(s,u,u_s,u_{ss})+q_1(s,u,u_s,u_{ss})u_{sss},
\]
where, with the same evaluation convention,
\[
 q_0=-\frac{v}{a}\frac{2B_0+K^3+K}{(1+K^2)^2},\qquad
 q_1=-\frac{v}{a}\frac{2B_1}{(1+K^2)^2}.
\]
These formulas prove the asserted affine dependence on $u_{sss}$.
The functions $B,B_0,B_1,q_0,q_1$ and $a_4$ are smooth in $s$
and real analytic in their finite-dimensional jet variables near the
zero jet, since they are built from trigonometric functions, positive
square roots, and quotients with nonvanishing denominators. On a
sufficiently small $C^2$
neighborhood of zero, $a\geq1/2$, $v$ is bounded above and away
from zero, and $\kappa[u]$ is bounded. Consequently $a_4\geq
c_\eta>0$, and all coefficient derivatives on the chosen jet range
are bounded. These bounds depend only on $\eta$ and that range.

Use the exponent $q=6$, $X_0=L^6$, $X_1=W^{4,6}$, with periodic
boundary conditions. The trace space is
$(X_0,X_1)_{1-1/6,6}=W^{10/3,6}$; the exponent is nonintegral.
In one dimension
\[
 W^{4,2}\hookrightarrow W^{10/3,6}\hookrightarrow C^{3,\beta}
 \quad(0<\beta<1/6).
\]
Fix the prescribed finite interval $I=[0,2\tau]$ and set
\[
 \mathbb X(I)=W^{1,6}(I;L^6)\cap L^6(I;W^{4,6}),\qquad
 \mathbb Y(I)=L^6(I;L^6)\times W^{10/3,6}.
\]
The trace embedding is
$\mathbb X(I)\hookrightarrow C(I;W^{10/3,6})$, hence also
$\mathbb X(I)\hookrightarrow C(I;C^3)$. For a sufficiently small
open neighborhood $\mathcal U$ of zero in $\mathbb X(I)$, define
\[
 \Phi(u)=\bigl(u_t+a_4(u)\partial_s^4u-Q(u),\ u(0)\bigr),
 \qquad a_4(u)=a_4(s,u,u_s,u_{ss}).
\]
We verify that $\Phi:\mathcal U\to\mathbb Y(I)$ is continuously
Fr\'echet differentiable. The finite-jet formulas give continuously
differentiable maps $u\mapsto a_4(u)$ and $u\mapsto q_i(u)$ into
$C(I\times\R/L_\eta\miwfZ)$, with derivatives given by the
finite-dimensional chain rule. For example,
\[
 Da_4(u)[h]=(a_4)_y h+(a_4)_p h_s+(a_4)_z h_{ss},
 \qquad \|Da_4(u)[h]\|_\infty\leq C\|h\|_{\mathbb X(I)}.
\]
Their derivative maps vary continuously in operator norm because
$\mathcal U$ has a common compact jet range. Multiplication from
$C(I\times\R/L_\eta\miwfZ)\times L^6(I;L^6)$ to $L^6(I;L^6)$
is bounded. Thus the principal nonlinear term has derivative
\[
 D\bigl[a_4(u)\partial_s^4u\bigr][h]
 =Da_4(u)[h]\partial_s^4u+a_4(u)\partial_s^4h,
\]
and the lower-order term has derivative
\[
 DQ(u)[h]=Dq_0(u)[h]+Dq_1(u)[h]u_{sss}+q_1(u)h_{sss}.
\]
These are bounded maps from $\mathbb X(I)$ to $L^6(I;L^6)$,
continuous in $u$. For the first term in $DQ$, the cylinder
$I\times\R/L_\eta\miwfZ$ has finite measure, so its continuous
functions embed continuously into $L^6$; the other terms use the
stated multiplication bound. The time derivative and initial
trace are bounded linear maps into their respective target spaces.
This proves the claim directly. The corresponding surface
functional-analytic construction is developed in
\cite{JakobRegularity}.

Since $\eta$ is elastic, $\Phi(0)=(0,0)$. Its derivative is
$D\Phi(0)h=(h_t+\mathcal Lh,h(0))$, where $\mathcal L$ has smooth
periodic coefficients and principal part
$2(1+\kappa_*^2)^{-2}\partial_s^4$.
It is an isomorphism by the autonomous maximal-$L^6$-regularity result
in \cite[Theorem 1.23(i) and Remark 1.31]{Amann}, with
spatial manifold $\R/L_\eta\miwfZ$, scalar bundle, order $r=4$,
$q=6$ and regularity index $s=0$. The circle is compact without
boundary, the coefficients are smooth, and the principal symbol is
positive for every nonzero covector, uniformly on the circle; these
verify the hypotheses of that theorem. The Banach inverse function
theorem supplies open neighborhoods $\mathcal U_1\subset\mathcal U$
and $\mathcal V\subset\mathbb Y(I)$ of zero on which $\Phi$ is a
$C^1$ diffeomorphism. The admissible initial data form the open
neighborhood
\[
 \mathcal V_0=\{h\in W^{10/3,6}:(0,h)\in\mathcal V\}
 \subset\operatorname{tr}_0(\mathcal U_1),
 \qquad \operatorname{tr}_0u=u(0).
\]
In particular there is $\varepsilon_{\eta,I}>0$ such that
$\|h\|_{W^{10/3,6}}<\varepsilon_{\eta,I}$ implies $h\in\mathcal V_0$.
Since the initial graph heights $h_j\to0$ in $W^{4,2}$ and hence
in $W^{10/3,6}$, the solutions
$u_j=\Phi^{-1}(0,h_j)$ exist on the same $I$ for all sufficiently
large $j$, satisfy $u_j(0)=h_j$, and obey
\begin{equation}\label{miwf:traceconvergence}
 u_j\longrightarrow0\quad\text{in }\mathbb X(I)
 \text{ and in }C(I;W^{10/3,6}).
\end{equation}
The initial-data radius $\varepsilon_{\eta,I}$ depends on the
prescribed finite interval $I$ and may shrink as $I$ grows.

The positive-time bounds are uniform. We first verify the spacetime
embedding used below, including at the initial time. Let
$0<\alpha<1/6$ and choose
\[
 \frac{19}{20}<\lambda<1-\frac{3\alpha}{10}.
\]
For functions on the fixed circle, fractional interpolation between
$L^6$ and $W^{10/3,6}$, followed by the one-dimensional Sobolev
embedding, gives
\begin{equation}\label{miwf:time-interpolation}
 \|h\|_{C^3}\leq C\|h\|_{L^6}^{1-\lambda}
                         \|h\|_{W^{10/3,6}}^\lambda.
\end{equation}
Indeed, choose a nonintegral $r$ with
$19/6<r<(10/3)\lambda$. Real interpolation gives
\[
 \|h\|_{(L^6,W^{10/3,6})_{\lambda,\infty}}
 \leq C\|h\|_{L^6}^{1-\lambda}\|h\|_{W^{10/3,6}}^\lambda,
\]
and this interpolation space embeds into $W^{r,6}$ by the strict
inequality $r<(10/3)\lambda$. The interpolation and embedding
properties on the compact circle follow by localization from their
Euclidean versions; see the function-space discussion in
\cite[Section~1]{Amann}. Finally use $W^{r,6}\hookrightarrow C^3$.
The strict inequalities avoid endpoint interpolation. For $u\in\mathbb X(I)$, the time trace
embedding and the fundamental theorem of calculus in $L^6$ give,
for $t,t'\in I$,
\[
 \sup_{r\in I}\|u(r)\|_{W^{10/3,6}}\leq C_I\|u\|_{\mathbb X(I)},
 \qquad
 \|u(t)-u(t')\|_{L^6}
 \leq |t-t'|^{5/6}\|u_t\|_{L^6(I;L^6)}.
\]
Applying \eqref{miwf:time-interpolation} to $h=u(t)-u(t')$
therefore bounds its $C^3$ norm by
$C_I\|u\|_{\mathbb X(I)}|t-t'|^{(5/6)(1-\lambda)}$.
The exponent is strictly larger than $\alpha/4$. The trace
embedding also gives $u\in C(I;C^{3,\alpha})$, since
$10/3-1/6>3+\alpha$. Combining the spatial and temporal estimates
by the triangle inequality proves
\begin{equation}\label{miwf:parabolic-embedding}
 \sum_{\ell=0}^3
 \|\partial_s^\ell u\|_{C^{\alpha,\alpha/4}(I\times\mathcal C_{L_\eta})}
 \leq C_{I,\alpha}\|u\|_{\mathbb X(I)}.
\end{equation}
All estimates extend to $t=0$ by the trace continuity; the constants
depend on the fixed circle, $I$, $\alpha$, and the interpolation
choice, not on $u$ or $j$.
Fix $0<\beta<\alpha$ and write
$b_j=a_4(s,u_j,\partial_su_j,\partial_s^2u_j)$,
$b_0=a_4(s,0,0,0)$, and $Q_j=Q(u_j)$. The finite-jet formulas,
\eqref{miwf:parabolic-embedding}, and the embedding from exponent
$\alpha$ into the little H\"older space at exponent $\beta$ give
\[
 b_j\geq c_\eta>0,\qquad
 \|b_j-b_0\|_{bc^{\beta,\beta/4}}
     +\|Q_j\|_{bc^{\beta,\beta/4}}
 \leq C\|u_j\|_{\mathbb X(I)}\longrightarrow0 .
\]
Here $Q(0)=0$ by the elastic equation for $\eta$.
With $b_j$ and $Q_j$ fixed, the linear isomorphism of
\cite[Theorem~1.23(ii)]{Amann} at $s=\bar s=\beta$ and the smooth
initial value $h_j$ give a solution in $bc^{4+\beta,1+\beta/4}$.
It belongs to $\mathbb X(I)$, so uniqueness in part~\textup{(i)}
at $s=0$, $q=6$ identifies it with $u_j$. This establishes the
basic regularity for each $j$.

For the uniform estimate take a smooth time cutoff $\chi$ vanishing
near zero and equal to one on $[\tau/2,2\tau]$. The function
$w_j=\chi u_j$ solves
\begin{equation}\label{miwf:cutoff-equation}
 \partial_tw_j+b_j\partial_s^4w_j=f_j,\qquad
 f_j=\chi Q_j+\chi'u_j,\qquad w_j(0)=0 .
\end{equation}
The preceding bounds give
$\|f_j\|_{bc^{\beta,\beta/4}}\leq C\|u_j\|_{\mathbb X(I)}\to0$.
The zero-trace estimate in \cite[Theorem~1.23(ii) and Remark~1.24]{Amann}
applied to \eqref{miwf:cutoff-equation} on the fixed cylinder gives
\[
 \|w_j\|_{bc^{4+\beta,1+\beta/4}}
 \leq C\|f_j\|_{bc^{\beta,\beta/4}}\longrightarrow0.
\]
Here $C$ is independent of $j$ by the common coefficient and
ellipticity bounds and the fixed length of $I$.

For the higher orders apply the linear Schauder estimate of
\cite[Proposition~3]{JakobShortTime}. Write $C^{r}_{\rm par}$ for
$C^{r,r/4}$. Let $J=[a,2\tau]$ with $0\leq a\leq\tau$, so that
$[\tau,2\tau]\subset J\subset I$. On the cylinder
$J\times\mathcal C_{L_\eta}$, suppose $b,F\in C^{k+\beta}_{\rm par}$
and $b\geq c>0$. For a solution $z\in C^{4+\beta}_{\rm par}$ of
$z_t+b\partial_s^4z=F$ with zero initial function at $t=a$,
the estimate is
\[
 \|z\|_{C^{4+k+\beta}_{\rm par}(J\times\mathcal C_{L_\eta})}
 \leq C_k\bigl(\|F\|_{C^{k+\beta}_{\rm par}(J\times\mathcal C_{L_\eta})}
       +\|z\|_{L^\infty(J\times\mathcal C_{L_\eta})}\bigr)
 \qquad (k\in\mathbb N_0).
\]
For each such fixed $J$, the constant depends on the fixed circle,
the length of $J$, $k$, $\beta$, a positive lower bound for $b$,
and its $C^{k+\beta}_{\rm par}$ bound on this cylinder.
In the surface formulation one extends $z$ and $F$ independently
of $y$ to a flat product torus with coordinates $(s,y)$ and uses the
principal tensor $E\otimes E$ with $E=\operatorname{diag}(\sqrt b,1)$.
On these functions the spatial operator is precisely $b\partial_s^4$.
The background metric is fixed throughout. To accommodate the
initial-immersion wording of that proposition, let $G_{\rm tor}$ be a fixed
smooth isometric product-torus immersion and apply it to the
vector-valued solution $(z,\varepsilon G_{\rm tor})$. Its initial value
$(0,\varepsilon G_{\rm tor})$ is an immersion for $\varepsilon>0$; the added
forcing and initial-value norms are $O(\varepsilon)$. The constant in
the estimate is independent of $\varepsilon$, since the background
immersion and the operator are fixed. Letting
$\varepsilon\downarrow0$ gives the stated zero-trace scalar estimate.

We now use induction with successive time cutoffs. Put
$a_k=\tau-\tau/2^{k+1}$ and $J_k=[a_k,2\tau]$. We have just
bounded $u_j$ in $C^{4+\beta}_{\rm par}(J_0)$. If it is uniformly
bounded in $C^{k+3+\beta}_{\rm par}(J_{k-1})$ for $k\geq1$, then
$b_j$ and $Q_j$ are bounded there in $C^{k+\beta}_{\rm par}$:
their jet orders are two and three, respectively. Choose a smooth
cutoff $\chi_k$ zero near $a_{k-1}$ and one on $J_k$. Then
$z_j=\chi_k u_j$ has zero trace at $a_{k-1}$ and solves
$\partial_tz_j+b_j\partial_s^4z_j=\chi_k Q_j+\chi_k'u_j$
on $J_{k-1}$. Its forcing is uniformly bounded in
$C^{k+\beta}_{\rm par}$ on $J_{k-1}$. The interval $J_{k-1}$ and cutoff are
fixed for each $k$, so the Schauder constant and cutoff bounds may
depend on $k$ but are independent of $j$. The Schauder estimate
therefore gives a uniform $C^{4+k+\beta}_{\rm par}$ bound on $J_k$
for the original solutions $u_j$. Thus, by induction, the spatial
partial derivatives of $u_j$ of every order are bounded on
$[\tau,2\tau]$, independently of $j$. Interpolating a higher spatial norm with the uniform
$C^0$ convergence in \eqref{miwf:traceconvergence} gives
$\sup_{\tau\leq t\leq2\tau}\|u_j(t)\|_{C^k}\to0$ for each $k$.

For each fixed $j$ the same induction on all of $I$ applies to
$z=u_j-h_j$, with zero trace at $t=0$ and forcing
$Q_j-b_j\partial_s^4h_j$.
The initial function $h_j$ is smooth and the spatial circle is closed.
This gives smoothness up to $t=0$ for each solution $u_j$. Higher norms of
$h_j$ enter only this last, individual bootstrap-regularity argument.

The graph solutions have the normal velocity prescribed by
\eqref{miwf:flow}. To remove their tangential velocity, write
$T_\Gamma=\Gamma_s/v$ for the unit tangent, where $v=|\Gamma_s|>0$,
and set $V_{\mathrm{tan}}=\langle\partial_t\Gamma,T_\Gamma\rangle$.
For $\widetilde\Gamma(t,x)=\Gamma(t,\rho(t,x))$, the chain rule gives
\[
 \langle\partial_t\widetilde\Gamma,
 T_\Gamma(t,\rho(t,x))\rangle
 =V_{\mathrm{tan}}(t,\rho(t,x))
  +v(t,\rho(t,x))\partial_t\rho(t,x).
\]
Thus the initial-value problem
\[
 \partial_t\rho(t,x)=-\frac{V_{\mathrm{tan}}(t,\rho(t,x))}
                              {v(t,\rho(t,x))},
 \qquad \rho(0,x)=x,
\]
removes the tangential velocity. The smooth time-dependent vector field
on the compact parameter circle generates diffeomorphisms on the whole
finite interval $I$. For the solution constructed from $h_j$, this
gives $\widetilde\Gamma_j(0)=\Gamma(h_j)$, because
$\rho_j(0)=\mathrm{id}$. The initial graph construction gave
$\eta_j=\Gamma(h_j)\circ\sigma_j$. Thus the normal-gauge solution
with initial function $\eta_j$ is explicitly
\[
 \gamma_j(t,x)=\widetilde\Gamma_j(t,\sigma_j(x)),\qquad
 \gamma_j(0,x)=\eta_j(x).
\]
The fixed map $\sigma_j$ accounts only for that initial graph
construction; being independent of time, it preserves the normal
velocity equation.

To verify uniqueness, let $u_1,u_2$ be any two smooth solutions of
\eqref{miwf:graphflow} on a common interval with the same initial
height $u_1(0)=u_2(0)$, whose jets remain in the chosen convex graph
neighborhood. Put $w=u_1-u_2$, so $w(0)=0$. Their difference satisfies
\[
 w_t+a_4(u_1)w_{ssss}
 =Q(u_1)-Q(u_2)-\bigl(a_4(u_1)-a_4(u_2)\bigr)\partial_s^4u_2.
\]
Apply the classical mean-value theorem in integral form along the
segment between the two jets at fixed $s$. Each difference of the
coefficients $a_4,q_0,q_1$ is a linear combination of $w,w_s,w_{ss}$,
with coefficients given by the averages of the corresponding jet
derivatives along this segment. Since $Q$ is affine in $u_{sss}$,
the equation becomes
\[
 w_t+b_4w_{ssss}=\sum_{\ell=0}^3 b_\ell\partial_s^\ell w,
 \qquad b_4=a_4(u_1)>0.
\]
Here $b_3=q_1(u_1)$. The resulting coefficients involve derivatives of
$u_2$ through order four and are therefore smooth and bounded on
every common compact time interval. On such an interval choose
$\beta_0>0$ with $b_4\geq\beta_0$.
All integrals and norms in the following calculation are taken on
the fixed parameter circle. Periodic integration by parts gives
\[
 \int b_4ww_{ssss}\dd s
 =\int\bigl(b_4w_{ss}^2+2(\partial_sb_4)w_sw_{ss}
                  +(\partial_s^2b_4)ww_{ss}\bigr)\dd s,
\]
and the third-order term on the right satisfies
\[
 \int b_3ww_{sss}\dd s
 =-\int\bigl(b_3w_s+(\partial_sb_3)w\bigr)w_{ss}\dd s.
\]
Young's inequality bounds the remaining terms, for every
$\varepsilon>0$, by
$\varepsilon\|w_{ss}\|_2^2+C_\varepsilon\|w\|_{H^1}^2$.
The periodic interpolation inequality
$\|w_s\|_2^2\leq\|w\|_2\|w_{ss}\|_2$, followed by another
application of Young's inequality, allows us to absorb these terms.
Taking the $L^2$ inner product of the equation with $w$ and
absorbing the lower-order terms gives
\[
 \frac{\dd}{\dd t}\|w\|_2^2
       +\beta_0\|w_{ss}\|_2^2\leq C\|w\|_2^2.
\]
The constant $C$ uses only the indicated coefficient bounds on this
compact time interval and the positive lower bound $\beta_0$.
Gronwall's inequality proves uniqueness for equal initial data.
Let $\gamma_t,\widetilde\gamma_t$ be two smooth normal-gauge solutions
with the same parametrized initial curve. At any time $t_0$ at which
the two maps agree, use their common curve as the reference curve,
after the same fixed arclength reparametrization. For a short time
after $t_0$, both maps have graphs on the inverse branches through
$(s,0)$. The two heights start at zero. The difference estimate above
applies to this smooth reference curve, since it only uses the
smooth graph coefficients and $a_4>0$. Thus the heights agree.
The two base-coordinate diffeomorphisms agree at $t_0$ and give
a smooth relative reparametrization $\varphi_t$ on the short interval
such that
\[
 \widetilde\gamma_t=\gamma_t\circ\varphi_t,\qquad
 \varphi_{t_0}=\mathrm{id}.
\]
Differentiation gives
\[
 \partial_t\widetilde\gamma_t
 =(\partial_t\gamma_t)\circ\varphi_t
   +D\gamma_t(\varphi_t)\partial_t\varphi_t.
\]
The first two velocities are normal to the curve. Tangential
projection therefore gives
$D\gamma_t(\varphi_t)\,\partial_t\varphi_t=0$.
Regularity implies $\partial_t\varphi_t=0$, hence
$\varphi_t=\mathrm{id}$ on this interval. Equality of the parametrized
maps is closed in time by continuity and extends to the right at each
time of agreement. Starting at zero gives equality on their entire
common interval, including for multiply covered initial curves.
\end{proof}

\begin{corollary}[Bounded-dissipation endpoint regularity and energy identity]
\label{miwf:endpoint-W42}
Assume the hypotheses of Theorem~\ref{thm:all-energy-profile}, and
let $\gamma_T$ be its unique anchored
endpoint profile.  If there are times $t_j\uparrow T$ and a constant
$D_0<\infty$ such that $\miwfD(\gamma_{t_j})\leq D_0$, then
$\gamma_T\in W^{4,2}(\Sph^1,\R^3)$ with values in $\Sph^2$.
Along a subsequence,
the anchored profiles converge weakly to $\gamma_T$ in
$W^{4,2}$ and strongly in $C^{3,\alpha}$ for every
$0<\alpha<\frac12$.  Moreover,
\begin{equation}\label{miwf:endpoint-energy-identity}
 \miwfE(\gamma_T)=\lim_{t\uparrow T}\miwfE(\gamma_t),\qquad
 \widetilde\gamma_t\longrightarrow\gamma_T
       \text{ strongly in }W^{2,2}\text{ as }t\uparrow T.
\end{equation}
The endpoint profile is regular; in the
low-energy setting of Theorem~\ref{thm:anchored-profile} it is also
embedded.  Bounded dissipation alone does not imply that it is elastic.
\end{corollary}
\begin{proof}
The energy along the profiles is bounded by $E_0$, so
Proposition~\ref{miwf:bounded-dissipation-compactness} applies to the
anchored constant-speed representatives at the times $t_j$.
It gives a subsequence converging weakly in $W^{4,2}$ and strongly in
$C^{3,\alpha}$ to a regular $W^{4,2}$ curve.  Theorem
\ref{thm:all-energy-profile} already gives full $L^2$ convergence of
these same anchored profiles to $\gamma_T$.  Thus the
subsequential limit is $\gamma_T$.  Regularity follows from that
theorem, and the stated low-energy embeddedness follows from
Theorem~\ref{thm:anchored-profile}.

For a constant-speed profile on $\R/\mathbb Z$, the frame equation
\eqref{miwf:frame} and the energy definition in \eqref{miwf:flow} give
$\miwfE(\gamma)=L^{-3}\int|\gamma_{xx}|^2\dd x$. Thus strong
$C^2$ convergence and positive limiting speed imply
$\miwfE(\gamma_{t_j})\to\miwfE(\gamma_T)$ along the subsequence.
Energy monotonicity identifies this subsequential limit with
$\lim_{t\uparrow T}\miwfE(\gamma_t)$. For the complete anchored trajectory,
\[
 \|\partial_x^2\widetilde\gamma_t\|_2^2
 =L_t^3\miwfE(\gamma_t)\longrightarrow L_T^3\miwfE(\gamma_T)
 =\|\partial_x^2\gamma_T\|_2^2.
\]
The full weak $W^{2,2}$ convergence from
Theorem~\ref{thm:all-energy-profile} and convergence of these Hilbert
norms imply strong $L^2$-convergence of the second derivatives.
The lower derivatives already converge in $C^0$, proving
\eqref{miwf:endpoint-energy-identity}.
\end{proof}

\begin{corollary}[Dissipation obstruction at finite maximal time]
\label{miwf:finite-max-dissipation}
Let $\gamma_t$ be the smooth maximal normal-gauge solution of
\eqref{miwf:flow} with $T_{\max}(\gamma_0)<\infty$. Then
\[
 \liminf_{t\uparrow T_{\max}}\miwfD(\gamma_t)>0.
\]
\end{corollary}
\begin{proof}
If the liminf were zero, choose $t_j\uparrow T_{\max}$ with
$\miwfD(\gamma_{t_j})\to0$.  Proposition~\ref{miwf:ps}, with
constant-speed reparametrizations, gives a subsequence converging
strongly in $W^{4,2}$ to a smooth regular elastic curve.  Apply
Lemma~\ref{miwf:smoothing} to this sequence.  It supplies a common
$\tau>0$ for which the solutions starting at the reparametrized
$\gamma_{t_j}$ exist through time $t_j+2\tau$.  For $j$ large,
$T_{\max}-t_j<\tau$. Let $\chi_j$ be the fixed diffeomorphism used
at time $t_j$, so the restarted graph solution initially represents
$\gamma_{t_j}\circ\chi_j$. Remove its tangential velocity on
$[0,2\tau]$ by the ODE above, with initial diffeomorphism the identity,
and then compose with $\chi_j^{-1}$. This gives a smooth normal-gauge
solution on the original parameter circle, initially equal to
$\gamma_{t_j}$. Normal-gauge uniqueness identifies it with the
original solution on their common time interval. It exists through
$t_j+2\tau>T_{\max}$ and therefore extends that solution beyond
$T_{\max}$, contradicting maximality.
\end{proof}

\begin{remark}[Finite-time dissipation]
At a finite maximal endpoint, finite total dissipation
$\int_0^{T_{\max}}\miwfD\,\dd t
=\miwfE(\gamma_0)-E_{T_{\max}}$ coexists with the positive
lower limiting dissipation in
Corollary~\ref{miwf:finite-max-dissipation}.
\end{remark}

\begin{samepage}
\subsection{Smooth subsequential convergence of global trajectories}

\begin{theorem}\label{miwf:global}
Every smooth global solution of \eqref{miwf:flow} with a smooth regular
closed initial curve has times $s_j\to\infty$, smooth diffeomorphisms
$\psi_j$ of the parameter circle, and a smooth regular closed elastic
curve $\gamma_\infty$ such that
\[
 \gamma_{s_j}\circ\psi_j\longrightarrow\gamma_\infty
 \quad\text{in }C^k\text{ for every integer }k\geq0.
\]
The same sequence can be used for all $k$.
\end{theorem}
\end{samepage}

\begin{proof}
Choose the times and the smooth elastic limit given by Corollary
\ref{miwf:goodtimes}. Apply Lemma \ref{miwf:smoothing} at those times with
$\eta=\gamma_\infty$. After the fixed reparametrizations that identify
the initial profiles in Lemma~\ref{miwf:smoothing}, uniqueness of
normal-gauge solutions of \eqref{miwf:flow} identifies the restarted
solutions with the time-translated trajectory in that gauge. Thus their
values at the common positive time $\tau$ are reparametrizations of
$\gamma_{t_j+\tau}$. Set $s_j=t_j+\tau$.
Then $s_j\to\infty$ and the claimed smooth convergence follows.
\end{proof}

\begin{remark}[The regularizing time shift]
Proposition~\ref{miwf:ps} supplies strong $W^{4,2}$ convergence
at the selected times. The common positive time shift in the
preceding proof gives uniform bounds of every spatial derivative
through Lemma~\ref{miwf:smoothing}.
\end{remark}

\subsection{A fixed-gauge \L ojasiewicz--Simon inequality and full convergence}
\label{miwf:LS-section}

Let $\eta$ be a regular closed elastic curve in $\Sph^2$, parametrized by
arclength $s$ on $\R/L_\eta\miwfZ$, and put
$T_\eta=\eta_s$ and $\nu_\eta=\eta\times T_\eta$. For a small function
$u\in H^2(\R/L_\eta\miwfZ)$ define its normal graph in $\Sph^2$ by
\[
 \Gamma_u(s)=\cos u(s)\,\eta(s)+\sin u(s)\,\nu_\eta(s),
 \qquad \mathcal E(u)=\miwfE(\Gamma_u).
\]
Write $\mathcal C_\eta=\R/L_\eta\miwfZ$. Throughout this argument,
\[
 H^2=H^2(\mathcal C_\eta;\R)=W^{2,2}(\mathcal C_\eta;\R),
 \qquad H^{-2}=(H^2)^*;
\]
the duality convention is induced by the scalar $L^2$ pairing.
The height $u$ is scalar, whereas $\Gamma_u$ is an $\R^3$-valued
map whose image lies in $\Sph^2$. Thus the graph construction is
the nonlinear map
\[
 \mathcal G:H^2(\mathcal C_\eta;\R)\longrightarrow
 W^{2,2}(\mathcal C_\eta;\R^3),\qquad
 \mathcal G(u)=\cos u\,\eta+\sin u\,\nu_\eta.
\]
Only a small neighborhood of zero is used as a chart of regular
curves. The coordinate extraction for a nearby curve $\gamma$ is
instead $\sigma=P_\alpha\circ\gamma$, $z=R_\alpha\circ\gamma$,
and $u=z\circ\sigma^{-1}$ on the selected inverse branches from
Lemma~\ref{miwf:smoothing}. It then gives
$\gamma=\mathcal G(u)\circ\sigma$.
After a reparametrization close to the identity, this graph chart fixes
the local reparametrization freedom.  It is valid for immersed curves,
including self-intersections and multiple covers, since it is defined
on the parameter circle rather than by a tubular neighborhood of the
image.

\begin{proposition}[Fixed-gauge \L ojasiewicz--Simon inequality]
\label{miwf:LS}
There are $\rho>0$, $C>0$, and $\theta\in(0,\frac12]$, depending on
$\eta$, such that
\begin{equation}\label{miwf:LS-ineq}
 |\mathcal E(u)-\mathcal E(0)|^{1-\theta}
 \leq C\|D\mathcal E(u)\|_{H^{-2}}
 \qquad (\|u\|_{H^2}<\rho).
\end{equation}
The dual norm is taken with respect to the fixed arclength coordinate
of $\eta$.  After decreasing $\rho$ and enlarging $C$, the inequality
also has the following form adapted to the reduced MIWF:
\begin{equation}\label{miwf:LS-weighted}
 |\mathcal E(u)-\mathcal E(0)|^{1-\theta}
 \leq C\sqrt{\miwfD(\Gamma_u)}
 \qquad (u\text{ smooth},\ \|u\|_{H^2}<\rho).
\end{equation}
Neither a bound on the curvature throughout the graph ball nor
nondegeneracy of $\eta$ is assumed.
\end{proposition}
\begin{proof}
Use the graph quantities from \eqref{miwf:graphcurvature}: set
$p=u_s$, $a=\cos u-\kappa_\eta\sin u$,
$c=\sin u+\kappa_\eta\cos u$, and $v=(a^2+p^2)^{1/2}$.  Write
\[
 B(s,u,p)=\frac{c(a^2+2p^2)+\kappa_\eta'(s)\sin u\,p}{v^3}.
\]
Then $\kappa_{\Gamma_u}=(a/v^3)u_{ss}+B$, and hence
\begin{equation}\label{miwf:LS-energy-form}
 \mathcal E(u)=\int_0^{L_\eta}
 \left(\frac{a^2}{v^5}u_{ss}^2
       +\frac{2aB}{v^2}u_{ss}+v(1+B^2)\right)\dd s.
\end{equation}
Take an open ball $U$ about zero in $V=H^2(\mathcal C_\eta;\R)$.
By $H^2\hookrightarrow C^1$ we may choose it so that
$\|v^2-1\|_{C^0}<1/2$. In the Banach algebra $C^0(\mathcal C_\eta)$,
sine and cosine have convergent power series, and the binomial series
for each power of $v^2$ in \eqref{miwf:LS-energy-form} converges
on this ball. Multiplication by the fixed smooth $\kappa_\eta$
and $\kappa_\eta'$ is bounded. Thus all three coefficients in the
integrand depend real analytically on $u$ with values in $C^0$.
The bounded maps $u\mapsto u_{ss}$ from $H^2$ to $L^2$ and
$L^2\times L^2\to L^1$ given by multiplication now prove analyticity
of $\mathcal E:U\to\R$. Its Fr\'echet derivative
$M=D\mathcal E:U\to W:=V^*=H^{-2}$ is analytic in the dual norm,
and $M(0)=0$ by the elastic equation and the first variation.

We apply \cite[Corollary~3.2]{MantegazzaPozzetta} with
$V=Z=H^2$ and $W=V^*=H^{-2}$. The gradient $M:U\to W$ is
analytic, $M(0)=0$, and its linearization
$L:=DM(0):H^2\to H^{-2}$ has principal part $2\partial_s^4$.
Indeed, the coefficient of $u_{ss}^2$ in
\eqref{miwf:LS-energy-form} equals one at the zero jet; the remaining
terms in the second variation have total differential order at most
three. The Fourier multiplier $2(\partial_s^4+1)$ is an isomorphism
from periodic $H^2$ to $H^{-2}$. Its difference from $L$ maps
$H^2$ boundedly to $H^{-1}$, which embeds compactly into $H^{-2}$.
Thus $L$ is Fredholm of index zero. The cited corollary yields
\eqref{miwf:LS-ineq} directly, with the dual norm on the full
parameter circle.
The higher-regularity realization in
\cite[Appendix~A]{DallAcquaPozziSpener} provides a related application;
for spherical elastic flow see \cite[Remark~1.4]{MantegazzaPozzetta}.

We now prove \eqref{miwf:LS-weighted}.  Shrink the graph ball so that
$|\mathcal E(u)-\mathcal E(0)|\leq1$ and
$\mathcal E(u)\leq\mathcal E(0)+1$.  For a smooth graph, vary $u$ in the direction $h$.  The variation
field is $h(-\sin u\,\eta+\cos u\,\nu_\eta)$.  Its scalar
component in the unit normal direction to $\Gamma_u$ in $\Sph^2$
is $(a/v)h$, with the normal convention of
\eqref{miwf:graphcurvature}.  Since
$\dd s_{\Gamma_u}=v\dd s$, the first variation formula is
\begin{equation}\label{miwf:LS-first-var}
 D\mathcal E(u)[h]=\int_0^{L_\eta}aG(\Gamma_u)h\,\dd s.
\end{equation}
If $\miwfD(\Gamma_u)\leq1$, the curvature estimate
\eqref{miwf:max}, applied with $E_0=\mathcal E(0)+1$, bounds
$\|\kappa_{\Gamma_u}\|_\infty$ uniformly.
Writing $\dd s_{\Gamma_u}=v\dd s$, Cauchy--Schwarz yields
\[
 |D\mathcal E(u)[h]|
 \leq\sqrt{\miwfD(\Gamma_u)}
 \left(\int_{\Gamma_u}
   (1+\kappa_{\Gamma_u}^2)^2\frac{a^2}{v^2}h^2\,
   \dd s_{\Gamma_u}\right)^{1/2}
 \leq C\sqrt{\miwfD(\Gamma_u)}\|h\|_{H^2}.
\]
The graph ball bounds $a$ and $v$ above and away from zero.  Applying
\eqref{miwf:LS-ineq} proves the weighted inequality in this case.
If $\miwfD(\Gamma_u)>1$, its right-hand side, with $C\geq1$, is at
least $1$, whereas its left-hand side is at most $1$.  This proves
\eqref{miwf:LS-weighted} in the remaining case.  This split according
to dissipation supplies the weight control on the $H^2$ neighborhood.
\end{proof}

The following estimate supplies the no-exit step and the convergence rate.
\begin{lemma}[Energy control in the normal graph chart]
\label{miwf:LS-compactness}
There are $r_*>0$ and $C_*>0$, depending on $\eta$, such that
\begin{equation}\label{miwf:graph-energy-control}
 \|u\|_{H^2}^2
 \leq C_*\bigl(|\mathcal E(u)-\mathcal E(0)|+\|u\|_{L^2}^2\bigr)
 \qquad (\|u\|_{H^2}\leq r_*).
\end{equation}
Consequently, for every $0<r\leq r_*$ there are
$\varepsilon,\delta>0$ such that
\[
 \|u\|_{H^2}\leq r,\quad \|u\|_{L^2}\leq\delta,\quad
 |\mathcal E(u)-\mathcal E(0)|\leq\varepsilon
 \quad\Longrightarrow\quad \|u\|_{H^2}<r/2.
\]
\end{lemma}
\begin{proof}
Write the integrand in \eqref{miwf:LS-energy-form} as
$\mathsf e_2(s,u,u_s)u_{ss}^2+\mathsf e_1(s,u,u_s)u_{ss}
 +\mathsf e_0(s,u,u_s)$.
For $r_*$ small, $\mathsf e_2\geq c_*>0$ uniformly.  Moreover,
$\mathsf e_1(s,0,0)=2\kappa_\eta(s)$, and smoothness on the compact
coefficient domain gives the pointwise estimates
\begin{align*}
 |\mathsf e_1(s,u,u_s)-\mathsf e_1(s,0,0)|
 &\leq C(|u|+|u_s|),\\
 |\mathsf e_0(s,u,u_s)-\mathsf e_0(s,0,0)
       -\mathsf e_{0,u}(s,0,0)u-\mathsf e_{0,p}(s,0,0)u_s|
 &\leq C(u^2+u_s^2).
\end{align*}
The integral of
$\mathsf e_1(s,0,0)u_{ss}+\mathsf e_{0,u}(s,0,0)u
 +\mathsf e_{0,p}(s,0,0)u_s$ is $D\mathcal E(0)[u]=0$.
Subtracting it from the energy difference and using Young's
inequality therefore gives
\[
 \mathcal E(u)-\mathcal E(0)
 \geq \frac{c_*}{2}\|u_{ss}\|_2^2-C\|u\|_{H^1}^2.
\]
Periodic integration by parts gives
$\|u_s\|_2^2\leq\|u\|_2\|u_{ss}\|_2$.
Absorbing a sufficiently small multiple of $\|u_{ss}\|_2^2$ yields
\[
 \mathcal E(u)-\mathcal E(0)
 \geq c\|u\|_{H^2}^2-C\|u\|_2^2
\]
with $c>0$, proving \eqref{miwf:graph-energy-control} after increasing
$C_*$. All constants depend only on the fixed curve $\eta$
and on the size of the graph ball.
Choose $\varepsilon,\delta>0$ so that
$C_*(\varepsilon+\delta^2)<r^2/4$; this proves the final assertion.
\end{proof}

\begin{lemma}[Identification of a limit in the continued graph gauge]
\label{miwf:graph-orbit}
Let $\eta$ be a smooth regular closed curve parametrized by arclength
on $\R/L_\eta\miwfZ$, and use its normal graphs $\Gamma_u$ as above.
There exists $r_0>0$ with the following property. If $u_j$ is bounded
by $r_0$ in $H^2$ and
\[
 \Gamma_{u_j}\circ\varphi_j\longrightarrow\eta\quad\text{in }C^1
\]
for orientation-preserving circle diffeomorphisms $\varphi_j$, then
$u_j\to0$ in $C^1$, and hence in $L^2$.
No embeddedness or single-cover assumption is needed.
\end{lemma}
\begin{proof}
We first prove a local uniqueness assertion: a sufficiently small
$C^1$ normal graph that parametrizes $\eta$ has zero height.
Otherwise there are nonzero $w_n\to0$ in $C^1$ and
orientation-preserving $C^1$ diffeomorphisms $\sigma_n$ such that
$\Gamma_{w_n}=\eta\circ\sigma_n$. Since $\eta$ has unit speed,
$\sigma_n'=|\partial_s\Gamma_{w_n}|\to1$ uniformly.
Choose lifts of $\sigma_n$ with $\sigma_n(0)\in[0,L_\eta]$ and pass
to a subsequence with $\sigma_n(0)\to b$. Integration gives
$\sigma_n(s)\to s+b$ in $C^1$, and hence $\eta(s+b)=\eta(s)$.
Thus $\widehat\sigma_n:=\sigma_n-b$ tends to the identity and
$\eta\circ\widehat\sigma_n=\Gamma_{w_n}$.
The map $(s,z)\mapsto\cos z\,\eta(s)+\sin z\,\nu_\eta(s)$
is a local diffeomorphism near every $(s,0)$. Compactness of the
parameter circle provides a uniform size for these local inverse
branches. Applying their uniqueness to $(s,w_n(s))$ and
$(\widehat\sigma_n(s),0)$ gives $w_n=0$ for large $n$, a contradiction.
Let $\epsilon_0>0$ be a $C^1$ radius for this orbit separation.

Choose $r_0$ so small that the closed $H^2$ ball of radius $r_0$
lies in $\{\|u\|_{C^1}<\epsilon_0/2\}$ and all its graph speeds
are bounded above and away from zero.
Given any subsequence, compact embedding supplies a further
subsequence with $u_j\to u$ in $C^1$.
Put $f_j=\Gamma_{u_j}\circ\varphi_j$. The identity
\[
 \varphi_j'
 =\frac{|f_j'|}{|\partial_s\Gamma_{u_j}|\circ\varphi_j}
\]
first bounds $\varphi_j'$ above and away from zero. After choosing
lifts and extracting, the maps converge uniformly to a circle
homeomorphism $\varphi$. The same identity, the uniform convergence
of the graph speeds, and $f_j\to\eta$ in $C^1$ then give
$C^1$ convergence to a diffeomorphism $\varphi$.
Consequently $\Gamma_u\circ\varphi=\eta$.
The orbit separation just proved gives $u=0$.
Every subsequence therefore has a further subsequence converging
to zero in $C^1$, proving the assertion.
\end{proof}

Lemma~\ref{miwf:graph-orbit} identifies the limiting graph height
under arbitrary orientation-preserving reparametrizations. The maps
themselves can retain a symmetry of the reference curve: on the
$m$-fold equator, $m\geq2$, the shift $s\mapsto s+2\pi$ preserves
the curve and its zero graph. More generally, if
$\eta(s+b)=\eta(s)$, then
$\Gamma_u(s+b)=\Gamma_{u(\cdot+b)}(s)$. Heights on distinct choices
of sheet can therefore differ. Lemma~\ref{miwf:graph-orbit} identifies
the zero height in the limiting orbit, which suffices once the initial
choice of sheet has been continued.

\begin{theorem}[Full smooth convergence of global reduced trajectories]
\label{miwf:full}
Every smooth global solution of \eqref{miwf:flow} with a smooth regular
closed initial curve has a smooth regular closed elastic limit $\eta$
and reparametrizations $\psi_t$, defined for all sufficiently large
$t$, such that
\[
 \gamma_t\circ\psi_t\longrightarrow\eta
 \quad\text{in }C^k\text{ for every integer }k\geq0
 \quad\text{as }t\to\infty.
\]
Moreover, the reduced dissipation has finite tail length:
\[
 \int_{t_0}^{\infty}\sqrt{\miwfD(\gamma_t)}\,\dd t<\infty
\]
for some $t_0<\infty$.
No uniform curvature bound along the trajectory is assumed.
With $\theta$ as in Proposition~\ref{miwf:LS}, and
$\Delta(t)=\miwfE(\gamma_t)-\miwfE(\eta)$, there are constants
$C,c>0$ and $t_0<\infty$ such that, for $t\geq t_0$,
\begin{equation}\label{miwf:full-energy-distance}
 \|\gamma_t\circ\psi_t-\eta\|_{W^{2,2}}\leq C\Delta(t)^\theta.
\end{equation}
If $0<\theta<1/2$, then
\begin{align}
 \Delta(t)&\leq C(1+t-t_0)^{-1/(1-2\theta)},\label{miwf:rate-energy}\\
 \|\gamma_t\circ\psi_t-\eta\|_{W^{2,2}}
 &\leq C(1+t-t_0)^{-\theta/(1-2\theta)}.\label{miwf:rate-curve}
\end{align}
If $\theta=1/2$, the respective bounds are
$Ce^{-c(t-t_0)}$ and $Ce^{-c(t-t_0)/2}$.
The comparison maps and norms in these estimates use the fixed
arclength circle $\R/L_\eta\miwfZ$ as their common domain.
\end{theorem}
\begin{proof}
Theorem~\ref{miwf:global} supplies times $t_j\to\infty$ and a smooth
regular elastic curve $\eta$ such that, after reparametrization,
$\gamma_{t_j}\to\eta$ smoothly.  In particular these curves are
normal graphs $u_j$ over $\eta$ with $u_j\to0$ in $H^2$.  Energy
monotonicity and this subsequence imply
\[
 \miwfE(\gamma_t)\downarrow\miwfE(\eta),\qquad
 \Delta(t):=\miwfE(\gamma_t)-\miwfE(\eta)\geq0,
 \quad \Delta(t)\to0.
\]
All reparametrizations in that construction preserve orientation.
After a fixed orientation-preserving linear change of material
parameter, regard $\gamma_t$ also as a map on $\mathcal C_\eta$.
Choose $z_0>0$ so small that
\[
 a(s,z):=\cos z-\kappa_\eta(s)\sin z\geq\tfrac12
 \quad(s\in\mathcal C_\eta,\ |z|\leq z_0).
\]
Fix $0<r<\min\{\rho,r_*,r_0\}$, with the radii from
Proposition~\ref{miwf:LS} and
Lemmas~\ref{miwf:LS-compactness}--\ref{miwf:graph-orbit}, small enough
that $\|u\|_{H^2}\leq r$ implies $\|u\|_\infty<z_0/2$.
The embedding $H^2\hookrightarrow C^1$ then supplies a constant
$V_r$, depending only on $\eta,r$, such that the graph speed
\[
 v_u(s)=\sqrt{a(s,u(s))^2+u_s(s)^2}
 \quad\hbox{satisfies}\quad \tfrac12\leq v_u(s)\leq V_r
 \quad\hbox{on this closed ball}.
\]
Choose $\varepsilon,\delta$ from
Lemma~\ref{miwf:LS-compactness} for that ball.
While the trajectory is represented
there by $u(t)$, Proposition~\ref{miwf:LS} gives directly
\begin{equation}\label{miwf:LS-dissipation}
 \Delta(t)^{1-\theta}\leq C\sqrt{\miwfD(\Gamma_{u(t)})}.
\end{equation}
In the normal graph gauge, \eqref{miwf:flow} gives
\[
 u_t=-\frac{v}{a}\frac{G(\Gamma_u)}{(1+\kappa_{\Gamma_u}^2)^2}.
\]
Using $\dd s_{\Gamma_u}=v\dd s$, $a\geq1/2$, and $v\leq V_r$
shows explicitly how this velocity is controlled by dissipation:
\[
 \|u_t\|_{L^2(\dd s)}^2
 =\int_{\Gamma_u}\frac{v}{a^2(1+\kappa_{\Gamma_u}^2)^2}
       \frac{G(\Gamma_u)^2}{(1+\kappa_{\Gamma_u}^2)^2}
       \dd s_{\Gamma_u}
 \leq4V_r\miwfD(\Gamma_u).
\]
The constant uses only the displayed bounds for $a$ and $v$.
By the energy identity
$\Delta'=-\miwfD$, inequality \eqref{miwf:LS-dissipation} implies
\begin{equation}\label{miwf:LS-length}
 \int_{t_1}^{t_2}\sqrt{\miwfD(\Gamma_{u(t)})}\,\dd t
 \leq C\Delta(t_1)^\theta
 \qquad (t_2\geq t_1),
\end{equation}
for every segment that stays in the graph ball and has $\Delta(t_1)\leq1$.
Indeed, wherever $\Delta>0$,
$-(\Delta^\theta)'=\theta\Delta^{\theta-1}\miwfD
 \geq(\theta/C)\sqrt{\miwfD}$, and integration proves the estimate.
If $\Delta$ reaches zero, monotonicity and the already identified
limiting energy give $\miwfD=0$ at every later time.  The normal flow
is then stationary, so the same estimate remains valid.

Take $t_1=t_j$ sufficiently large that $u_j$ is in the inner half of
the graph ball, $\Delta(t_j)<\min\{\varepsilon,1\}$, and
$\|u_j\|_{L^2}+C\Delta(t_j)^\theta<\delta$, where $C$ also
includes the constant in the estimate for $\|u_t\|_2$.
Starting with the graph representation at $t_j$, continue its chosen
local inverse branches while $\|u(t)\|_{H^2}<r$, and write
\[
 \gamma_t=\Gamma_{u(t)}\circ\sigma_t,
 \qquad \sigma_t:\mathcal C_\eta\longrightarrow\mathcal C_\eta.
\]
The maps $\sigma_t$ are orientation-preserving circle diffeomorphisms.
Let $[t_j,t_*)$ be the maximal interval of this continuation inside
the open graph ball. For every $t<t_*$, the velocity estimate and
\eqref{miwf:LS-length} give
\[
 \|u(t)\|_2\leq\|u_j\|_2+C\Delta(t_j)^\theta<\delta,
 \qquad 0\leq\Delta(t)\leq\Delta(t_j)<\varepsilon.
\]
Lemma~\ref{miwf:LS-compactness} therefore gives
$\|u(t)\|_{H^2}<r/2$ throughout this interval. We prove explicitly
that a finite $t_*$ cannot be a failure of the coordinate
representation either.

Use the normal-coordinate map and its local inverses from
Lemma~\ref{miwf:smoothing}, now written as
\[
 \begin{aligned}
 \Xi(s,z)&=\cos z\,\eta(s)+\sin z\,\nu_\eta(s),\\
 b(s,z)&=-\sin z\,\eta(s)+\cos z\,\nu_\eta(s),\\
 (P_\alpha,R_\alpha)&=\Xi_\alpha^{-1}.
 \end{aligned}
\]
Since $\Xi_s=aT_\eta$, $\Xi_z=b$, and $T_\eta$ and $b$ are
orthonormal, the bound $a\geq1/2$ makes $\Xi$ a local
diffeomorphism on $|z|<z_0$. On a chosen branch the two coordinates
of the given curve are
\[
 P_\alpha(\gamma_t(x))=\sigma_t(x),\qquad
 R_\alpha(\gamma_t(x))=u(t,\sigma_t(x)).
\]
Thus $P_\alpha$ extracts the base parameter on the selected
local branch. Differentiating
$\gamma_t=\Gamma_{u(t)}\circ\sigma_t$ and using
$\partial_s\Gamma_u=aT_\eta+u_s b$ gives
\[
 \partial_x\gamma_t
   =(aT_\eta+u_s b)\circ\sigma_t\,\sigma_t',\qquad
 (P_\alpha\circ\gamma_t)'(x)=\sigma_t'(x)
   =\frac{|\partial_x\gamma_t(x)|}{v_{u(t)}(\sigma_t(x))}>0.
\]
Here the sign is fixed by the orientation of $\sigma_t$.
The graph also satisfies a uniform angle estimate. The coordinate
fiber through $\Xi(s,z)$ is the curve $z\mapsto\Xi(s,z)$ with
$s$ fixed, and its unit tangent is $b(s,z)$. If
$T_u=\partial_s\Gamma_u/v_u$ and $\alpha_u(s)\in[0,\pi/2]$
is the smaller angle between the graph tangent line and this fiber
tangent line at $\Gamma_u(s)$, then
\[
 \sin\alpha_u(s)=|\langle T_u(s),T_\eta(s)\rangle|
       =\frac{a(s,u(s))}{v_u(s)}\geq\frac1{2V_r}>0.
\]
Hence the tangent line of the graph is uniformly transverse to
the normal-coordinate fiber on each branch, throughout the chosen
$H^2$ ball.

Suppose now that $t_*<\infty$. The original global trajectory is
smooth and regular on a compact time interval containing
$[t_j,t_*]$. Hence, for constants $m_*,M_*,B_*>0$ depending on that
interval and the original trajectory,
\[
 m_*\leq|\partial_x\gamma_t|\leq M_*,\qquad
 |\partial_t\gamma_t|\leq B_*.
\]
The preceding derivative formula implies
$m_*/V_r\leq\sigma_t'\leq2M_*$ for $t<t_*$.
Write $z_t(x)=u(t,\sigma_t(x))$. At each fixed material point $x$,
time differentiation of $\gamma_t(x)=\Xi(\sigma_t(x),z_t(x))$ gives
\[
 \partial_t\sigma_t(x)=
 \frac{\langle\partial_t\gamma_t(x),T_\eta(\sigma_t(x))\rangle}
      {a(\sigma_t(x),z_t(x))},\qquad
 \partial_t z_t(x)=
 \langle\partial_t\gamma_t(x),b(\sigma_t(x),z_t(x))\rangle .
\]
Both coordinate speeds are bounded by $2B_*$. Lift the material
coordinate $x$ to $\R$ and continue a degree-one lift of $\sigma_{t_j}$
in time. The lifted base coordinate and the height
are uniformly Cauchy as $t\uparrow t_*$. Their limits $(\sigma_*,z_*)$
satisfy $\Xi(\sigma_*,z_*)=\gamma_{t_*}$ and
$|z_*|\leq z_0/2$. The lift of $\sigma_*$ has increment $L_\eta$
over one period and inherits the two-sided Lipschitz bounds.

For each material point choose an inverse of $\Xi$ at its limiting
pair $(\sigma_*,z_*)$. Choose finitely many source intervals covering the circle such that
the limiting pairs on the closure of each interval lie in its chosen
inverse chart. For $t$ close to
$t_*$ evaluate these inverses on $\gamma_t$ to obtain smooth local
pairs $(\sigma(t,x),z(t,x))$ on both sides of $t_*$. They agree on
overlaps because they approach the same limiting pairs uniformly
on the closures. Indeed, two sufficiently close points in the compact strip
$|z|\leq 3z_0/4$ with the same $\Xi$-image coincide. The latter
assertion follows by compactness and local injectivity at any common
limit of two such points. The same uniqueness identifies the pairs
with $(\sigma_t,z_t)$ before $t_*$.

Thus $\sigma_*$ is smooth, with $\partial_x\sigma_*\geq m_*/V_r$.
By compactness, the continued base maps still have positive
derivative near $t_*$ and retain degree one by continuity. Their
smooth inverses give the continued height
$\widehat u(t)=z(t,\cdot)\circ\sigma(t,\cdot)^{-1}$.
It is smooth across $t_*$, hence continuous in $H^2$. Since
$\|u(t)\|_{H^2}<r/2$ before $t_*$, its value at $t_*$ has norm
at most $r/2$. The continued graphs therefore stay inside the open
ball of radius $r$ beyond $t_*$, contradicting maximality.
This proves $t_*=\infty$ for the sheet chosen at $t_j$.

Consequently \eqref{miwf:LS-length} proves finite tail dissipation
length and existence of a full $L^2$ limit for $u(t)$. At the
subsequence times $t_j$, the continued graph and the representatives
converging smoothly to $\eta$ differ by orientation-preserving
reparametrizations. Lemma~\ref{miwf:graph-orbit} gives
$u(t_j)\to0$ in $L^2$, identifying the full limit:
\[
 L^2(\mathcal C_\eta;\R)\text{-}\lim_{t\to\infty}u(t)=0.
\]
Letting $t_2\to\infty$ in \eqref{miwf:LS-length} now gives
\[
 \|u(t)\|_2\leq\int_t^\infty\|u_r\|_2\,\dd r
 \leq C\Delta(t)^\theta.
\]
The quantitative estimate \eqref{miwf:graph-energy-control} yields
\[
 \|u(t)\|_{H^2}^2
 \leq C\bigl(\Delta(t)+\Delta(t)^{2\theta}\bigr)
 \leq C\Delta(t)^{2\theta},
\]
since $\Delta(t)\leq1$ and $2\theta\leq1$.

To pass from the scalar height to the parametrized curve, use the
graph-construction operator $\mathcal G$ defined at the start of
this subsection. We prove its precise local Lipschitz estimate:
for every $R>0$ there is $C_{\eta,R}<\infty$ such that
\[
 \|\mathcal G(u)-\mathcal G(w)\|_{W^{2,2}(\mathcal C_\eta;\R^3)}
 \leq C_{\eta,R}\|u-w\|_{H^2(\mathcal C_\eta;\R)}
 \quad\text{if }\|u\|_{H^2},\|w\|_{H^2}\leq R.
\]
Put $h=u-w$, $q_\lambda=w+\lambda h$, and
$b_q=-\sin q\,\eta+\cos q\,\nu_\eta$. Pointwise integration
in $\lambda$ gives
\begin{equation}\label{miwf:graph-integral}
 \mathcal G(u)-\mathcal G(w)
      =\int_0^1 b_{q_\lambda}h\,\dd\lambda.
\end{equation}
For $\|q\|_{H^2}\leq R$, the embedding $H^2\hookrightarrow C^1$
and the fixed smoothness of $\eta,\nu_\eta$ give
\[
 \|b_q\|_\infty+\|\partial_s b_q\|_\infty
       +\|\partial_s^2 b_q\|_2\leq C_{\eta,R}.
\]
Indeed the first derivative uses only $q_s$ and fixed first
derivatives, while the second derivative is a sum of bounded
trigonometric coefficients multiplying $q_{ss}$, $q_s^2$, $q_s$,
or fixed smooth functions. Here
$\|q_s^2\|_2\leq\|q_s\|_\infty\|q_s\|_2$.
The product identities
\[
 (b_qh)_s=(b_q)_sh+b_qh_s,\qquad
 (b_qh)_{ss}=(b_q)_{ss}h+2(b_q)_sh_s+b_qh_{ss}
\]
therefore imply
$\|b_qh\|_{W^{2,2}(\mathcal C_\eta;\R^3)}
 \leq C_{\eta,R}\|h\|_{H^2}$.
The product identities hold for $H^2$ functions by approximation
with smooth functions. Moreover, $H^2$ is a Banach algebra on the
circle, and the sine and cosine maps are continuous on $H^2$.
Since $\lambda\mapsto q_\lambda$ is continuous in $H^2$, so is
$\lambda\mapsto b_{q_\lambda}h$ as a $W^{2,2}$-valued map.
Thus the integral in \eqref{miwf:graph-integral} is a well-defined
Bochner integral in $W^{2,2}$, and the pointwise identity holds in
that space.
Since $\|q_\lambda\|_{H^2}\leq R$, integration proves the
claimed estimate. Taking $w=0$, $u=u(t)$, and
$\psi_t=\sigma_t^{-1}$ gives $\gamma_t\circ\psi_t=\mathcal G(u(t))$
and proves full strong $W^{2,2}$ convergence together with
\eqref{miwf:full-energy-distance}.

We next upgrade this convergence to every $C^k$ norm in the same
continued graph gauge. First, at any times $q_j\to\infty$ for which
$\miwfD(\gamma_{q_j})\to0$, we claim that $u(q_j)\to0$ in $W^{4,2}$.
Write $u_j=u(q_j)$ and denote the curvature of $\Gamma_{u_j}$ by
$\kappa_j$. In the fixed graph coordinate $s$ it is given by
$k_j(s)=\kappa[u_j](s)$. The three estimates \eqref{miwf:first},
\eqref{miwf:max}, and \eqref{miwf:second} bound, respectively,
the first arclength derivative of the curvature in $L^2$, the curvature
itself in $L^\infty$, and its second arclength derivative in $L^2$.
The first and third, together with the length bounds and the
one-dimensional Morrey embedding, give a uniform $C^0$ bound for
$\partial_{s_{\Gamma}}\kappa_j$.
Formula~\eqref{miwf:graphcurvature} takes the form
\[
 (u_j)_{ss}=\frac{v_j^3}{a_j}\bigl(k_j-B(s,u_j,(u_j)_s)\bigr),
 \qquad (k_j)_s=v_j\partial_{s_{\Gamma}}\kappa_j.
\]
The already proved $H^2$ convergence bounds $u_j$ in $C^1$ and
keeps $a_j,v_j$ bounded above and away from zero. Since
$\dd s_{\Gamma_{u_j}}=v_j\dd s$, these speed bounds make the
$L^2$ norms in graph arclength and in the fixed coordinate $s$
uniformly equivalent throughout the following estimates. The first formula
bounds $(u_j)_{ss}$ in $L^\infty$. Since its coefficients depend only
on $s,u_j,(u_j)_s$, their first derivatives are now bounded.
The second identity and the upper bound for $v_j$ give
$\|(k_j)_s\|_{L^\infty}\leq C$. Differentiating the first identity
once now bounds $(u_j)_{sss}$ in $L^\infty$. Moreover,
\[
 (k_j)_{ss}
 =v_j^2\partial_{s_{\Gamma}}^2\kappa_j
   +(v_j)_s\partial_{s_{\Gamma}}\kappa_j
\]
is uniformly bounded in $L^2$. The coefficient derivatives through
order two involve only the derivatives of $u_j$ already bounded.
Differentiating the first formula twice therefore bounds
$(u_j)_{ssss}$ in $L^2$.
The compact embedding $W^{4,2}\Subset C^3$ and the unique $H^2$
limit $0$ now give $u_j\to0$ in $C^3$. The curvature formula and
the definition of graph speed then give $k_j\to\kappa_\eta$
in $C^1$, $v_j\to1$ in $C^2$, and the residual estimate
\eqref{miwf:residual} gives
\[
 (k_j)_{ss}
 =\frac{v_j^2}{2}\bigl(G_j-k_j^3-k_j\bigr)
   +\frac{(v_j)_s}{v_j}(k_j)_s
 \longrightarrow-\frac12(\kappa_\eta^3+\kappa_\eta)
 =\partial_s^2\kappa_\eta\quad\text{in }L^2.
\]
To recover the fourth derivative of the height explicitly, put
$C_j=a_j/v_j^3$ and $B_j=B(s,u_j,(u_j)_s)$. Then
$k_j=C_j(u_j)_{ss}+B_j$, so
\[
 C_j(u_j)_{ssss}
  =(k_j)_{ss}-\partial_s^2B_j
       -2(\partial_sC_j)(u_j)_{sss}
       -(\partial_s^2C_j)(u_j)_{ss}.
\]
The $C^3$ convergence of $u_j$ gives $C_j\to1$ and
$B_j\to\kappa_\eta$ in $C^2$. The last two terms therefore
tend to zero uniformly, while the first two tend to zero together
in $L^2$ by the residual limit just proved. Since $C_j$ is bounded
away from zero, $(u_j)_{ssss}\to0$ in $L^2$. Together with the
already established convergence of the lower derivatives, this
proves $u_j\to0$ in $W^{4,2}$.

Given an arbitrary sequence $t_n\to\infty$, choose
$q_n\in[t_n-1,t_n-1/2]$ with
\[
 \miwfD(\gamma_{q_n})
 \leq2\int_{t_n-1}^{t_n-1/2}\miwfD(\gamma_t)\,\dd t
 \longrightarrow0.
\]
Such a choice follows from continuity and the mean-value bound;
the integral tends to zero by finite total dissipation.
The preceding claim gives $u(q_n)\to0$ in $W^{4,2}$.
Apply the graph construction of Lemma~\ref{miwf:smoothing} on
$[0,1]$ with initial height exactly $u(q_n)$. Its solution and the
continued height have the same initial data in the same parameter
coordinate. The smooth graph uniqueness argument in that lemma gives
agreement first near zero and then on $[0,1]$: at any last time of
agreement their common jets lie in the open coefficient neighborhood
of the constructed solution, so uniqueness applies again. Since $t_n-q_n\in[1/2,1]$, the uniform
positive-time conclusion gives $u(t_n)\to0$ in every $C^k$ norm.
The sequence $t_n$ was arbitrary, proving full smooth convergence
with the same reparametrizations $\psi_t$.

Finally \eqref{miwf:LS-dissipation} and $\Delta'=-\miwfD$ give
$\Delta'\leq-c\Delta^{2(1-\theta)}$.
For $0<\theta<1/2$, put $p=1-2\theta>0$.  As long as $\Delta>0$,
integration gives
\[
 \Delta(t)\leq
 \bigl(\Delta(t_0)^{-p}+cp(t-t_0)\bigr)^{-1/p}.
\]
For $\theta=1/2$ it gives
$\Delta(t)\leq\Delta(t_0)e^{-c(t-t_0)}$.
If $\Delta(t_0)=0$, the trajectory is stationary and equals the
identified limit in graph coordinates.  The preceding bounds thus
hold in all cases.  Together with \eqref{miwf:full-energy-distance}
they prove the asserted rates.
\end{proof}

\begin{remark}[Convergence norms]
The energy estimate gives the quantitative $W^{2,2}$ rate.
Positive-time regularization supplies smooth convergence and
uniform bounds of every fixed spatial order at sufficiently
large times.
\end{remark}

\begin{corollary}[Infinite-time limit of global Hopf-torus trajectories]
\label{miwf:hopf-global}
Let $F_t$ be a smooth global MIWF trajectory starting from a Hopf
immersion in the convention of \eqref{eq:hopf-cover-energy}, with any
finite covering degree and any finite initial energy.  Then its
multiplicity-counting varifolds converge to a unique nonzero integral
two-varifold as $t\to\infty$.  If $\eta$ is the parametrized elastic
limit of the reduced trajectory and $d$ is the additional torus
covering degree, the limit is $d\Var(X_\eta)$.  The marked conformal
classes for a fixed transported marking converge as well.  In the
standard pullback marking, for some $C>0$ and all sufficiently large
$t$,
\begin{equation}\label{eq:global-modulus-rate}
 d_h(\tau_t,\tau_\infty)
 \leq C\bigl(\miwfE(\gamma_t)-\miwfE(\eta)\bigr)^\theta,
\end{equation}
where $\theta$ is the exponent in Theorem~\ref{miwf:full}.
The same bound holds after every fixed finite torus covering.
\end{corollary}
\begin{proof}
Write $F_0=X_{\gamma_0}\circ q$ with covering degree $d$; see the
introduction of Section~\ref{sec:profiles} and
\eqref{eq:hopf-cover-energy}.
A connected finite covering of a torus is regular, since its
fundamental group is abelian.  Every deck transformation preserves
$F_0$, and uniqueness and domain reparametrization invariance of the
MIWF imply that it preserves $F_t$ for every $t$.  Hence $F_t$ descends
through the same covering $q$ to a smooth regular torus immersion
$X_t$. On this quotient torus, the fibre-equivariance argument in
the proof of Corollary~\ref{cor:hopf} applies to $X_t$ at every
energy: uniqueness transports the initial fibre action, and the
rank-two immersion descends under the Hopf map to a regular curve
on the fixed base circle. Removing its tangential velocity on
each compact time interval gives a global normal-gauge profile
$\gamma_t$. By \eqref{eq:sheetwise-hopf-gradient} and
$D\pi(N)=2\nu$ for the horizontal unit normal $N$, its velocity is
$-G\nu/(1+\kappa^2)^2$ on every parameter sheet. Integration with
\eqref{eq:hopf-cover-energy} consequently gives
\[
 \Will(F_t)=d\pi\miwfE(\gamma_t),\qquad
 -W'(t)=d\pi\miwfD(\gamma_t).
\]
Theorem~\ref{miwf:full} and smoothness on compact time intervals imply
\[
 \int_0^\infty\sqrt{-W'(t)}\,\dd t
 =\sqrt{d\pi}\int_0^\infty\sqrt{\miwfD(\gamma_t)}\,\dd t<\infty.
\]
Also $|A_t^0|^2\geq2$.  Proposition~\ref{prop:infinite-finite-length}
therefore gives the unique integral-varifold limit.  It is nonzero:
\eqref{miwf:lengthbound} and the energy bound
$\miwfE(\gamma_t)\leq E_0$ give
$\Area(F_t)=d\pi L(\gamma_t)\geq4d\pi^3/E_0>0$, and convergence of
the weights preserves this lower bound on their total mass.

To identify the limit, represent the reparametrized profiles as the
smoothly convergent graphs over $\eta$ from Theorem~\ref{miwf:full}.
The horizontal-lift construction in the proof of
Theorem~\ref{thm:all-energy-profile} supplies Hopf representatives
converging in $C^1$ on one fixed torus.  Their varifolds converge to
$\Var(X_\eta)$, and the degree-$d$ covering multiplies this by $d$.
Varifolds are invariant under source reparametrization, so uniqueness
identifies the original trajectory's endpoint with $d\Var(X_\eta)$.

By \eqref{eq:marked-speed}, the hyperbolic speed is at most
$\sqrt{E_0}\sqrt{\miwfD}/(2\pi)$.  Finite dissipation length thus
makes the marked conformal path Cauchy in the complete upper
half-plane.  Its limit represents the area class and length of
$\eta$, by smooth profile convergence and continuous dependence of
horizontal holonomy.  The tail estimate
\eqref{miwf:LS-length} gives \eqref{eq:global-modulus-rate}.
The fixed-cover fractional-linear map in the proof of
Theorem~\ref{thm:all-energy-modulus} is a hyperbolic isometry, so it
preserves both convergence and this bound.
\end{proof}

\begin{remark}[Global trajectories and multiplicity]
Theorem~\ref{miwf:full} applies to smooth global reduced
trajectories and retains parametrized multiplicity. In particular,
an $m$-fold equator is an admissible elastic limit; its Hopf
Willmore energy is $2\pi^2m$ before any additional torus covering.
\end{remark}

\section{Conformal endpoints, rank, and circular rigidity}
\label{sec:conformal}
The marked conformal-class map on a fixed oriented torus is
$P(F)=[F^*g_{\Sph^3}]_{\mathcal T_1}$, where
$\mathcal T_1\simeq\HH$ and
$\mathcal M_1=\HH/\mathrm{PSL}_2(\miwfZ)$.
Here two smooth metrics represent the same marked point when they
differ by a positive conformal factor and the pullback by an
orientation-preserving diffeomorphism isotopic to the identity.
A positively oriented ordered homology basis fixes the identification
with $\HH$.
The projection $t\mapsto P(F_t)$ of a smooth MIWF trajectory
$\{F_t\}_{0\leq t<T_{\max}}$ into $\mathcal T_1$ is smooth,
including at rank-deficient slices. Convergence of this projected path
at an endpoint and maximal rank of the differential $DP_{F_t}$ at
every interior time $t<T_{\max}$ are separate properties.
The operator norms for the right inverses below are specified in
\eqref{tp:eq:inverse-bounds} and \eqref{tp:eq:full-Hopf-inverse}.
A profile is called \emph{circular} when its image is the
intersection of $\Sph^2$ with an affine plane; its parametrization
may traverse that circle several times. An $m$-fold equator is a
great circle traversed $m$ times, with $m$ a positive integer.

\subsection{The intrinsic low-energy conformal endpoint}
\begin{corollary}[Unique Pinkall-modulus limit]
\label{cor:pinkall-modulus}
Under the hypotheses of Theorem~\ref{thm:anchored-profile}, orient the
profiles consistently with the reduced flow.  Let
\(\mathcal A_t\in\R/4\pi\mathbb Z\) be the oriented area class
enclosed by \(\gamma_t\), and let \(\mathcal A_T\) be the corresponding 
area class for \(\gamma_T\).  Then
\begin{equation}
 \label{eq:area-class-limit}
 \mathcal A_t\longrightarrow\mathcal A_T
 \qquad\text{in }\R/4\pi\mathbb Z .
\end{equation}
After fixing a real representative at \(t=0\), choose the unique
continuous real lift on \([0,T]\) of this area-class path, and denote
its values again by \(\mathcal A_t\). Then
\(\mathcal A_t\to\mathcal A_T\) in \(\R\) as \(t\nearrow T\).
For these representatives, Pinkall's lattice is
\begin{equation}
 \label{eq:pinkall-lattice}
 \Gamma_t=
 \operatorname{span}_{\mathbb Z}
 \left\{(2\pi,0),
 \left(\frac{\mathcal A_t}{2},\frac{L(\gamma_t)}2\right)\right\}.
\end{equation}
The ordered generators represent the positive fibre cycle and a
cycle projecting to one traversal of the parametrized profile.
We use this ordered homology basis as the \emph{standard Hopf
marking}. Its continuous continuation in time is the transported
marking. The corresponding marked modulus is
\begin{equation}
 \label{eq:pinkall-modulus}
 \tau_t=\frac{\mathcal A_t+iL(\gamma_t)}{4\pi}\in\mathbb H .
\end{equation}
It satisfies
\[
 \tau_t\longrightarrow
 \tau_T:=\frac{\mathcal A_T+iL_T}{4\pi}.
\]
Its imaginary part satisfies the quantitative bounds
\begin{equation}
 \label{eq:modulus-imaginary-bounds}
 \frac{\pi}{8}\leq\operatorname{Im}\tau_T\leq\frac2\pi.
\end{equation}
Consequently the conformal classes
\[
 [\tau_t]\longrightarrow[\tau_T]
 \qquad\text{in }\mathcal M_1=\mathbb H/\mathrm{PSL}_2(\mathbb Z)
\]
have a unique finite-time limit, independent of the approaching time
sequence.  Orient each Hopf torus by the ordered pair consisting of the
positive Hopf-fibre tangent and the horizontal lift of the oriented
profile tangent.  With this orientation, the class \([\tau_T]\) is the
intrinsic oriented conformal class of the limiting \(C^1\) Hopf torus.
In particular, no conformal
degeneration occurs at the finite endpoint.
\end{corollary}

\begin{proof}
We fix the connection and phase conventions explicitly.  Use the
quaternionic Hopf map $\pi(q)=\widetilde q\,q$ of \cite{Pinkall,Jakob},
where $q\mapsto\widetilde q$ is the real-linear involution fixing
$1,j,k$ and sending $i$ to $-i$.
The positive fibre action is $q\mapsto e^{i\theta}q$.  Its connection
one-form is
\[
 \alpha_q(v)=\langle iq,v\rangle_{\R^4},\qquad
 \alpha_q(iq)=1,
 \qquad \mathrm d\alpha=\tfrac12\pi^*\omega,
\]
where $\omega$ is the area form on the unit sphere in
$\operatorname{span}\{1,j,k\}$ with its boundary orientation induced
by the ordered basis $(1,j,k)$. The horizontal space is
$\ker\alpha_q$, the orthogonal complement of the fibre tangent $iq$
in $T_q\Sph^3$. For the curvature identity, at $q=1$
the horizontal vectors $j,k$ satisfy
$\mathrm d\alpha(j,k)=2$ and $(\pi^*\omega)(j,k)=4$.
Right multiplication preserves the connection, acts transitively on
$\Sph^3$, and intertwines $\pi$ with orientation-preserving
isometries of the base sphere. Hence the identity holds everywhere.
Choose an open set $U\subset\Sph^2$ admitting a smooth local
section $\sigma:U\to\Sph^3$,
$\pi\circ\sigma=\mathrm{id}_U$, and set
$b=\sigma^*\alpha\in\Omega^1(U)$. The associated trivialization is
\[
 \Phi:U\times\Sph^1\longrightarrow\pi^{-1}(U),\qquad
 \Phi(x,e^{i\theta})=e^{i\theta}\sigma(x).
\]
Invariance of $\alpha$ under the fibre action and $\alpha(iq)=1$ give
\[
 \Phi^*\alpha=\operatorname{pr}_{\Sph^1}^*(\mathrm d\theta)
                +\operatorname{pr}_U^*b,\qquad
 \mathrm db=\sigma^*(\mathrm d\alpha)=\tfrac12\omega|_U.
\]
Here $\operatorname{pr}_U$ and $\operatorname{pr}_{\Sph^1}$ are
the two projections from $U\times\Sph^1$, and $\mathrm d\theta$
is the angular one-form on $\Sph^1$;
the last identity is on the base $U$, whereas
$\mathrm d\alpha=\pi^*\omega/2$ is on $\Sph^3$.
In particular $b$ is not closed. For a lift written locally as
$\eta(s)=e^{i\theta(s)}\sigma(\gamma(s))$, the same formula gives
\[
 \alpha_{\eta(s)}(\eta'(s))
 =\theta'(s)+b_{\gamma(s)}(\gamma'(s)).
\]
Thus horizontality is exactly the scalar equation
$\theta'=-b(\gamma')$.

For a smooth closed profile choose a smooth spanning map
$u:D^2\to\Sph^2$ whose oriented boundary is $\gamma$.
Here $D^2$ is the closed unit disc with its standard orientation.
Such a map exists since $\Sph^2$ is simply connected. The pullback
of the Hopf bundle to the contractible disc is trivial. Choose a
smooth section, represented by a lift $\widehat u:D^2\to\Sph^3$
with $\pi\circ\widehat u=u$, and put
$\beta=\widehat u^*\alpha$. Then
$\mathrm d\beta=u^*\omega/2$. Parametrize the oriented boundary
by $c:[0,1]\to\partial D^2$, with $c(0)=c(1)$, and write its
horizontal lift as
$\eta(s)=e^{i\theta(s)}\widehat u(c(s))$, choosing a real angle on
$[0,1]$. Horizontality gives $\theta'=-\beta(c')$; no equality of
the endpoint angles is imposed. Integration and Stokes' theorem give
the terminal transport phase
\[
 h_\gamma
 =\exp\!\left(-i\int_{\partial D^2}\beta\right)
 =\exp\!\left(-\frac i2\int_{D^2}u^*\omega\right).
\]
Two spanning maps, glued with opposite orientations along the
boundary, give a sphere map of integer degree. Their area integrals
therefore differ by $4\pi$ times that degree. Consequently
$\mathcal A_\gamma=[\int_{D^2}u^*\omega]\in\R/4\pi\mathbb Z$
is well defined and
\[
 h_\gamma=\exp\!\left(-\frac{i\mathcal A_\gamma}{2}\right).
\]
The angle in the positive fibre coordinate needed to close the lift is
$\delta_\gamma=\mathcal A_\gamma/2\pmod{2\pi}$.
Thus the terminal lift point is $h_\gamma$ times the initial point,
and multiplication by $h_\gamma^{-1}$ closes the lift.

Horizontal transport depends continuously on a path in the $C^1$
topology. To see this, subdivide the parameter interval into finitely
many closed subintervals whose limiting path segments lie in section
domains. Every sufficiently $C^1$-close path uses the same domains.
The coefficients $b(\gamma')$ then converge uniformly on each
subinterval, as do the section transition functions at the subdivision
points. Integration of the scalar equation, with initial lift points
chosen continuously, proves convergence of the transport phase.
The holonomy identity also holds for the regular $C^1$ limiting
profile: approximate it in $C^1$ by smooth closed curves. The
pointwise short-geodesic homotopy between each approximation and
the limiting curve has area tending to zero, since its transverse
derivative tends uniformly to zero and its tangential derivative
stays bounded. Stokes' theorem therefore gives convergence of their
oriented area classes, while the transport phases converge by the
preceding argument. Consequently
\eqref{eq:profile-strong} gives
\[
 \exp\!\left(\frac{i\mathcal A_t}{2}\right)
 \longrightarrow
 \exp\!\left(\frac{i\mathcal A_T}{2}\right).
\]
The map $[\mathcal A]\mapsto\exp(i\mathcal A/2)$ is a homeomorphism from
$\R/4\pi\mathbb Z$ onto the unit circle.  This proves
\eqref{eq:area-class-limit}. The area-class path, already continuous
for $t<T$, thus extends continuously to all of $[0,T]$.
Since $[0,T]$ is simply connected, the covering
$p:\R\to\R/4\pi\mathbb Z$ has a unique continuous lift of this
path after one fixes its real value at $t=0$.
Using the notation $\mathcal A_t$ for the lifted values gives
$\mathcal A_t\to\mathcal A_T$ in $\R$ as $t\nearrow T$.

To identify the intrinsic modulus of the regular $C^1$ limiting
profile, let $\eta_T(y)$ be a horizontal
lift satisfying $\pi(\eta_T(y))=\gamma_T(y/P_T)$, where
$P_T=L_T/2$ and $\gamma_T$ is defined on $\R/\miwfZ$.
The horizontal dilation of $\pi$ is $2$, so $|\eta_T'(y)|=1$.
Fix its initial point and
put $\delta_T=\mathcal A_T/2\in\R$, the real representative of
the closing-angle class modulo $2\pi$ determined by the chosen
area lift. Uniqueness of
horizontal lifting and its equivariance under the fibre action give
\[
 \eta_T(y+P_T)=e^{-i\delta_T}\eta_T(y)
 \qquad\text{for every }y\in\R.
\]
In the ordered Euclidean coordinates $(\theta,y)$ the map
\[
 X_T(\theta,y)=e^{i\theta}\eta_T(y)
\]
is a $C^1$ local isometry: its two coordinate derivatives are
orthogonal unit vectors in Euclidean $\R^4$.
Since the profile traverses its embedded
closed curve once, two parameter pairs have the same image precisely when
their difference belongs to
\[
 \Gamma_T=\operatorname{span}_{\mathbb Z}
 \{(2\pi,0),(\delta_T,P_T)\}.
\]
Indeed, equality of their projections first forces the difference in
$y$ to be an integer multiple of $P_T$, after which the lift identity
determines the difference in $\theta$ modulo $2\pi$.
This proves that $X_T$ descends to an isometry from
$\R^2/\Gamma_T$ onto the limiting Hopf torus, with the orientation
specified in the statement.  The identical calculation for each
smooth profile gives \eqref{eq:pinkall-lattice}, in agreement with
\cite[Proposition~1]{Pinkall}.

Identifying $(\theta,y)$ with $\theta+iy$ and dividing the second
lattice generator by the first gives
$\tau_t=(\mathcal A_t+iL(\gamma_t))/(4\pi)$ and
$\tau_T=(\mathcal A_T+iL_T)/(4\pi)$.
Changing $\mathcal A$ by $4\pi$ changes $\tau$ by $1$ and hence
changes the chosen marking while preserving its class in
$\mathcal M_1$. The initial real area representative fixes this
integer ambiguity throughout the time interval. The convergence of
the lifted area and of the length proves the stated modulus convergence.
The bound $L_T\in[\pi^2/2,8]$ from Theorem~\ref{thm:anchored-profile}
gives \eqref{eq:modulus-imaginary-bounds}. The covering map $X_T$
constructed above induces an isometry from $\R^2/\Gamma_T$ onto the
limiting $C^1$ Hopf torus. Hence $[\tau_T]$ is the intrinsic oriented
conformal class of its induced metric.
\end{proof}

\subsection{Marked conformal estimates without a rank hypothesis}

The finite-time marked-path estimates follow directly from the
reduced energy identity and the length bounds in
Section~\ref{sec:profiles}. At infinite time,
Corollaries~\ref{miwf:hopf-global} and
\ref{cor:global-coordinate-tails} combine these estimates with
the finite dissipation length supplied by Theorem~\ref{miwf:full}.

\begin{theorem}[Dissipation controls the marked conformal path]
\label{thm:all-energy-modulus}
Let $\gamma_t$, $0\leq t<T\leq\infty$, be a smooth reduced MIWF of
regular closed profiles, counting multiplicity. Choose one initial
real lift $\mathcal A_0$ of its spherical area class and continue it
along the time interval. Write $\mathcal A_t$ for this lift,
$L_t=L(\gamma_t)$, $E_t=\miwfE(\gamma_t)$, $D_t=\miwfD(\gamma_t)$,
and $\tau_t=(\mathcal A_t+iL_t)/(4\pi)$ for the standard marked
Hopf modulus. With $h=(dx^2+dy^2)/y^2$ on $\HH$, one has
\begin{equation}\label{eq:marked-speed}
 |\dot\tau_t|_h^2\leq\frac{D_t}{L_t}
                  \leq\frac{E_0}{4\pi^2}D_t,
 \qquad E_0=\miwfE(\gamma_0).
\end{equation}
For $0\leq s<t<T$,
\begin{equation}\label{eq:marked-length}
 \operatorname{Length}_h(\tau|_{[s,t]})
 \leq\frac{\sqrt{E_0}}{2\pi}\sqrt{(t-s)(E_s-E_t)}.
\end{equation}
The length satisfies the sharper scalar bounds
\begin{align}
 |L_t'|^2
 &\leq D_t\int_{\gamma_t}
       \frac{\kappa_t^2}{(1+\kappa_t^2)^2}\dd s_{\gamma_t}
 \leq\min\{L_t/4,E_t-L_t\}\,D_t,
 \label{eq:sharp-length-speed}\\
 |\sqrt{L_t}-\sqrt{L_s}|
 &\leq\frac14\int_s^t\sqrt{D_r}\dd r
 \leq\frac14\sqrt{(t-s)(E_s-E_t)}.
 \label{eq:sqrt-length-increment}
\end{align}
The coefficient $1/4$ multiplying $L_tD_t$ in
\eqref{eq:sharp-length-speed} is optimal among universal constants.

If $T<\infty$, the path has a unique endpoint $\tau_T\in\HH$,
$(\mathcal A,L)\in W^{1,2}((0,T);\R^2)$, and
\begin{equation}\label{eq:marked-tail}
 d_h(\tau_t,\tau_T)
 \leq\frac{\sqrt{E_0}}{2\pi}\sqrt{(T-t)(E_t-E_T)},\qquad
 \frac{\pi}{E_0}\leq\operatorname{Im}\tau_T\leq\frac{E_0}{4\pi},
\end{equation}
where $E_T=\lim_{t\uparrow T}E_t$. Moreover,
\begin{equation}\label{eq:sqrt-length-tail}
 |\sqrt{L_t}-\sqrt{L_T}|
 \leq\frac14\sqrt{(T-t)(E_t-E_T)},\qquad
 L_T=4\pi\operatorname{Im}\tau_T.
\end{equation}
The limit is the marked conformal type of the parametrized endpoint
in Theorem~\ref{thm:all-energy-profile}, with its transported marking.
Every fixed finite torus cover has the corresponding marked endpoint
and the same hyperbolic length and distance estimates.
\end{theorem}
\begin{proof}
The continuous area-class path has a unique real lift after fixing
$\mathcal A_0$, because $[0,T)$ is simply connected, also when
$T=\infty$. In the setting of Corollary~\ref{cor:pinkall-modulus},
this is the restriction of its lift on $[0,T]$ with the same initial
value. For the normal speed $f_t$ the oriented cylinder
swept by the profile gives
$\mathcal A_t'=-\int f_t\dd s_{\gamma_t}$, and length variation gives
$L_t'=-\int\kappa_t f_t\dd s_{\gamma_t}$.
Indeed the unit-sphere area form satisfies $\omega(f_t\nu_t,T_t)=-f_t$.
The cylinder formula also shows that the continued lift is smooth.
For every $(z_1,z_2)\in\Sph^1$, weighted Cauchy--Schwarz yields
\begin{align*}
 |z_1\mathcal A_t'+z_2L_t'|^2
 &=\left|\int_{\gamma_t}
   \frac{G(\gamma_t)(z_1+z_2\kappa_t)}{(1+\kappa_t^2)^2}
          \dd s_{\gamma_t}\right|^2\\
 &\leq D_t\int_{\gamma_t}
          \frac{(z_1+z_2\kappa_t)^2}{(1+\kappa_t^2)^2}
          \dd s_{\gamma_t}\\
 &\leq D_t\int_{\gamma_t}\frac{1}{1+\kappa_t^2}
          \dd s_{\gamma_t}\leq L_tD_t.
\end{align*}
The penultimate inequality uses
$(z_1+z_2\kappa_t)^2\leq1+\kappa_t^2$.
Taking the supremum over all $(z_1,z_2)\in\Sph^1$ proves
$(\mathcal A_t')^2+(L_t')^2\leq L_tD_t$.
Since $|\dot\tau_t|_h=\sqrt{(\mathcal A_t')^2+(L_t')^2}/L_t$,
Lemma~\ref{miwf:length} proves \eqref{eq:marked-speed}.
Integration and $E_t'=-D_t$ give \eqref{eq:marked-length}.

Applying weighted Cauchy--Schwarz to the length component alone gives
the first inequality in \eqref{eq:sharp-length-speed}.
For every real $k$,
\[
 \frac{k^2}{(1+k^2)^2}\leq\frac14,
 \qquad \frac{k^2}{(1+k^2)^2}\leq k^2.
\]
Integrating these two pointwise bounds and using
$\int\kappa_t^2\dd s_{\gamma_t}=E_t-L_t$ proves its second inequality.
Since $L_t>0$, division by $4L_t$ gives
$|\partial_t\sqrt{L_t}|\leq\sqrt{D_t}/4$.
Integration and time Cauchy--Schwarz prove
\eqref{eq:sqrt-length-increment}.
For sharpness, at a covered latitude with constant geodesic curvature
$\kappa=1$, one has $G=2$, $D=L$, and $L'=L/2$.
Such a curve is admissible smooth initial data for the local reduced
flow, so $|L'|^2=LD/4$ at its initial time.

If $T<\infty$, \eqref{eq:marked-length} bounds total hyperbolic length.
For $s<t<T$, it also gives
\[
 d_h(\tau_s,\tau_t)
 \leq\frac{\sqrt{E_0}}{2\pi}\sqrt{(T-s)(E_s-E_T)}.
\]
This tends to zero as $s\uparrow T$, uniformly for $s<t<T$.
Completeness of $\HH$ therefore gives one common endpoint
$\tau_T\in\HH$. Letting the upper time tend to $T$ in
\eqref{eq:marked-length} proves \eqref{eq:marked-tail}; its
imaginary-part bounds follow from $4\pi^2/E_0\leq L_t\leq E_0$.
The same limit in \eqref{eq:sqrt-length-increment} proves
\eqref{eq:sqrt-length-tail}. Also
\[
 \int_0^T\bigl((\mathcal A_t')^2+(L_t')^2\bigr)\dd t
 \leq E_0(E_0-E_T).
\]
Both coordinates are bounded by their finite limits and continuity
on compact time intervals, giving the stated $W^{1,2}$ membership.
Tangential velocity can first be removed by a source diffeomorphism;
this does not change the area lift, length, or marked modulus.
The $C^1$ endpoint-profile convergence in
Theorem~\ref{thm:all-energy-profile} gives convergence of horizontal
holonomy, as in Corollary~\ref{cor:pinkall-modulus}. Thus
$\mathcal A_T$ represents the area class of $\gamma_T$.
For an immersed or multiply traversed profile, use its
orientation-preserving constant-speed representative on
$\R/\miwfZ$ and perform the same construction on $P_\gamma$,
defined in Section~\ref{sec:profiles}.
If $P=L/2$ and $\eta$ is the horizontal lift, its isometric
parametrization on the covering plane is
\[
 (\theta,y)\longmapsto
 \bigl(y/P\bmod\miwfZ,\ e^{i\theta}\eta(y)\bigr)\in P_\gamma.
\]
Equality in this parameter bundle forces the first coordinates
to agree, so the period lattice is generated by
$(2\pi,0)$ and $(\mathcal A/2,P)$ even when the image in
$\Sph^3$ has self-intersections. Its marked modulus is therefore
$(\mathcal A+iL)/(4\pi)$. Applied to the $C^1$ endpoint, this
identifies the limit with $\tau_T$.

For a fixed orientation-preserving finite covering of the parameter
torus, choose a positively oriented basis of its finite-index
sublattice. After division by $2\pi$, this basis is
$a+b\tau$, $c+e\tau$, with fixed integers and
$d:=ae-bc>0$ equal to the covering degree. Its modulus is
\begin{equation}\label{eq:cover-modulus-map}
 \tau_{\rm cov}=\mathfrak c(\tau):=\frac{c+e\tau}{a+b\tau}.
\end{equation}
This formula applies to the pullback of every marked conformal
structure on the base torus. Indeed, its normalized holomorphic
one-form has periods $1,\tau$ on the ordered base cycles. Its pullback
has periods $a+b\tau,c+e\tau$ on the ordered covering cycles.
Dividing by the first period gives \eqref{eq:cover-modulus-map}.
Here $a+b\tau\ne0$ since $ae-bc=d>0$ and $\Im\tau>0$. One has
\[
 \operatorname{Im}\mathfrak c(\tau)
   =\frac{d\operatorname{Im}\tau}{|a+b\tau|^2},\qquad
 \mathfrak c'(\tau)=\frac{d}{(a+b\tau)^2}.
\]
Consequently $\mathfrak c$ is a hyperbolic isometry of $\HH$.
The homology classes of the covering cycles are fixed during the
evolution, so their integer coefficients $a,b,c,e$ remain constant.
The same isometry $\mathfrak c$ consequently applies at every time,
proving the covering assertions.
\end{proof}

For a Hopf lift with additional torus covering degree $d$, the energy
normalization in Lemma~\ref{tp:lem:Hopf-reconstruction} reads
$W(t)=d\pi E_t$ and $-W'(t)=d\pi D_t$. In particular,
\begin{equation}\label{eq:length-surface-normalization}
 |L_t'|\leq\frac{\sqrt{L_t}}{2\sqrt{d\pi}}\sqrt{-W'(t)},\qquad
 |\sqrt{L_t}-\sqrt{L_s}|
 \leq\frac{\sqrt{(t-s)(W(s)-W(t))}}{4\sqrt{d\pi}}.
\end{equation}
These are estimates for profile length; an additional torus covering
does not multiply that length. It multiplies surface area and energy.

\begin{corollary}[Global coordinate variation and geodesic length rate]
\label{cor:global-coordinate-tails}
Let $\gamma_t$ be a smooth global reduced trajectory, and let $\eta$
and $\theta$ be supplied by Theorem~\ref{miwf:full}. Put
$E_t=\miwfE(\gamma_t)$, $D_t=\miwfD(\gamma_t)$, and
$\Delta(t)=E_t-\miwfE(\eta)$. Use the standard Hopf marking and a
continued real area representative to define
$\tau_t=(\mathcal A_t+iL(\gamma_t))/(4\pi)$.
There are constants $C,c>0$ and
$t_0\geq0$ such that the following assertions hold.
For every transported area lift there
is a real limit $\mathcal A_\infty$, and, with $L_\infty=L(\eta)$,
\begin{equation}\label{eq:global-coordinate-tail}
 \int_t^\infty\bigl(|\mathcal A_r'|+|L_r'|\bigr)\dd r
 \leq C\Delta(t)^\theta
 \qquad(t\geq t_0).
\end{equation}
Both derivatives belong to $L^1(0,\infty)\cap L^2(0,\infty)$, and
\begin{equation}\label{eq:global-speed-decay}
 D_t\longrightarrow0,\qquad\mathcal A_t'\longrightarrow0,\qquad
 L_t'\longrightarrow0,\qquad |\dot\tau_t|_h\longrightarrow0.
\end{equation}
If $\eta$ is an $m$-fold equator, then additionally
\begin{equation}\label{eq:geodesic-length-rate}
 |L_t-2\pi m|\leq C\Delta(t)^{2\theta}
 \qquad(t\geq t_0).
\end{equation}
Thus this length difference is bounded by
$C(1+t-t_0)^{-2\theta/(1-2\theta)}$ when $\theta<1/2$ and by
$Ce^{-c(t-t_0)}$ when $\theta=1/2$, after adjusting the constants.
\end{corollary}
\begin{proof}
Theorem~\ref{miwf:full} and its tail estimate
\eqref{miwf:LS-length} give
$\int_t^\infty\sqrt{D_r}\dd r\leq C\Delta(t)^\theta$ for large $t$.
The proof of Theorem~\ref{thm:all-energy-modulus} gives
$|\mathcal A_t'|\leq\sqrt{E_0D_t}$ and
$|L_t'|\leq\sqrt{E_0D_t}/2$. These imply
\eqref{eq:global-coordinate-tail} and finite total variation of both
coordinates. Their limits therefore exist as $t\to\infty$, and
smooth profile convergence identifies
$L_\infty=L(\eta)$ and the area class of $\mathcal A_\infty$.
Integrating $(\mathcal A_t')^2+(L_t')^2\leq E_0D_t$ over
$(0,\infty)$ proves the $L^2$ assertion.

The smooth convergence modulo reparametrization in
Theorem~\ref{miwf:full} makes $G(\gamma_t)$ converge uniformly to
$G(\eta)=0$ in the convergent graph coordinates. The arclength
densities converge uniformly and the lengths stay bounded.
Reparametrization invariance of $D_t$ therefore gives $D_t\to0$.
The estimates of Theorem~\ref{thm:all-energy-modulus} now imply
all the remaining limits in \eqref{eq:global-speed-decay}.

Suppose finally that $\eta$ is an $m$-fold equator, parametrized by
arclength on $\R/(2\pi m)\miwfZ$. The continued normal graph in
the proof of Theorem~\ref{miwf:full} has
$\|u(t)\|_{H^2}\leq C\Delta(t)^\theta$ and $\kappa_\eta=0$.
Its exact length is
\[
 L_t=\int_0^{2\pi m}\sqrt{\cos^2u+u_s^2}\dd s.
\]
For $u$ small in $C^1$, the speed is bounded away from zero and
\[
 \left|\sqrt{\cos^2u+u_s^2}-1\right|
 =\frac{|u_s^2-\sin^2u|}{\sqrt{\cos^2u+u_s^2}+1}
 \leq u_s^2+u^2.
\]
Integration proves
$|L_t-2\pi m|\leq\|u(t)\|_{H^1}^2\leq C\Delta(t)^{2\theta}$.
The energy rates in Theorem~\ref{miwf:full} prove the final assertions.
\end{proof}

The geodesic length rate \eqref{eq:geodesic-length-rate} applies
when the limiting curve is an $m$-fold equator. The loop number $m$
is part of that limiting curve and is not determined by the rate.

For example, the scalar functions $\ell(t)=4\pi+e^{-t}$ and
$e(t)=4\pi+2e^{-t}$ satisfy
\[
 |\ell'(t)|^2\leq
 \min\{\ell(t)/4,e(t)-\ell(t)\}\,(-e'(t)),
\]
while $\ell(t)\to4\pi$. This illustrates a limitation of the scalar
inequality; it is not a construction of a reduced-flow trajectory.

\subsection{A sufficient endpoint criterion}
Use the complete hyperbolic metric $h=(dx^2+dy^2)/y^2$ on $\HH$.
\begin{proposition}[Finite projected length]\label{tp:prop:projected-length}
Let $\tau:[0,T)\to\HH$ be $C^1$, with $T$ finite or infinite.
If $\int_0^T|\dot\tau(t)|_h\dd t<\infty$, then $\tau$ has a
unique limit in $\HH$, and its moduli class has the corresponding
limit in $\mathcal M_1$. No rank hypothesis is required.
\end{proposition}
\begin{proof}
For $0\leq s<t<T$, the distance is bounded by the remaining path length:
\[
 d_h(\tau(s),\tau(t))
 \leq\int_s^t|\dot\tau(r)|_h\dd r
 \leq\int_s^T|\dot\tau(r)|_h\dd r.
\]
The tail tends to zero as $s\to T$, also for $T=\infty$.
The path is therefore Cauchy in the complete space $\HH$ and has a
unique limit. Continuity of the quotient map gives its moduli limit.
\end{proof}

\subsection{A differential criterion}
For the metric $g=F^*g_{\Sph^3}$, let $\TT(g)$ denote the space of
smooth symmetric, trace-free, divergence-free covariant two-tensors.
On a torus its real dimension is two. It represents the cotangent
space of marked Teichm\"uller space by the metric pairing.

\begin{proposition}[Differential of the marked conformal class]\label{tp:prop:rank-general}
For a smooth immersion $F:\Sigma\to\Sph^3$ of the fixed oriented
torus, choose a unit normal $N$ and write $A_F^0$ for the scalar
trace-free second fundamental form relative to $N$. For every
$f\in C^\infty(\Sigma;\R)$, the normal variation $fN$ induces the
cotangent pairing
\begin{equation}\label{tp:eq:rank-pairing}
 \langle DP_F(fN),q\rangle
 =-2\int_\Sigma f\langle A_F^0,q\rangle_g\dd\mu_g,
 \qquad q\in\TT(g),
\end{equation}
where the cotangent identification is chosen by this metric pairing.
The differential has rank less than two if and only if a nonzero
$q\in\TT(g)$ satisfies $\langle A_F^0,q\rangle_g=0$ everywhere.
\end{proposition}
\begin{proof}
The first variation of the induced metric is $\dot g=-2fA$.
Its trace part pairs trivially with $q$. An infinitesimal source
reparametrization gives a Lie derivative of $g$, whose pairing also
vanishes by integration by parts and $\operatorname{div}_g q=0$.
The identification of the differential of metric projection with these
pairings is the Teichm\"uller metric-variation formula
\cite[Section~3, equations~(3.2)--(3.4)]{KuwertSchaetzleConformal}.
Those equations apply to smooth closed marked surfaces. Regard
$F$ as an immersion into $\R^4$; it has the same induced metric.
The Euclidean and spherical trace-free second fundamental forms
agree, since the radial second fundamental form of $\Sph^3$
is pure trace. The Euclidean metric-variation formula therefore
gives exactly \eqref{tp:eq:rank-pairing} for the spherical normal
variation $fN$.

Failure of surjectivity into a two-dimensional space is equivalent
to the existence of a nonzero cotangent annihilator. By the displayed
formula and arbitrary smooth test functions $f$, this is precisely
the pointwise condition in the statement.
\end{proof}

\subsection{The isothermic condition}
In this paper an immersion is called \emph{globally isothermic} if
there is a nonzero $q\in\TT(g)$ satisfying
$\langle A_F^0,q\rangle_g=0$ everywhere; it is called
\emph{nonisothermic} if no such tensor exists. In conformal
coordinates, a transverse-traceless tensor is the real part of a
holomorphic quadratic differential. Thus the condition requires
one globally defined differential.
Proposition~\ref{tp:prop:rank-general} identifies it with failure
of surjectivity of $DP_F$. The implication to local conformal
principal-curvature coordinates away from the zeros of that
differential, and the bound of one on the rank deficiency, are
proved in \cite[Proposition~4.1]{KuwertSchaetzleConformal}, applied
to the same immersion in $\R^4$. Radial variations alter the metric
only by a conformal factor, so they add no Teichm\"uller direction.
The condition is invariant under fixed ambient conformal maps.
For Hopf immersions, $|A^0|^2=2(1+\kappa^2)>0$ everywhere;
the next proposition characterizes their isothermic condition.

\subsection{Rank loss and its quantitative size}
\begin{proposition}[Hopf rank criterion]\label{tp:prop:Hopf-rank}
For a smooth parametrized Hopf immersion, including finite covers,
$\rank DP_F=1$ if $\kappa$ is constant and $\rank DP_F=2$ otherwise.
For the profile moment operator
\[
 R_\gamma h=\left(\int\kappa h\dd s,\int h\dd s\right),
 \qquad R_\gamma:L^2(\dd s)\to\R^2,
\]
the Gram matrix and its determinant are
\begin{equation}\label{tp:eq:Gram}
 M_\gamma=R_\gamma R_\gamma^*
 =\begin{pmatrix}\int\kappa^2\dd s&\int\kappa\dd s\\
                 \int\kappa\dd s&L\end{pmatrix},\qquad
 \det M_\gamma=L\int(\kappa-\bar\kappa)^2\dd s,
 \quad \bar\kappa=L^{-1}\int\kappa\dd s.
\end{equation}
All operator norms of $R_\gamma$ and its right inverse use the
arclength $L^2$ norm and the Euclidean norm on $\R^2$.
For every profile put
$\sigma_\gamma:=\int(\kappa-\bar\kappa)^2\dd s$; thus
$\sigma_\gamma=0$ exactly when $\kappa$ is constant.
For nonconstant $\kappa$ and every $z\in\R^2$, there is a unique
$h_z\in(\ker R_\gamma)^\perp$ satisfying $R_\gamma h_z=z$.
It has the smallest $L^2(\dd s)$ norm among all preimages of $z$.
The least-norm right inverse is therefore
$R_\gamma^\dagger:=R_\gamma^*M_\gamma^{-1}$, with
$h_z=R_\gamma^\dagger z$; in particular, $h_0=0$. Its operator norm satisfies
\begin{equation}\label{tp:eq:inverse-bounds}
 \|R_\gamma^\dagger\|^2
 =\frac{\miwfE+\sqrt{\miwfE^2-4L\sigma_\gamma}}{2L\sigma_\gamma},\qquad
 \frac{\miwfE}{2L\sigma_\gamma}\leq\|R_\gamma^\dagger\|^2
                         \leq\frac{\miwfE}{L\sigma_\gamma}.
\end{equation}
Thus $\miwfE\leq E_0$, $L\geq\ell_0>0$, and
$\sigma_\gamma\geq\sigma_0>0$ imply
$\|R_\gamma^\dagger\|\leq(E_0/(\ell_0\sigma_0))^{1/2}$.
\end{proposition}
\begin{proof}
On the flat Hopf torus every transverse-traceless tensor has constant
coefficients in the parallel frame and has the form
$q=\left(\begin{smallmatrix}\alpha&\beta\\\beta&-\alpha\end{smallmatrix}\right)$.
Indeed, the complex coefficient lifts to a doubly periodic entire
function on $\mathbb C$. It is bounded on a compact fundamental
parallelogram and hence on $\mathbb C$, so it is constant by
Liouville's theorem.
Formula \eqref{tp:eq:Hopf-geometry} gives
$\langle A^0,q\rangle=2(\alpha\kappa+\beta)$.
Thus a nonzero annihilator exists exactly when $\kappa$ is constant;
in that case its dimension is one. This proves the rank statement
using Proposition~\ref{tp:prop:rank-general}. On a finite cover the
same parallel frame and holomorphic-constancy argument apply, so
the conclusion is unchanged.

The $L^2$ adjoint of $R_\gamma$ is
$R_\gamma^*(\alpha,\beta)=\alpha\kappa+\beta$.
Matrix multiplication and expansion about the mean give
\eqref{tp:eq:Gram}. Strict positivity in the nonconstant case proves
the right-inverse formula: $R_\gamma R_\gamma^*M_\gamma^{-1}=I$,
and $\operatorname{Ran}R_\gamma^*=(\ker R_\gamma)^\perp$.
Every preimage of $z$ has the form $h_z+k$ with
$k\in\ker R_\gamma$, and orthogonality gives
\[
 \|h_z+k\|_{L^2(\dd s)}^2
 =\|h_z\|_{L^2(\dd s)}^2+\|k\|_{L^2(\dd s)}^2.
\]
Thus $h_z$ is the unique preimage of smallest norm, also when
$z=0$. Since $R_\gamma R_\gamma^*=M_\gamma$, one has
\[
 \|h_z\|_{L^2(\dd s)}^2
 =z^TM_\gamma^{-1}R_\gamma R_\gamma^*M_\gamma^{-1}z
 =z^TM_\gamma^{-1}z.
\]
The norm of this least-norm preimage is therefore
$\|h_z\|_{L^2(\dd s)}=\sqrt{z^TM_\gamma^{-1}z}$ for every
$z\in\R^2$. The maximum of its square over $|z|=1$ is
$\lambda_{\min}(M_\gamma)^{-1}$. Since
$\operatorname{tr}M_\gamma=\miwfE$ and $\det M_\gamma=L\sigma_\gamma$,
the quadratic eigenvalue formula gives \eqref{tp:eq:inverse-bounds};
its numerator lies between $\miwfE$ and $2\miwfE$. Substituting the three
stated bounds gives the uniform estimate.
\end{proof}

\begin{corollary}[Asymptotic degeneration]\label{tp:cor:asymptotic-rank}
If a sequence of noncircular constant-speed profiles converges in
$C^2$, with positive limiting length, to a covered circle, then
the least-norm right inverses in Proposition~\ref{tp:prop:Hopf-rank}
have norms tending to infinity. In particular this holds at a
double-equator limit.
\end{corollary}
\begin{proof}
Write $L_j=L(\gamma_j)$ and $E_j=\miwfE(\gamma_j)$.
Curvature converges uniformly to a constant and $L_j\to L_\infty>0$.
Consequently $\det M_{\gamma_j}=L_j\sigma_{\gamma_j}\to0$, while
$E_j\to E_\infty>0$. Now apply \eqref{tp:eq:inverse-bounds}:
\[
 \|R_{\gamma_j}^\dagger\|^2
 \geq\frac{E_j}{2\det M_{\gamma_j}}\longrightarrow\infty.
\]
\end{proof}
The moment calculation also determines the full normal differential in
specified surface norms. For the standard pullback torus, extend $DP_F$
to scalar normal variations $fN$ with $f\in L^2(\dd\mu_F)$, and
identify $T_{P(F)}\mathcal T_1$ with $\mathbb C\simeq\R^2$ using
the standard marked coordinate $\tau$. If $\kappa$
is nonconstant, its least-norm right inverse $(DP_F)^\dagger$ satisfies
\begin{equation}\label{tp:eq:full-Hopf-inverse}
 \|(DP_F)^\dagger\|_{\mathrm{Eucl}}=2\pi\sqrt\pi\,
                                  \|R_\gamma^\dagger\|,\qquad
 \|(DP_F)^\dagger\|_h=\frac{L\sqrt\pi}{2}
                                  \|R_\gamma^\dagger\|.
\end{equation}
For either target norm, $(DP_F)^\dagger z$ is the preimage with
smallest surface $L^2$ norm. The first operator norm in
\eqref{tp:eq:full-Hopf-inverse} uses the Euclidean target norm;
the second uses $|z|_h=|z|/\operatorname{Im}\tau$. Both use the
unweighted surface $L^2$ norm for the preimage.
To verify this, average a normal scalar speed $v$ over each Hopf fibre
and denote its average by $\bar v(s)$. The cotangent formula
\eqref{tp:eq:rank-pairing} shows that $DP$ annihilates $v-\bar v$,
since the coefficients of $A^0$ and every $q\in\TT$ are fibre
invariant. Averaging is an orthogonal projection, so a least-norm
preimage is fibre invariant. For smooth $\bar v$, vary the profile
with normal velocity $h\nu$, where $h=2\bar v$. Its Hopf lift on the
fixed torus has normal component $\bar v N$ because $D\pi(N)=2\nu$,
as in Lemma~\ref{tp:lem:Hopf-reconstruction}. The tangential part
pairs to zero in \eqref{tp:eq:rank-pairing}. Integrating the local
area element $\frac12\dd s\,\dd\varphi$ over a fibre gives
\[
 \|\bar v\|_{L^2(\dd\mu_F)}^2
       =\pi\int|\bar v|^2\dd s
       =\frac\pi4\|h\|_{L^2(\dd s)}^2.
\]
The area and length variations, proved by the cylinder calculation
in Theorem~\ref{thm:all-energy-modulus}, are
$(-\int h\dd s,-\int\kappa h\dd s)$.
For a prescribed modulus variation $z\in\mathbb C$, the profile
constraint is $R_\gamma h=-4\pi(\Im z,\Re z)$. Its least-norm
solution is explicitly
\[
 h=-4\pi R_\gamma^\dagger(\Im z,\Re z)
   =-4\pi R_\gamma^*M_\gamma^{-1}(\Im z,\Re z).
\]
The preceding surface/profile
norm conversion multiplies its norm by $\sqrt\pi/2$, giving the
factor $2\pi\sqrt\pi$ in the Euclidean target norm. Since a
hyperbolic unit vector has Euclidean norm $\Im\tau=L/(4\pi)$,
the second formula follows. The pairing formula also extends the
differential to a bounded operator on the surface $L^2$ space. By
Proposition~\ref{tp:prop:Hopf-rank}, $R_\gamma^*(a,b)=a\kappa+b$
and the unique least-norm profile preimage belongs to
$\operatorname{Ran}R_\gamma^*=\operatorname{span}\{1,\kappa\}$.
Consequently each least-norm surface preimage is smooth and is a
fibre-invariant linear combination of $1$ and the lifted curvature $\kappa$.
Thus Corollary~\ref{tp:cor:asymptotic-rank} is also a statement about
the full normal differential in these specified norms. For uniformly equivalent domain and target norms, the divergence
conclusion persists: for $T=DP_F$ the quantity
$\sup_{z\ne0}\inf_{Tf=z}\|f\|/\|z\|$ changes by at most the uniform
norm-equivalence factors. In the specified Hilbert norms this quantity
equals $\|T^\dagger\|$.

For a noncircular profile, let $q$ be any fixed additional torus
covering of degree $d$, and use the transported marked coordinate
$\tau_{\rm cov}=\mathfrak c(\tau)$ from
\eqref{eq:cover-modulus-map}. Then
\begin{equation}\label{eq:covered-inverse-norm}
 \|(DP_{F\circ q})^\dagger\|_h
 =\sqrt d\,\|(DP_F)^\dagger\|_h
 =\frac{L\sqrt{d\pi}}2\|R_\gamma^\dagger\|.
\end{equation}
To prove this, write the base and covering tori as
$\Sigma=\R^2/\Lambda$ and $\widehat\Sigma=\R^2/\widehat\Lambda$,
where $\widehat\Lambda\subset\Lambda$ has index $d$ and
$q([x]_{\widehat\Lambda})=[x]_\Lambda$. The deck group of $q$ is
$G=\Lambda/\widehat\Lambda$, of order $d$; its action on the covering
torus $\widehat\Sigma$ is
\[
 [\lambda]\cdot[x]_{\widehat\Lambda}=[x+\lambda]_{\widehat\Lambda}.
\]
Since $q([\lambda]\cdot\widehat x)=q(\widehat x)$, this action is
isometric with respect to the metric
$(F\circ q)^*g_{\Sph^3}=q^*(F^*g_{\Sph^3})$ on $\widehat\Sigma$. Averaging scalar
normal variations over $G$ is therefore the orthogonal projection onto
the deck-invariant functions, which are precisely the pullbacks of
functions on $\Sigma$.

The tensor $A^0_{F\circ q}$ is pulled back from the base, and every
TT tensor on the flat covering torus has constant coefficients in the
pulled-back parallel frame. Thus the integrand coefficients in
\eqref{tp:eq:rank-pairing} are deck-invariant, and the differential
annihilates the orthogonal complement of the invariant functions.
Consequently a preimage of least $L^2$ norm is deck-invariant and
descends. For every smooth variation $F_\varepsilon$ of the base
immersion, $(F_\varepsilon\circ q)^*g_{\Sph^3}
=q^*(F_\varepsilon^*g_{\Sph^3})$. The period calculation in
\eqref{eq:cover-modulus-map} therefore gives
$P(F_\varepsilon\circ q)=\mathfrak c(P(F_\varepsilon))$
also for variations through general immersions.
Differentiating at a variation with velocity $fN$ gives, with
$\widehat F=F\circ q$,
\[
 DP_{\widehat F}\bigl((fN)\circ q\bigr)
 =D\mathfrak c_{P(F)}\bigl(DP_F(fN)\bigr),\qquad
 \|f\circ q\|_{L^2(\dd\mu_{\widehat F})}
 =\sqrt d\,\|f\|_{L^2(\dd\mu_F)}.
\]
By continuity these identities also hold for $f\in L^2(\dd\mu_F)$.
Set $C:=D\mathfrak c_{P(F)}$. For a prescribed target vector
$\widehat z\in T_{P(\widehat F)}\HH$, define
$z:=C^{-1}\widehat z\in T_{P(F)}\HH$ and the smooth function
$f_z$ on $\Sigma$ by $f_zN:=(DP_F)^\dagger z$.
For every $f\in L^2(\dd\mu_F)$, the first identity gives
\[
 DP_{\widehat F}\bigl((fN)\circ q\bigr)=\widehat z
 \quad\Longleftrightarrow\quad DP_F(fN)=z.
\]
The second identity multiplies the norm of every such base preimage
$fN$ by the same factor $\sqrt d$ upon pullback. Since every
least-norm preimage on the covering torus is deck-invariant, these
identities and uniqueness of the minimizers yield
\[
 (DP_{\widehat F})^\dagger\widehat z=(f_zN)\circ q,
\]
and hence the pointwise norm identity
\[
 \bigl\|(DP_{\widehat F})^\dagger\widehat z\bigr\|_{L^2(\dd\mu_{\widehat F})}
 =\sqrt d\,\bigl\|(DP_F)^\dagger(C^{-1}\widehat z)\bigr\|_{L^2(\dd\mu_F)}.
\]
As $C$ is an isometry for the hyperbolic tangent norms,
$C^{-1}$ maps the unit circle at $P(\widehat F)$ onto the unit circle
at $P(F)$. Taking the supremum over unit target vectors
$\widehat z$ on both sides of the last identity proves the first
equality in \eqref{eq:covered-inverse-norm}; the second follows from
\eqref{tp:eq:full-Hopf-inverse}.

On an energy window $\miwfE\leq E_0$, the length lower bound and
\eqref{tp:eq:inverse-bounds} also show that
$\sigma_\gamma\geq\sigma_0>0$ implies
$\|R_\gamma^\dagger\|\leq E_0/(2\pi\sqrt{\sigma_0})$.
The marked-path estimates of Theorem~\ref{thm:all-energy-modulus}
hold independently of this curvature-fluctuation bound.

\begin{samepage}
\subsection{Circular profiles and Pinkall's semicircles}
\begin{proposition}[Exact circular trajectories]\label{tp:prop:circles}
If a smooth reduced trajectory is circular at a time
$t_*\in[0,T_{\max})$, it has a circular continuation for all
$t\geq t_*$. Its covering multiplicity, say
$m\in\mathbb N$ with $m\geq1$, remains constant for $t\geq t_*$.
After a fixed rotation of the sphere and a choice of orientation,
write this continuation, up to reparametrization, as
\[
 \gamma_t(\theta)=(\sin r(t)\cos(m\theta),
 \sin r(t)\sin(m\theta),\cos r(t)),\qquad 0<r(t)<\pi.
\]
Then $K(t)=\cot r(t)$ satisfies
\begin{equation}\label{tp:eq:circle-ODE}
 K'=-K,\qquad K(t)=K(t_*)e^{-(t-t_*)},\qquad
 r'=\sin r\cos r.
\end{equation}
For the area lift $\mathcal A=2\pi m(1-\cos r)$ and $L=2\pi m\sin r$,
\begin{equation}\label{tp:eq:semicircle}
 (\mathcal A/2-m\pi)^2+(L/2)^2=(m\pi)^2.
\end{equation}
The circular continuation exists for all future times and converges
to the $m$-fold equator, with
\begin{equation}\label{tp:eq:circle-limit}
 (\mathcal A/2,L/2)\to(m\pi,m\pi),\quad
 \tau\to\tfrac m2(1+i),\quad
 \miwfE\to2\pi m,\quad \Will(X_{\gamma_t})\to2\pi^2m.
\end{equation}
\end{proposition}
\end{samepage}
\begin{proof}
First, a regular closed curve of constant geodesic curvature $K$
is a covered circle. Indeed, by \eqref{miwf:frame},
$(K\gamma+\nu)_s=0$. The curve therefore lies in the plane
$(K\gamma+\nu)\cdot\gamma=K$; its regular closed parametrization
is a finite cover of the resulting geometric circle.

For the displayed latitude parametrization, $K=\cot r$ and
$\nu=-\partial_r\gamma$. Its reduced speed is
$f=-K/(1+K^2)$, so \eqref{miwf:flow} becomes
$r'=K/(1+K^2)=\sin r\cos r$. Differentiating $\cot r$ gives
$K'=-K$. Conversely these scalar equations construct a smooth
regular circular solution for every future time. The uniqueness
argument at the end of the proof of Lemma~\ref{miwf:smoothing}
uses only a smooth regular reference curve and positivity of the
graph principal coefficient, so it also applies to this circular
reference without requiring it to be elastic. Together with the
ordinary differential equation removing tangential velocity, it
identifies the constructed solution, up to parametrization, with
the given solution after $t_*$; compare the short-time construction
in \cite[Proposition~4.5]{Jakob} for local existence. The constructed
solution uses the same parametrized $m$-fold cover at every time,
so uniqueness preserves this integer on $[t_*,\infty)$.

The area and length formulas follow by integration over a spherical
cap and by traversing its boundary $m$ times. Substitution proves
\eqref{tp:eq:semicircle}. The exponential formula gives $r\to\pi/2$,
also when $K(t_*)=0$, and then all limits follow from
\eqref{eq:hopf-cover-energy} and \eqref{eq:pinkall-modulus}.
\end{proof}

For $m=1$, \eqref{tp:eq:semicircle} is precisely the equality case of
Pinkall's spherical isoperimetric inequality
$L^2\geq4\pi \mathcal A-\mathcal A^2$ \cite[p.~382, equations~(14)--(15)]{Pinkall}.
For $m>1$, the chosen cap-area representative and the length are
$m$ times those of one traversal, which gives
\eqref{tp:eq:semicircle} directly. Other area lifts
translate the semicircle horizontally by $2\pi h$, $h\in\miwfZ$.
The displayed semicircles use the orientation and area representative
specified in Proposition~\ref{tp:prop:circles}.

\begin{figure}[htbp]
\centering
\begin{tikzpicture}
\begin{axis}[width=.94\linewidth,height=6.4cm,xmin=0,xmax=8.2,
 ymin=0,ymax=4.4,xlabel={$\mathcal A/(2\pi)$},ylabel={$L/(2\pi)$},
 axis equal image,axis lines=left,grid=major,legend style={at={(.98,.98)},anchor=north east,
 font=\small},samples=151,domain=0:180]
\addplot[blue,thick] ({1-cos(x)},{sin(x)});\addlegendentry{$m=1$}
\addplot[teal,thick] ({2*(1-cos(x))},{2*sin(x)});\addlegendentry{$m=2$}
\addplot[orange,thick] ({3*(1-cos(x))},{3*sin(x)});\addlegendentry{$m=3$}
\addplot[red,thick] ({4*(1-cos(x))},{4*sin(x)});\addlegendentry{$m=4$}
\addplot[only marks,mark=*,black] coordinates {(1,1)(2,2)(3,3)(4,4)};
\end{axis}
\end{tikzpicture}
\caption{Circular profiles of multiplicity $m$ for the specified area lifts.
Each dot marks the corresponding equatorial endpoint.}
\end{figure}
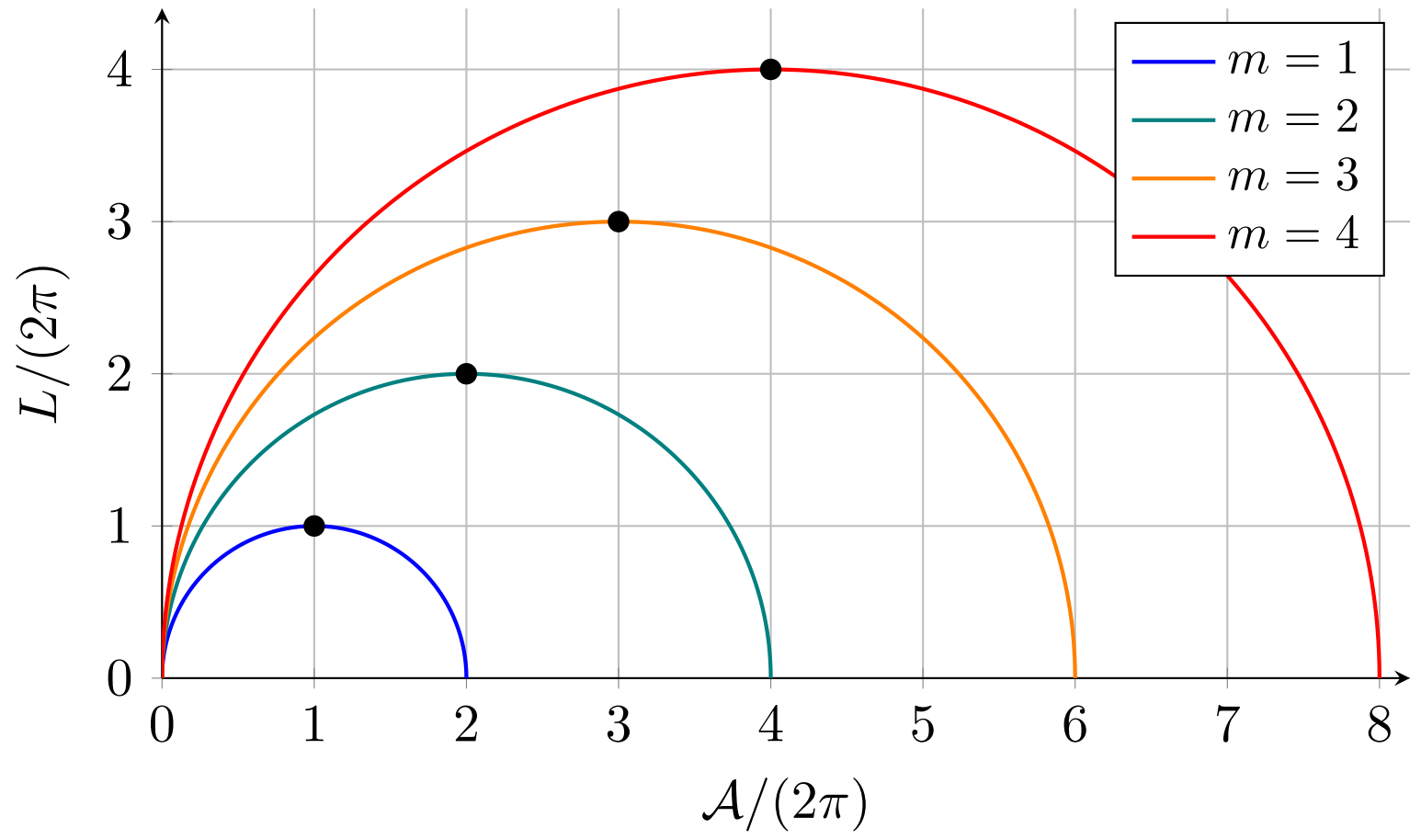

\begin{remark}[Circular solutions for the classical Willmore flow]
For the unweighted curve flow $f=-c(2\kappa_{ss}+\kappa^3+\kappa)$
with $c>0$, circular profiles instead satisfy
$K'=-cK(1+K^2)^2$, by substituting constant curvature into
\eqref{tp:eq:curve-variation}. For the classical surface equation
$\partial_tF=-\nabla_{L^2}\Will(F)$, the Hopf gradient relation in
\eqref{eq:sheetwise-hopf-gradient} and the horizontal dilation
\eqref{eq:hopf-horizontal-dilation} give $f=-4G_\gamma$, so $c=4$. For the reduced MIWF the weight
$(1+K^2)^{-2}$ instead gives $K'=-K$, as in
\eqref{tp:eq:circle-ODE}.
\end{remark}

\subsection{A finite-time circular-rigidity theorem}
Positive-time analyticity yields finite-time circular rigidity.

\begin{lemma}[Analyticity of the scaled mean-curvature variance]\label{tp:lem:analytic-diagnostic}
For any smooth normal MIWF of a closed oriented torus in $\Sph^3$
on a smooth existence interval $[0,T)$, choose a smooth unit normal
$N_t$ compatible with the orientation and write
$\mathsf H_t=\langle H_t,N_t\rangle$ for its signed scalar trace
mean curvature. The function
\begin{equation}\label{tp:eq:CMC-diagnostic}
 \mathcal Q(t)=\Area(F_t)\int_\Sigma \mathsf H_t^2\dd\mu_t
                 -\left(\int_\Sigma \mathsf H_t\dd\mu_t\right)^2
\end{equation}
is real analytic on $(0,T)$. It is nonnegative and vanishes exactly
when the slice has constant mean curvature (CMC).
\end{lemma}
\begin{proof}
Writing $\overline{\mathsf H}_t
=\Area(F_t)^{-1}\int_\Sigma\mathsf H_t\dd\mu_t$, one has
\[
 \mathcal Q(t)=\Area(F_t)\int_\Sigma
       (\mathsf H_t-\overline{\mathsf H}_t)^2\dd\mu_t.
\]
This proves nonnegativity and the equality criterion.
The functional is unchanged by source reparametrization and by
reversal of the chosen unit normal.

Fix an interior time $t_0$. The compact image of $F_{t_0}$ misses
some point of $\Sph^3$, since a smooth surface has zero
three-dimensional measure. Smooth time continuity therefore gives
one stereographic projection valid on a time neighborhood of $t_0$.
By fixed conformal covariance its Euclidean representatives form a
smooth umbilic-free MIWF on this neighborhood.

Endow the parameter torus with a complex structure and fix a
compatible real analytic reference immersion into $\R^3$, for
example a standard torus of revolution after identifying the domain
with $\R^2/(2\pi\miwfZ)^2$. Fix $p>3$ and use the Euclidean
slice at $t_0$ as the center $U_0$ in
\cite[Theorem~2.2(1)]{JakobRegularity}. It is smooth and
umbilic-free. That theorem, with the reference immersion fixed,
supplies one positive existence time for all initial immersions in
a $W^{4-4/p,p}$ neighborhood of $U_0$. Smooth time continuity puts
every sufficiently nearby earlier slice in this neighborhood.
Choose such a time $a<t_0$ with $t_0-a$ smaller than half this
existence time. The modified solution restarted at $a$ therefore
exists past $t_0$. Restrict to an open interval $J$ containing
$t_0$ inside both its positive-time existence interval and the
chosen stereographic time neighborhood. By
\cite[Theorem~2.3(ii)]{JakobRegularity}, this DeTurck-reparametrized
family is jointly real analytic in space and time on $J$; its initial
immersion at $a$ need only be smooth. The DeTurck correspondence and
smooth uniqueness identify it, up to source reparametrization, with
the original Euclidean flow on $J$.

Inverse stereographic projection is real analytic. In oriented local
coordinates, the spherical unit normal is obtained by normalizing
the vector of signed $3\times3$ minors of the matrix whose rows are
the position vector and the two coordinate tangents in $\R^4$.
Its sign is chosen to agree with the orientation. Its denominator
is positive, so it is analytic as well.
Thus the spherical metric, signed scalar mean curvature, and area
density of the resulting representative are analytic. On compact
subintervals of $J$, finitely many spatial charts give a common
radius and bound for the time Taylor expansions. Termwise
integration over the compact parameter torus shows that each of
the three integrals in
\eqref{tp:eq:CMC-diagnostic}, namely area, total signed mean curvature,
and total squared mean curvature, is real analytic on $J$.
Transport the chosen normal under the source diffeomorphisms.
Change of variables then identifies all three integrals with those
for $F_t$, so $\mathcal Q$ is real analytic near $t_0$.
Since $t_0$ was arbitrary, $\mathcal Q$ is
analytic on $(0,T)$.
\end{proof}

\begin{theorem}[Preservation of nonconstant profile curvature]\label{tp:thm:no-circle}
Let $\gamma_t$ be a smooth reduced MIWF with a smooth regular closed
initial profile. If $\gamma_0$ has nonconstant geodesic curvature,
then $\gamma_t$ has nonconstant geodesic curvature for every
$0<t<T_{\max}(\gamma_0)$. The conformal-class differential of its
Hopf torus consequently has rank two at every such time.
\end{theorem}
\begin{proof}
By Lemma~\ref{tp:lem:Hopf-reconstruction}, the given reduced trajectory
has a smooth normal Hopf representative on a fixed torus. The
functional $\mathcal Q$ is analytic by
Lemma~\ref{tp:lem:analytic-diagnostic}.
In the standard pullback model,
\eqref{tp:eq:Hopf-geometry}--\eqref{eq:hopf-cover-energy} give the exact identity
\begin{equation}\label{tp:eq:Q-Gram}
 \mathcal Q(t)=4\pi^2\left(L\int\kappa^2\dd s-
                              \left(\int\kappa\dd s\right)^2\right)
     =4\pi^2\det M_{\gamma_t}.
\end{equation}
An additional torus covering of degree $d$ multiplies $\mathcal Q$ by $d^2$.
Consequently $\mathcal Q(t)=0$ exactly when $\kappa_t$ is constant.

If $\kappa_{t_*}$ were constant at an interior time, Proposition
\ref{tp:prop:circles} would give $\mathcal Q(t)=0$ on a nonempty open interval
immediately after $t_*$. The identity theorem for real analytic
functions gives $\mathcal Q\equiv0$ on the connected interval
$(0,T_{\max})$. Smooth continuity at time zero then makes $\mathsf H_0$,
and hence $\kappa_0$, constant, contrary to the hypothesis.
Proposition~\ref{tp:prop:Hopf-rank} proves the last assertion.
\end{proof}

\begin{remark}[Asymptotic circular limits]\label{tp:rem:rigidity-limit}
For an initially noncircular trajectory converging smoothly to a
covered circle, Corollary~\ref{tp:cor:asymptotic-rank} gives
divergence of the least-norm right inverses as $t\to\infty$.
If this trajectory is global, Theorem~\ref{miwf:full} makes the
limit elastic. Its constant curvature $K$ then satisfies
$K(1+K^2)=0$, so that the limit is a covered equator.
The circular trajectories of Proposition~\ref{tp:prop:circles}
converge to an equator with the same covering multiplicity.
\end{remark}

\section{Obstructions to autonomous Teichm\"uller projection}\label{tp:sec:nonprojectability}

\subsection{The projected velocity retains profile information}
Let
\[
 \mathcal X(F)=-|A_F^0|^{-4}\nabla_{L^2}\Will(F)
\]
denote the normal MIWF velocity. An autonomous projected equation
on a collection of immersions would be a single vector field
$b:U\to\mathbb C$, with $U\subset\HH$ containing their marked
points, such that $DP_F\mathcal X(F)=b(P(F))$ for every immersion
in that collection. Tangent vectors in the marked coordinate are
identified with complex numbers by their two real components.
For equatorial graphs, write $e_3=(0,0,1)$.
By the area and length variations in the proof of Theorem
\ref{thm:all-energy-modulus}, along any
smooth reduced trajectory,
\begin{equation}\label{tp:eq:projected-velocity}
 \mathcal A'=\int\frac{G_\gamma}{(1+\kappa^2)^2}\dd s,\qquad
 L'=\int\frac{\kappa G_\gamma}{(1+\kappa^2)^2}\dd s.
\end{equation}
Their combination $(\mathcal A'+iL')/(4\pi)$ is
$DP_{X_\gamma}\mathcal X(X_\gamma)$ in the standard Hopf marking.
Throughout this section, primes on $\mathcal A$ and $L$ denote
derivatives with respect to flow time at the stated initial profile.
Amplitude derivatives are written explicitly as $\partial_\epsilon$.
Each smooth regular initial profile below generates a local reduced
flow by \cite[Proposition~4.5]{Jakob}. Lemma~\ref{tp:lem:Hopf-reconstruction}
lifts it to a normal surface flow on the marked parameter torus.
For the finite-dimensional families below, an $O(\epsilon^r)$
remainder in a fixed spatial $C^j$ norm is bounded by
$C_j|\epsilon|^r$ on a sufficiently small common amplitude interval;
$C_j$ may depend on the family and on $j$.

\begin{theorem}[Failure of autonomous projection]\label{tp:thm:nonprojectable}
There are pairs of smooth embedded noncircular profiles whose standard
Hopf tori have Willmore energy below $8\pi$, the same marked
modulus in \eqref{eq:pinkall-modulus}, and different instantaneous marked
modulus velocities under MIWF. Each immersion has rank two for its
conformal-class differential. Consequently, the Teichm\"uller
projection does not transform the MIWF into a common autonomous
first-order equation $\dot\tau=b(\tau)$ valid for all such Hopf
trajectories: their projected velocities are not determined by the
marked conformal type alone.
\end{theorem}
\begin{proof}
On the once-traversed equator $C(\theta)=(\cos\theta,\sin\theta,0)$,
$0\leq\theta\leq2\pi$, let
\[
 u_{n,\epsilon}(\theta)=\epsilon\cos(n\theta),\qquad
 \gamma_{n,\epsilon}=\cos u_{n,\epsilon}\,C
                         +\sin u_{n,\epsilon}\,e_3,
 \qquad n=2,3.
\]
For sufficiently small $|\epsilon|$ these are embedded regular
graphs. In latitude-longitude coordinates, the component toward
the north pole has oriented area
\[
 \mathcal A(u)=\int_0^{2\pi}(1-\sin u(\theta))\dd\theta\in(0,4\pi).
\]
We use this real representative for all equatorial graphs below.
It gives exactly
\[
 \mathcal A_{n,\epsilon}=\int_0^{2\pi}(1-\sin u_{n,\epsilon})\dd\theta=2\pi,
\]
because a shift of $\theta$ by $\pi/n$ changes the sign of $u$.
Their length is an even analytic function of $\epsilon$ and has expansion
\begin{equation}\label{tp:eq:length-expansion}
 L_{n,\epsilon}
 =\int_0^{2\pi}\sqrt{\cos^2u_{n,\epsilon}
                            +(u_{n,\epsilon})_\theta^2}\dd\theta
 =2\pi+c_n\epsilon^2+O(\epsilon^4),\qquad
 c_n=\tfrac\pi2(n^2-1)>0.
\end{equation}
Indeed the quadratic integrand is
$\tfrac12(u_\theta^2-u^2)$ and the cosine and sine squares have
integral $\pi$. Write $z=\epsilon^2$. The inverse function theorem
applied to this analytic length function of $z$ gives, for every
sufficiently small common length excess $\delta>0$, a unique
positive $\epsilon_n(\delta)$ with
$L_{n,\epsilon_n}=2\pi+\delta$ and
$\epsilon_n^2=\delta/c_n+O(\delta^2)$.

We compute the length velocity in the fixed graph gauge
$\Gamma(u)=\cos u\,C+\sin u\,e_3$. The reference curve is the
equator, so $\kappa_*=0$ and its arclength parameter is $\theta$.
Thus the quantities in \eqref{miwf:graphcurvature} specialize to
$a=\cos u$, $c=\sin u$, $p=u_\theta$, and
$v=(\cos^2u+u_\theta^2)^{1/2}$. In particular,
\[
 \kappa[u]
 =\frac{\cos u\,u_{\theta\theta}
          +\sin u\,(\cos^2u+2u_\theta^2)}{v^3},
 \qquad \partial_{s_\Gamma}=v^{-1}\partial_\theta.
\]
Let $V[u]$ denote the scalar reduced normal speed from
\eqref{miwf:flow}, with respect to the oriented unit normal of the
graph. The normal component of $\partial_t\Gamma(u)$ is
$(\cos u/v)u_t$, so the graph-height evolution equation is explicitly
\[
 V[u]=-\frac{2(v^{-1}\partial_\theta)^2\kappa[u]
                    +\kappa[u]^3+\kappa[u]}{(1+\kappa[u]^2)^2},
 \qquad u_t=\mathcal F(u):=\frac{v}{\cos u}\,V[u].
\]
Here $(v^{-1}\partial_\theta)^2$ means composition of the
arclength derivative with itself, including differentiation of
$v^{-1}$.

The displayed formulas show that $\mathcal F$ is a real analytic
function of $(u,u_\theta,u_{\theta\theta},u_{\theta\theta\theta},
u_{\theta\theta\theta\theta})$ near zero. Under $u\mapsto-u$,
$v$ and $\cos u$ are unchanged, whereas $\kappa[u]$ changes sign.
The arclength derivative is therefore unchanged and $V[u]$ changes
sign. Consequently $\mathcal F(-u)=-\mathcal F(u)$: the quadratic
Taylor term of the graph-height velocity $\mathcal F$ at zero
vanishes. Since $D\kappa(0)[w]=w_{\theta\theta}+w$ and $V[0]=0$,
\[
 D\mathcal F(0)[w]=-(2\partial_\theta^2+1)
                         (\partial_\theta^2+1)w=-\mathcal Hw,
 \qquad \mathcal H=2\partial_\theta^4+3\partial_\theta^2+1.
\]
For $u=u_{n,\epsilon}$, Taylor expansion now gives
\[
 \left.u_t\right|_{t=0}=\mathcal F(u_{n,\epsilon})
       =-\mathcal H u_{n,\epsilon}+O(\epsilon^3),\qquad
 \mathcal H\cos(n\theta)=\lambda_n\cos(n\theta),
 \quad \lambda_n=2n^4-3n^2+1.
\]
The remainder holds in every fixed spatial $C^j$ norm, uniformly
for $n=2,3$ on a sufficiently small common amplitude interval,
because these are fixed smooth analytic families and the
denominators above stay bounded away from zero.
Differentiation of the length functional gives
\[
 DL(u)[w]=\int_0^{2\pi}
 \frac{u_\theta w_\theta-\sin u\cos u\,w}
      {\sqrt{\cos^2u+u_\theta^2}}\dd\theta
 =\int_0^{2\pi}(u_\theta w_\theta-uw)\dd\theta
       +O(\|u\|_{C^1}^3\|w\|_{C^1}).
\]
The remainder estimate is uniform on a fixed sufficiently small
$C^1$ ball: the numerator errors are cubic in $u,u_\theta$, and
the denominator is $1+O(\|u\|_{C^1}^2)$ and bounded away from
zero. Substituting the initial velocity gives
\[
 L'_{n,\epsilon}=-2c_n\lambda_n\epsilon^2+O(\epsilon^4).
\]
Here $L'_{n,\epsilon}$ is the time derivative of length along the
flow starting at $\gamma_{n,\epsilon}$. For the matched profiles,
$\delta=L_{n,\epsilon_n(\delta)}-2\pi$ and
$\epsilon_n(\delta)^2=\delta/c_n+O(\delta^2)$, so
$\epsilon_n(\delta)^4=O(\delta^2)$. Substitution gives
$L'_{n,\epsilon_n(\delta)}=-2\lambda_n\delta+O(\delta^2)$.
With $\lambda_2=21$ and $\lambda_3=136$, this yields
\begin{equation}\label{tp:eq:different-speeds}
 L'_{2,\epsilon_2(\delta)}=-42\delta+O(\delta^2),\qquad
 L'_{3,\epsilon_3(\delta)}=-272\delta+O(\delta^2).
\end{equation}
The two velocities differ for all sufficiently small $\delta>0$.
The profiles have exactly the same pair $(\mathcal A,L)=(2\pi,2\pi+\delta)$,
so they have the same marked $\tau$. However, their imaginary modulus
velocities, $L'/(4\pi)$, differ by \eqref{tp:eq:different-speeds}.

Their linearized curvatures are $(1-n^2)\epsilon\cos(n\theta)$,
so both have nonconstant curvature for sufficiently small nonzero
amplitude. Proposition~\ref{tp:prop:Hopf-rank} gives rank two.
Finally $\miwfE(\gamma_{n,\epsilon})\to2\pi$ and
$\Will(X_{\gamma_{n,\epsilon}})\to2\pi^2<8\pi$.
Continuity makes both Hopf Willmore energies strictly below $8\pi$ for sufficiently
small $\delta$. A projected autonomous vector field would assign
one velocity to their common marked conformal class in $\mathcal T_1$,
contrary to the above computation.
\end{proof}

\subsection{Noncollinearity and recalibration of time}
Two velocities at the same marked point are called noncollinear
when they are linearly independent over $\R$ in
$T_\tau\HH\simeq\R^2$. Positive changes of time multiply a velocity
by a positive scalar. We construct two noncollinear velocities by
computing the cubic area term for a mixed Fourier profile.

\begin{lemma}[Equatorial graphs and integrated velocities]
\label{tp:lem:graph-velocities}
Let $C(\theta)=(\cos\theta,\sin\theta,0)$ and
$\gamma_u=\cos u\,C+\sin u\,e_3$, where $u$ is smooth,
$2\pi$-periodic, and $|u|<\pi/2$. Use the real area
representative $\mathcal A(u)$ defined above, and put $p=u_\theta$,
$a=\cos u$, $c=\sin u$, and $v=(a^2+p^2)^{1/2}$.
Then
\begin{equation}\label{tp:eq:equatorial-curvature}
 \dd s=v\dd\theta,\qquad
 \kappa=\frac{a u_{\theta\theta}+c(a^2+2p^2)}{v^3},\qquad
 \int\kappa\dd s=\int_0^{2\pi}\sin u\dd\theta=2\pi-\mathcal A.
\end{equation}
At any smooth closed profile the instantaneous reduced MIWF
velocities satisfy
\begin{align}
 \mathcal A'&=\int\left(\frac{\kappa}{1+\kappa^2}
        +\frac{8\kappa\kappa_s^2}{(1+\kappa^2)^3}\right)\dd s,
        \label{tp:eq:area-velocity-integrated}\\
 L'&=\int\left(\frac{\kappa^2}{1+\kappa^2}
        +\frac{2(3\kappa^2-1)\kappa_s^2}{(1+\kappa^2)^3}\right)\dd s.
        \label{tp:eq:length-velocity-integrated}
\end{align}
For an equatorial graph with $\mathcal A(u)=2\pi$, the
area-velocity integrand in \eqref{tp:eq:area-velocity-integrated}
can equivalently be replaced by
$-\kappa^3/(1+\kappa^2)+8\kappa\kappa_s^2/(1+\kappa^2)^3$.
\end{lemma}
\begin{proof}
Writing $n_u=-\sin u\,C+\cos u\,e_3$ gives
$\partial_\theta\gamma_u=a C_\theta+p n_u$ and unit normal
$(a n_u-p C_\theta)/v$. Differentiating once more and pairing
with this normal gives the curvature formula. Moreover
\[
 \partial_\theta\arctan\frac{p}{a}
      =\frac{a u_{\theta\theta}+c p^2}{a^2+p^2},\qquad
 \kappa v=\partial_\theta\arctan\frac{p}{a}+\sin u.
\]
The real arctangent is periodic since $a>0$. Its derivative
integrates to zero, proving \eqref{tp:eq:equatorial-curvature}.
The area formula is the cap-area formula already used above.

In \eqref{tp:eq:projected-velocity}, integrate the $2\kappa_{ss}$
term once by parts. Differentiating $(1+\kappa^2)^{-2}$ yields
$-4\kappa\kappa_s(1+\kappa^2)^{-3}$, proving the area identity.
For the length identity differentiate
$\kappa(1+\kappa^2)^{-2}$, whose derivative is
$(1-3\kappa^2)\kappa_s(1+\kappa^2)^{-3}$.
There are no boundary terms. Finally use $\int\kappa\dd s=0$
and $\kappa/(1+\kappa^2)=\kappa-\kappa^3/(1+\kappa^2)$.
\end{proof}

\begin{theorem}[Noncollinearity at a fixed marked conformal type]
\label{tp:thm:noncollinear}
For all sufficiently small $\delta>0$, two embedded noncircular
profiles with $(\mathcal A,L)=(2\pi,2\pi+\delta)$ can be chosen so that
their instantaneous reduced MIWF velocities obey
\begin{align}
 \mathcal A'_1&=0,& L'_1&=-42\delta+O(\delta^2),\label{tp:eq:pure-velocity}\\
 \mathcal A'_2&=-\frac{465}{2\sqrt\pi}\delta^{3/2}+O(\delta^{5/2}),&
 L'_2&=-782\delta+O(\delta^2).\label{tp:eq:mixed-velocity}
\end{align}
The subscript $i=1,2$ labels the two initial profiles, and
$\mathcal A_i',L_i'$ are their flow-time derivatives at time zero.
Their standard Hopf immersions have Willmore energy below $8\pi$
and rank two for $DP$, and converge smoothly to the Clifford torus
as $\delta\downarrow0$ after smooth identifications of the domains.
In particular, their projected velocities are noncollinear.
\end{theorem}
\begin{proof}
Take the first profile to be $\gamma_{2,\epsilon_2(\delta)}$
from Theorem~\ref{tp:thm:nonprojectable}. Its curvature changes sign
under $\theta\mapsto\theta+\pi/2$, while its arclength density
is unchanged. Both $\kappa_{ss}$ and $G_\gamma$ change sign as
well. Formula \eqref{tp:eq:projected-velocity} therefore gives
$\mathcal A'_1=0$ exactly. The matched length velocity in
\eqref{tp:eq:different-speeds} gives the second assertion in
\eqref{tp:eq:pure-velocity}.

For the second profile put $w=\cos2\theta+\cos4\theta$ and
choose $u_\epsilon=c(\epsilon)+\epsilon w$. The equation
\[
 \int_0^{2\pi}\sin(c+\epsilon w)\dd\theta=0
\]
has a unique small real analytic solution $c(\epsilon)$ with
$c(0)=0$, since its partial derivative with respect to $c$ at the
origin has the value $2\pi$.
Uniqueness of this solution and oddness of the sine give
$c(-\epsilon)=-c(\epsilon)$; its linear term vanishes because
$\int w=0$. Taylor expansion, using $\int w^3=3\pi/2$, gives
$c(\epsilon)=\epsilon^3/8+O(\epsilon^5)$.
By \eqref{tp:eq:equatorial-curvature}, this correction yields
$\mathcal A=2\pi$ exactly.

The length is even and analytic in $\epsilon$. Its quadratic
term is $\tfrac12\int(w_\theta^2-w^2)\dd\theta=9\pi$, so
\begin{equation}\label{tp:eq:mixed-length}
 L(\gamma_{u_\epsilon})=2\pi+9\pi\epsilon^2+O(\epsilon^4).
\end{equation}
Inverting its analytic expansion in $z=\epsilon^2$ supplies a
positive amplitude for each sufficiently small $\delta>0$, with
$\epsilon^2=\delta/(9\pi)+O(\delta^2)$.

Let $q=w_{\theta\theta}+w=-3\cos2\theta-15\cos4\theta$.
Formula \eqref{tp:eq:equatorial-curvature} gives, in every fixed
spatial $C^j$ norm,
\[
 \kappa=\epsilon q+O(\epsilon^3),\qquad
 \kappa_s=\epsilon q_\theta+O(\epsilon^3),\qquad
 \dd s=(1+O(\epsilon^2))\dd\theta.
\]
Here $\partial_s=v^{-1}\partial_\theta$. The differentiated
remainders follow by Taylor expansion of the curvature formula
and this arclength conversion: $c$ and $u$ are odd
analytic families, $v$ is even, and its denominator is bounded away
from zero on the compact parameter circle.
For $q=a\cos2\theta+b\cos4\theta$, trigonometric orthogonality gives
\begin{align*}
 \int q^3\dd\theta &=\tfrac{3\pi}{2}a^2b=-\tfrac{405\pi}{2},&
 \int q q_\theta^2\dd\theta &=6\pi a^2b=-810\pi,\\
 \int q^2\dd\theta &=\pi(a^2+b^2)=234\pi,&
 \int q_\theta^2\dd\theta &=\pi(4a^2+16b^2)=3636\pi .
\end{align*}
In the cubic terms only products with two factors of frequency $2$
and one of frequency $4$ contribute a constant term. The numerical
values use $a=-3$ and $b=-15$.
Lemma~\ref{tp:lem:graph-velocities} now gives
\begin{align*}
 \mathcal A'_2&=\epsilon^3\int(-q^3+8q q_\theta^2)\dd\theta
                                      +O(\epsilon^5)
       =-\frac{12555\pi}{2}\epsilon^3+O(\epsilon^5),\\
 L'_2&=\epsilon^2\int(q^2-2q_\theta^2)\dd\theta
                                      +O(\epsilon^4)
       =-7038\pi\epsilon^2+O(\epsilon^4).
\end{align*}
Substituting the amplitude from \eqref{tp:eq:mixed-length} proves
\eqref{tp:eq:mixed-velocity}.

Both graph families are smooth embeddings for small amplitude,
and their leading curvature terms are nonconstant. Proposition
\ref{tp:prop:Hopf-rank} gives rank two. Their Hopf Willmore energies tend to
$2\pi^2<8\pi$, so after decreasing a single threshold $\delta_0$
both lie strictly below $8\pi$ for every $0<\delta<\delta_0$.
All remainder constants concern these fixed analytic families and
are uniform on a sufficiently small amplitude interval.
The profiles converge in every spatial $C^j$ norm to $C$.
The pullback Hopf bundles can be smoothly trivialized in the amplitude
parameter, as in Lemma~\ref{tp:lem:Hopf-reconstruction}, which gives
the stated convergence on a fixed torus.

The first velocity has zero area component and nonzero length
component; the second has nonzero area component. More quantitatively,
\[
 \det\begin{pmatrix}\mathcal A'_1&\mathcal A'_2\\ L'_1&L'_2\end{pmatrix}
 =-\frac{9765}{\sqrt\pi}\delta^{5/2}+O(\delta^{7/2})\ne0
\]
for every sufficiently small positive $\delta$. In the ordered
coordinates $(\Re\tau,\Im\tau)$ the determinant is divided by
$(4\pi)^2$ and hence is nonzero.
\end{proof}

\begin{corollary}[Obstruction under monotone recalibration of time]
\label{tp:cor:no-time-change}
For the two immersions $F_1,F_2$ in
Theorem~\ref{tp:thm:noncollinear}, with their common standard
marked point, there is
no vector field $b$ near their common marked point and no assignment
of positive numbers $\rho(F)$ such that
\[
 \rho(F)DP_F\mathcal X(F)=b(P(F))
\]
holds for both immersions. In particular, no separate $C^1$
recalibrations of time with strictly positive finite derivative make all such
projected trajectories solve one autonomous equation.
\end{corollary}
\begin{proof}
The right-hand side would be one vector at the common marked point.
The left-hand sides lie on the two distinct nonzero velocity lines
from Theorem~\ref{tp:thm:noncollinear}; multiplying by positive numbers
cannot make them equal. Let $F_i(t)$ denote the normal MIWF with
$F_i(0)=F_i$, for $i=1,2$. For a change of time $t=\chi_i(s)$ with $\chi_i(s_i)=0$ and
$0<\chi_i'(s_i)<\infty$, the chain rule gives
\[
 \left.\frac{\mathrm d}{\mathrm ds}P(F_i(\chi_i(s)))\right|_{s=s_i}
 =\chi_i'(s_i)\,DP_{F_i}\mathcal X(F_i).
\]
Thus every admissible change of time gives exactly the positive
multiplier excluded above. If the same marked finite cover is taken
for both tori, its projected velocities are obtained by applying the
same invertible $D\mathfrak c_\tau$ to the two vectors. Their linear
independence, and hence the time-change obstruction, persists on the
cover. The energy on each cover is its degree times the base energy.
\end{proof}

\clearpage
\section{Local coupled coordinates and endpoint regularity}\label{tp:sec:coupled}

\subsection{Local coordinates supplied by full rank}
\begin{proposition}[A local submersion chart]\label{tp:prop:chart}
Let $X$ be a Banach space, $U\subset X$ open, $P:U\to\R^2$ a
$C^\ell$ map, where $\ell\geq1$ is an integer, and $x_0\in U$ with
$DP(x_0)$ onto.
Choose a bounded linear map $R:\R^2\to X$ with $DP(x_0)R=I$, put
$K=\ker DP(x_0)$, and $\Pi_K=I-RDP(x_0)$. Then
\[
 x\longmapsto(P(x),\Pi_K(x-x_0))
\]
is a $C^\ell$ coordinate map near $x_0$, with values near
$(P(x_0),0)$ in $\R^2\times K$. Write its coordinates as
$\tau=P(x)$ and $k=\Pi_K(x-x_0)$. Along a smooth trajectory
$\dot x=\mathcal X(x)$ within this chart, the exact transformed
equations are
\begin{equation}\label{tp:eq:coupled}
 \dot\tau=DP(x)\mathcal X(x),\qquad
 \dot k=\Pi_K\mathcal X(x),\qquad x=\Psi(\tau,k).
\end{equation}
\end{proposition}
\begin{proof}
An onto map to $\R^2$ has a bounded right inverse: choose two
preimages of its standard basis. The displayed derivative at $x_0$
is the isomorphism $h\mapsto(DP(x_0)h,\Pi_Kh)$, with inverse
$(z,w)\mapsto Rz+w$. The Banach inverse function theorem gives
the chart. The chain rule gives \eqref{tp:eq:coupled} along every
smooth trajectory contained in the coordinate neighborhood.
\end{proof}

In the Hopf sector this statement has a concrete Banach realization.
Near a fixed smooth regular profile $\eta$, use normal graphs
$\Gamma_u=\cos u\,\eta+\sin u\,\nu_\eta$, with
$u\in H^\ell$ for an integer $\ell\geq4$, on its fixed arclength circle. Local inverse
branches over that parameter circle give a graph gauge for nearby
parametrized curves, also for immersed and multiply covered profiles.
For a multiple cover this is a branchwise statement; deck symmetries
can give different parametrizations of the same curve. The maps $u\mapsto L(\Gamma_u)$
and a local real lift $u\mapsto \mathcal A(\Gamma_u)$ are smooth.
For length this follows from integrating the smooth positive speed.
For area, use the oriented short-geodesic cylinder between $\eta$
and its graph to define the difference from $\mathcal A(\eta)$; its coefficients
are smooth functions of $u,u_s$. Their first variations at zero are
\begin{equation}\label{tp:eq:AL-differential}
 D\mathcal A(0)[h]=-\int h\dd s,\qquad
 DL(0)[h]=-\int\kappa_\eta h\dd s.
\end{equation}
Identify $\mathbb C$ with $\R^2$ and take
$P(u)=(\mathcal A(\Gamma_u)+iL(\Gamma_u))/(4\pi)$,
so that the coordinate $\tau=P(u)$ is the marked modulus throughout.
Assume for this chart and its evolution equations that $\kappa_\eta$
is nonconstant. Then $DP(0)$ is onto by \eqref{tp:eq:Gram}.
A right inverse of $DP(0)$ can be chosen with
smooth image spanned by $1$ and $\kappa_\eta$. More explicitly, if
$J(x,y)=-(y+ix)/(4\pi)$, then $DP(0)=JR_\eta$, so
$R:=R_\eta^\dagger J^{-1}$ is a bounded right inverse of $DP(0)$ from $\R^2$ into each
$H^\ell$ used here. The inverse coordinate map $\Psi(\tau,0)$ supplies
a smooth local section of $P$.

\subsection{Dependence on the fibre variable}
Write \eqref{tp:eq:coupled} as
$\dot\tau=b(\tau,k)$, $\partial_t k=c(\tau,k)$.
The profile variable $k$ belongs to the infinite-dimensional space $K$.
Thus the second equation retains a spatial evolution equation,
while the first is an ordinary differential equation coupled to that
profile. The obstruction in Theorem~\ref{tp:thm:nonprojectable}
concerns a common equation $\dot\tau=b(\tau)$ for the projected
trajectories alone; it is compatible with this coupled description.

The fourth-order graph equation loses derivatives, so the following
interpretation specifies the domain and target of the coupled field.
Fix an integer $\ell\geq4$ and the right inverse $R$ constructed above.
On a sufficiently small $H^\ell$ graph neighborhood, the graph velocity
$\mathcal X$ given by \eqref{miwf:graphflow} is a smooth map
from $H^\ell$ to $H^{\ell-4}$. To verify these assertions, $H^{\ell-2}$
is a multiplication algebra in dimension one, and the graph speed
is bounded away from zero. The finite-jet coefficients
$a_4(u),q_0(u),q_1(u)$ are therefore smooth functions of $u$ with
values in $H^{\ell-2}$. Multiplication of $u_{ssss}\in H^{\ell-4}$
by $a_4(u)\in H^{\ell-2}$ is continuous in $H^{\ell-4}$; the
term $q_1(u)u_{sss}$ has at least this regularity as well.
These multiplication statements also hold at $\ell=4$, where the
target is $L^2$ and $H^2\hookrightarrow L^\infty$.

There is an explicit extension of the differential. In the fixed
arclength coordinate of $\eta$, put
$a(u)=\cos u-\kappa_\eta\sin u$. A graph variation $h$ has
normal component $a(u)h/v_u$, and $\dd s_{\Gamma_u}=v_u\dd s$.
The area and length first variations consequently give
\begin{equation}\label{tp:eq:graph-modulus-differential}
 DP(u)[h]=-\frac1{4\pi}\int_{\mathcal C_\eta}
             a(u)\bigl(1+i\kappa[u]\bigr)h\,\dd s .
\end{equation}
The coefficient belongs to $H^{\ell-2}$ and depends smoothly on
$u\in H^\ell$. Cauchy--Schwarz therefore makes this a bounded
functional even on $L^2$, with smooth dependence in operator norm.
By density the formula agrees with $DP(u)$ on $H^\ell$ and extends
to $H^{\ell-4}$. At $u=0$ it gives the bounded extension of
$\Pi_K=I-RDP(0)$ to $H^{\ell-4}$. Hence, in these coordinates, $b$ is $\R^2$-valued and $c$ is a smooth
map from an open subset of $\R^2\times(K\cap H^\ell)$ into
$\Pi_K H^{\ell-4}$. The equations hold along smooth trajectories by
the chain rule; explicitly,
$b(\tau,k)=DP(\Psi(\tau,k))\mathcal X(\Psi(\tau,k))$ and
$c(\tau,k)=\Pi_K\mathcal X(\Psi(\tau,k))$.
In particular, writing $u=\Psi(\tau,k)$ and
$v_u=|\partial_s\Gamma_u|$ for its graph speed, substitution of the graph normal velocity into
\eqref{tp:eq:graph-modulus-differential} gives
\begin{equation}\label{tp:eq:coupled-base-velocity}
 b(\tau,k)=\frac1{4\pi}\int_{\mathcal C_\eta}
       \frac{v_u\bigl(1+i\kappa[u]\bigr)G(\Gamma_u)}
            {(1+\kappa[u]^2)^2}\,\dd s.
\end{equation}
Since $\dd s_{\Gamma_u}=v_u\dd s$, its real and imaginary
components agree with \eqref{tp:eq:projected-velocity}, each divided
by $4\pi$.
The space $\Pi_KH^{\ell-4}$ is the closed kernel of the extended
$DP(0):H^{\ell-4}\to\R^2$, because $DP(0)R=I$ still holds.
The local section $k=0$ is tangent to the vector field precisely
where $c(\tau,0)=0$.

\subsection{Finite endpoints and continuation}\label{tp:sec:discussion}
In the low-energy genus-one case, \cite[Theorem~1.1(3)]{Jakob}
obtains strong $W^{2,2}$ convergence along a reparametrized
subsequence under an endpoint energy identity. For a finite-time
reduced trajectory, Corollary~\ref{miwf:endpoint-W42} proves the
profile energy identity and full strong $W^{2,2}$ convergence if
dissipation stays bounded along an approaching sequence. Smooth
continuation also requires a compatible Cauchy trace, as discussed
in Section~\ref{sec:topology}.

Weak convergence alone does not justify passing nonlinear curvature
products to products of weak limits: for example, on $\R/(2\pi\miwfZ)$,
$\sin(jx)$ converges weakly in $L^2$
to zero while its square converges weakly to $1/2$.
A divergence-form formulation of the Willmore equation, as in
\cite{RiviereAnalysis}, is interpreted as an identity of distributions:
one tests against smooth vector fields and integrates by parts.
Suppose smooth immersions $F_j:\Sigma\to\Sph^3$ converge weakly in
$W^{2,2}$ to an immersion $F$. To derive an equation for $F$ from
those for $F_j$, one must justify passage to the limit in the nonlinear
terms of this tested identity. This requires suitable compactness or
compensation estimates; writing the equation in divergence form does
not replace that argument. If the equations have a vector-valued
forcing term, its limit must also be identified in the same
parametrization. Bounds on the absolute normal speed alone do not
determine this vector-valued limit.

\subsection{Equatorial multiplicity and energy}
For an $m$-fold equatorial profile with additional torus covering
degree $d$, the Hopf identities give
\[
 \miwfE=2\pi m,\qquad \Will=2\pi^2dm.
\]
Thus the surface bound $\Will\leq8\pi$ corresponds to
$d\miwfE\leq8$, whereas a fourfold equator in the standard
pullback model has energy $8\pi^2$. By endpoint lower
semicontinuity, an equatorial endpoint under that surface bound is
consequently once traversed with $d=m=1$.
At arbitrary energy, the circular trajectories in
Proposition~\ref{tp:prop:circles} retain every positive integer
covering multiplicity.

For example, consider the perturbed fourfold equator
\[
 Z_\epsilon(\theta)=\cos(\epsilon\sin3\theta)
       (\cos4\theta,\sin4\theta,0)+\sin(\epsilon\sin3\theta)e_3.
\]
In reference arclength $s=4\theta$, its linearized geodesic curvature is
$(7/16)\epsilon\sin(3s/4)$. Hence the curvature is nonconstant for all
sufficiently small nonzero $\epsilon$. Theorem~\ref{tp:thm:no-circle}
therefore gives full rank at each finite smooth time; it does not determine
the maximal lifetime or exclude an infinite-time double-equator limit.

For this example take the standard pullback Hopf immersion
$X_{Z_\epsilon}$, so that $d=1$ in \eqref{eq:hopf-cover-energy}.
Since $Z_\epsilon$ converges smoothly to the fourfold equator,
\[
 \Will(X_{Z_\epsilon})=\pi\miwfE(Z_\epsilon)
       \longrightarrow 8\pi^2>8\pi
       \qquad\text{as }\epsilon\to0.
\]
Thus these perturbations lie outside the surface-energy range
$\Will\leq8\pi$ of the finite-time embedded-endpoint results.

A spherical trajectory may reach a fixed stereographic pole while remaining
smooth in the sphere, so a planar rotation-number argument requires control
of the projection and the limiting curve.

For regular homotopies on the $2$-sphere, the intrinsic tangent lift
takes values in $T^1\Sph^2\simeq SO(3)$ and records its homotopy
class in $\pi_1(SO(3))\simeq\miwfZ/2\miwfZ$; see
\cite[p.~9]{JakobHopfCorrection}. An $m$-fold equator has class $m$ modulo
$2$: its tangent frame makes $m$ full turns, and the lift to the double
cover $\Sph^3\to SO(3)$ has endpoint $(-1)^m$ times its initial point.
This parity identifies the fourfold and double equators with each other,
so distinguishing them requires additional information beyond regular
homotopy on the $2$-sphere.

\clearpage
\section*{Acknowledgements}
The authors used ChatGPT (OpenAI) to assist with exposition, consistency
checks, symbolic and asymptotic calculations, and the formulation and
preliminary checking of proof arguments, including the reduced-flow
estimates, fixed-gauge convergence argument, and projection examples.
The authors take full responsibility for the mathematical contents,
source verification, and final validation of the manuscript.


\begin{thebibliography}{99}

\bibitem{AhmedouJakobUnique}
M.~Ahmedou and R.~Jakob,
\emph{Uniqueness of Finite-Time Varifold Limits for the
M\"obius-Invariant Willmore Flow}, arXiv:2609.00501 [math.AP], 2026.
\href{https://arxiv.org/abs/2609.00501}{arXiv:2609.00501}.

\bibitem{Ahmedou.JakobHopf}
M.~Ahmedou, R.~Jakob, 
\emph{Hopf Reduction and Multiplicity-Resolved Endpoints
for the Classical Willmore Flow of Hopf Tori}, arXiv:2609.26997 [math.AP], 2026.

\bibitem{Allard}
W.~K.~Allard,
\emph{On the first variation of a varifold},
Ann. of Math. (2) \textbf{95} (1972), 417--491.

\bibitem{Amann}
H.~Amann,
\emph{Cauchy problems for parabolic equations in Sobolev--Slobodeckii
and H\"older spaces on uniformly regular Riemannian manifolds},
J. Evol. Equ. \textbf{17} (2017), 51--100.
\href{https://arxiv.org/abs/1602.02232}{arXiv:1602.02232}.

\bibitem{BernardRiviere}
Y.~Bernard and T.~Rivi\`ere,
\emph{Energy quantization for Willmore surfaces and applications},
Ann. of Math. (2) \textbf{180} (2014), 87--136.

\bibitem{Blatt}
S.~Blatt,
\emph{A singular example for the Willmore flow},
Analysis (Munich) \textbf{29} (2009), no.~4, 407--430.

\bibitem{ChillFasangovaSchaetzle}
R.~Chill, E.~Fa\v{s}angov\'a, and R.~Sch\"atzle,
\emph{Willmore blowups are never compact},
Duke Math. J. \textbf{147} (2009), no.~2, 345--376.
\href{https://doi.org/10.1215/00127094-2009-014}{doi:10.1215/00127094-2009-014}.

\bibitem{DallAcquaSphere}
A.~Dall'Acqua, T.~Laux, C.-C.~Lin, P.~Pozzi,
and A.~Spener, \emph{The elastic flow of curves on the sphere},
Geom. Flows \textbf{3} (2018), 1--13.
\href{https://doi.org/10.1515/geofl-2018-0001}{doi:10.1515/geofl-2018-0001}.

\bibitem{DallAcquaMuellerSchaetzleSpener}
A.~Dall'Acqua, M.~M\"uller, R.~Sch\"atzle, and A.~Spener,
\emph{The Willmore flow of tori of revolution},
Anal. PDE \textbf{17} (2024), no.~9, 3079--3124.

\bibitem{DallAcquaPozziSpener}
A.~Dall'Acqua, P.~Pozzi, and A.~Spener,
\emph{The \L ojasiewicz--Simon gradient inequality for open elastic
curves}, J. Differential Equations \textbf{261} (2016), 2168--2209.
\href{https://doi.org/10.1016/j.jde.2016.04.027}{doi:10.1016/j.jde.2016.04.027}.

\bibitem{Hatcher}
A.~Hatcher,
\emph{Algebraic Topology},
Cambridge University Press, Cambridge, 2002.

\bibitem{JakobShortTime}
R.~Jakob,
\emph{Short-time existence of the M\"obius-invariant Willmore flow},
J. Geom. Anal. \textbf{28} (2018), no.~2, 1151--1181.
Published online 2 May 2017.
\href{https://doi.org/10.1007/s12220-017-9857-5}{doi:10.1007/s12220-017-9857-5}.

\bibitem{JakobRegularity}
R.~Jakob,
\emph{Functional analytic properties and regularity of the
M\"obius-invariant Willmore flow in $\R^n$},
Partial Differ. Equ. Appl. \textbf{3} (2022), Paper No.~67.
\href{https://doi.org/10.1007/s42985-022-00205-z}{doi:10.1007/s42985-022-00205-z}.
\href{https://arxiv.org/abs/2011.03832v10}{arXiv:2011.03832v10}.

\bibitem{JakobHopf}
R.~Jakob,
\emph{The Willmore flow of Hopf-tori in the $3$-sphere},
J. Evol. Equ. \textbf{23} (2023), Paper No.~72.

\bibitem{JakobGlobal}
R.~Jakob,
\emph{Global existence and full convergence of the M\"obius-invariant
Willmore flow in the $3$-sphere},
J. Geom. Anal. \textbf{34} (2024), Paper No.~24.

\bibitem{Jakob}
R.~Jakob,
\emph{Singularities and full convergence of the M\"obius-invariant
Willmore flow in the $3$-sphere},
Asian J. Math. \textbf{29} (2025), no.~5, 635--702.

\bibitem{JakobHopfCorrection}
R.~Jakob,
\emph{Corrections to: ``The Willmore flow of Hopf-tori in the
$3$-sphere''}, arXiv:2602.13228v1, 2026.

\bibitem{KuwertLi}
E.~Kuwert and Y.~Li,
\emph{$W^{2,2}$-conformal immersions of a closed Riemann surface into
$\R^n$},
Comm. Anal. Geom. \textbf{20} (2012), no.~2, 313--340.

\bibitem{KuwertSchaetzleSmall}
E.~Kuwert and R.~Sch\"atzle,
\emph{The Willmore flow with small initial energy},
J. Differential Geom. \textbf{57} (2001), no.~3, 409--441.

\bibitem{KuwertSchaetzleGradient}
E.~Kuwert and R.~Sch\"atzle,
\emph{Gradient flow for the Willmore functional},
Comm. Anal. Geom. \textbf{10} (2002), no.~2, 307--339.

\bibitem{KuwertSchaetzleBounds}
E.~Kuwert and R.~Sch\"atzle,
\emph{Closed surfaces with bounds on their Willmore energy},
Ann. Sc. Norm. Super. Pisa Cl. Sci. (5) \textbf{11} (2012), 605--634.

\bibitem{KuwertSchaetzleConformal}
E.~Kuwert and R.~Sch\"atzle,
\emph{Minimizers of the Willmore functional under fixed conformal class},
J. Differential Geom. \textbf{93} (2013), 471--530.
\href{https://arxiv.org/abs/1009.6168}{arXiv:1009.6168}.

\bibitem{KuwertScheuer}
E.~Kuwert and J.~Scheuer,
\emph{Asymptotic estimates for the Willmore flow with small energy},
Int. Math. Res. Not. IMRN (2021), no.~18, 14252--14266.

\bibitem{LanMartinoRiviere}
T.~Lan, D.~Martino, and T.~Rivi\`ere,
\emph{The analysis of Willmore surfaces and its generalizations in
higher dimensions},
J. Math. Study \textbf{59} (2026), no.~1, 80--188.
\href{https://doi.org/10.4208/jms.v59n1.26.05}{doi:10.4208/jms.v59n1.26.05}.
\href{https://arxiv.org/abs/2511.01777}{arXiv:2511.01777}.

\bibitem{LiYau}
P.~Li and S.-T.~Yau,
\emph{A new conformal invariant and its applications to the Willmore
conjecture and the first eigenvalue of compact surfaces},
Invent. Math. \textbf{69} (1982), no.~2, 269--291.

\bibitem{MantegazzaPozzetta}
C.~Mantegazza and M.~Pozzetta,
\emph{The \L ojasiewicz--Simon inequality for the elastic flow},
Calc. Var. Partial Differential Equations \textbf{60} (2021), no.~1,
Paper No.~56, 17~pp.
\href{https://doi.org/10.1007/s00526-020-01916-0}{doi:10.1007/s00526-020-01916-0}.

\bibitem{MarquesNeves}
F.~C.~Marques and A.~Neves,
\emph{Min-max theory and the Willmore conjecture},
Ann. of Math. (2) \textbf{179} (2014), no.~2, 683--782.

\bibitem{MayerSimonett}
U.~F.~Mayer and G.~Simonett,
\emph{Self-intersections for Willmore flow},
in: \emph{Evolution equations: applications to physics, industry,
life sciences and economics}, Progr. Nonlinear Differential Equations
Appl., vol.~55, Birkh\"auser, Basel, 2003, pp.~341--348.

\bibitem{PalmurellaRiviereApproach}
F.~Palmurella and T.~Rivi\`ere,
\emph{The parametric approach to the Willmore flow},
Adv. Math. \textbf{400} (2022), Paper No.~108257, 48 pp.

\bibitem{PalmurellaRiviereFlow}
F.~Palmurella and T.~Rivi\`ere,
\emph{The parametric Willmore flow},
J. Reine Angew. Math. \textbf{811} (2024), 1--91.

\bibitem{Pinkall}
U.~Pinkall,
\emph{Hopf tori in $S^3$},
Invent. Math. \textbf{81} (1985), 379--386.

\bibitem{RiviereAnalysis}
T.~Rivi\`ere,
\emph{Analysis aspects of Willmore surfaces},
Invent. Math. \textbf{174} (2008), no.~1, 1--45.

\bibitem{RiviereDegenerating}
T.~Rivi\`ere,
\emph{Lipschitz conformal immersions from degenerating Riemann
surfaces with $L^2$-bounded second fundamental forms},
Adv. Calc. Var. \textbf{6} (2013), no.~1, 1--31.

\bibitem{RuppSpener}
F.~Rupp and A.~Spener,
\emph{Existence and convergence of the length-preserving elastic flow
of clamped curves},
J. Evol. Equ. \textbf{24} (2024), no.~3, Paper No.~59, 41~pp.

\bibitem{SchaetzleLowerSemicontinuity}
R.~Sch\"atzle,
\emph{Lower semicontinuity of the Willmore functional for currents},
J. Differential Geom. \textbf{81} (2009), no.~2, 437--456.

\bibitem{SchaetzleConformalFactor}
R.~Sch\"atzle,
\emph{Estimation of the conformal factor under bounded Willmore
energy},
Math. Z. \textbf{274} (2013), no.~3--4, 1341--1383.

\bibitem{SimonGMT}
L.~Simon,
\emph{Introduction to Geometric Measure Theory},
revised lecture notes, 2018.
\url{https://math.stanford.edu/~lms/ntu-gmt-text.pdf}.

\bibitem{Simonett}
G.~Simonett,
\emph{The Willmore flow near spheres},
Differential Integral Equations \textbf{14} (2001), no.~8, 1005--1014.

\bibitem{Skorzinski}
F.~Skorzinski,
\emph{Local minimizers of the Willmore functional},
Analysis \textbf{35} (2015), no.~2, 93--115.
\href{https://doi.org/10.1515/anly-2012-1274}{doi:10.1515/anly-2012-1274}.

\bibitem{Weiner}
J.~L.~Weiner,
\emph{On a problem of Chen, Willmore, et al.},
Indiana Univ. Math. J. \textbf{27} (1978), no.~1, 19--35.
\end{thebibliography}
\end{document}